\documentclass[pdflatex,sn-mathphys-num]{sn-jnl}

\usepackage{amsmath,amssymb,mathtools}
\numberwithin{equation}{section}
\allowdisplaybreaks[2]

\usepackage{enumitem}
\usepackage{booktabs,array}
\usepackage{xcolor}
\usepackage{tikz}
\makeatletter
\renewcommand\paragraph{\@startsection{paragraph}{4}{\z@}%
  {6pt \@plus 2pt \@minus 1pt}{-0.5em}{\paragraphfont}}
\makeatother

\theoremstyle{thmstyleone}
\newtheorem{mainthm}{Theorem} %
\newtheorem{theorem}{Theorem}[section] %
\newtheorem{lemma}[theorem]{Lemma}
\newtheorem{proposition}[theorem]{Proposition}
\newtheorem{corollary}[theorem]{Corollary}
\theoremstyle{thmstylethree}
\newtheorem{definition}[theorem]{Definition}
\theoremstyle{thmstyletwo}
\newtheorem{remark}[theorem]{Remark}

\newcommand{\RR}{\mathbb{R}}
\newcommand{\CC}{\mathbb{C}}
\newcommand{\ZZ}{\mathbb{Z}}

\newcommand{\TT}{\mathbb{T}}
\renewcommand{\Re}{\operatorname{Re}}
\renewcommand{\Im}{\operatorname{Im}}

\DeclareMathOperator{\diag}{diag}
\newcommand{\eps}{\varepsilon}
\newcommand{\epsp}{\varepsilon'}
\newcommand{\bom}{\bar\omega}
\newcommand{\cD}{\mathcal{D}} %
\newcommand{\sD}{\mathsf{D}} %
\newcommand{\cL}{\mathcal{L}} %
\newcommand{\cG}{\mathcal{G}}
\newcommand{\cU}{\mathcal{U}}
\newcommand{\cJ}{\mathcal{J}}
\newcommand{\hE}{\hat{E}} %
\newcommand{\hY}{\hat{Y}}
\newcommand{\hQ}{\hat{Q}}
\newcommand{\fc}[2]{#1^{[#2]}} %
\urldef{\urlGallavottiProblem}\url{https://web.math.princeton.edu/~aizenman/OpenProblems_MathPhys/9805.Separatrix.html}
\begin{document}

\title[Separatrix splitting in Gallavotti's problem]{Asymptotics of separatrix
splitting\texorpdfstring{\linebreak}{ } in Gallavotti's problem}

\author*[1]{\fnm{Aleksandr} \sur{Gulkanov}}\email{Aleksandr.Gulkanov@rice.edu}
\affil*[1]{\orgdiv{Applied Physics}, \orgname{Rice University},
  \orgaddress{\street{6100 Main St.}, \city{Houston}, \state{Texas} \postcode{77005},
  \country{USA}}}

\abstract{\unboldmath
We consider a pendulum under quasi-periodic forcing with frequencies of orders $\sqrt\eta$ and
$1/\sqrt\eta$. The splitting of the separatrices of its invariant torus is measured by the
Jacobian determinant $\Delta$ of the difference of the actions on the unstable and stable
manifolds at a homoclinic point. Gallavotti (1998) asked for the asymptotics of $\Delta$ as
$\eta\to0$ when the forcing contains a term of order one depending only on the fast angle, while
the rest has amplitude $\eta^c$ with a large fixed $c$. We prove that for $c$ above an absolute
constant $\Delta=[C+o(1)]\eta^{2c-1/2}e^{-b/\sqrt\eta}$ with the exponent $b$ of first-order
perturbation theory, uniformly in real amplitudes of the fast term on compact sets. No
Diophantine condition is needed. The coefficient $C$ is an explicit
constant times an entire function of the amplitude of the fast term, namely the response of the
first Stokes constant of the inner equation with this term as background to adding a profile
determined by the rest of the forcing. We express it as a convergent series of values of
explicit Borel transforms, whose terms are nonnegative when the Fourier coefficients of the
background and of the added profile are. For Gallavotti's example this function is at least
$2\pi$, so the coefficient is positive for every real amplitude. The same conclusion holds for
the other reading of the example, where for negative amplitudes the proof uses finitely many
polynomial inequalities verified in exact arithmetic. We use results of Baldom\'a, Fontich,
Guardia and Seara on the periodically forced pendulum.}


\maketitle

\section{Introduction}\label{sec:intro}

We study the splitting of the separatrices of a pendulum under a quasi-periodic forcing with one
slow and one fast frequency. Gallavotti~\cite{Gallavotti98} posed the problem of finding the
asymptotics of this splitting when the part of the forcing that depends only on the fast angle,
called the strong term below, is not small. If the whole forcing is small, the splitting is
exponentially small in the fast frequency, and for a generic forcing its leading term is the one
given by first-order perturbation theory~\cite{GGM99}. For a pendulum with a fast periodic
forcing the asymptotics is known also when the forcing is of order one, and the leading
coefficient is then the Stokes constant of an inner equation rather than the Melnikov
integral~\cite{BFGS12}. In this paper we prove an asymptotic formula for the splitting in the
problem of~\cite{Gallavotti98}, assuming that the amplitude of the remaining part of the forcing
is a sufficiently high power of the small parameter (Theorem~\ref{thm:main}). The strong term
changes neither the exponent nor the power of the small parameter. It changes the leading
coefficient, which becomes an entire function of the amplitude of the strong term, and we
express this function as a convergent series. For the example proposed
in~\cite{Gallavotti98} this function has no real zeros (Theorems~\ref{thm:example}
and~\ref{thm:cos2}).

\subsection{Gallavotti's problem}\label{ss:problem}

We consider the Hamiltonian system with three degrees of freedom and Hamiltonian
\begin{equation}\label{eq:H}
H=\omega\cdot A+\frac{I^2}{2J}+\bigl[Jg^2+\eps f(\alpha)+\epsp\,m(2\alpha_2)\bigr](\cos\varphi-1),
\qquad \omega=\bigl(\sqrt\eta\,\bom_1,\ \bom_2/\sqrt\eta\bigr).
\end{equation}
Here $(A,\alpha)\in\RR^2\times\TT^2$ and $(I,\varphi)\in\RR\times\TT$ are pairs of canonically
conjugate variables, with $\TT=\RR/2\pi\ZZ$. The constants $J,g>0$ and $\bom_1,\bom_2\ne0$ are
fixed, while $\eta\in(0,1]$, $\eps$ and $\epsp$ are real parameters. The functions
$f(\alpha_1,\alpha_2)$ and $m(\theta)$ are real trigonometric polynomials with zero mean. We
assume throughout that they are even, that is, $f(-\alpha)=f(\alpha)$ and
$m(-\theta)=m(\theta)$.

The variables $(I,\varphi)$ describe a pendulum with moment of inertia $J$. For
$\eps=\epsp=0$ the frequency of its small oscillations is $g/2\pi$, and the equilibrium
$\varphi=0$ is unstable. The angles rotate uniformly, $\alpha(t)=\alpha+\omega t$, and the
equation of motion of the pendulum reads
\begin{equation}\label{eq:pendulum}
J\ddot\varphi=\bigl[Jg^2+\eps f(\alpha+\omega t)+\epsp\,m(2\alpha_2+2\omega_2t)\bigr]\sin\varphi .
\end{equation}
This is a pendulum with a quasi-periodic parametric forcing that has one slow and one fast
frequency. As $\eta\to0$ the system has three frequency scales, namely $\sqrt\eta$ (the angle
$\alpha_1$), $1$ (the pendulum) and $1/\sqrt\eta$ (the angle $\alpha_2$). We call the term with
$\eps$ the weak perturbation and the term with $\epsp$ the strong term, or the background. The
strong term depends only on the fast angle and contains only even harmonics in $\alpha_2$.
According to~\cite{Gallavotti98}, Hamiltonians of exactly this kind arise in problems of
celestial mechanics.

For all $\eps$ and $\epsp$ the system has the solutions $I=0$, $\varphi=0$, $A=\mathrm{const}$,
$\alpha(t)=\alpha+\omega t$, which fill two-dimensional invariant tori. Since the Hamiltonian is
linear in $A$, it suffices to consider the torus $\mathcal T=\{I=0,\ \varphi=0,\ A=0\}$. For
$\eps=\epsp=0$ its unstable and stable manifolds coincide and form the product of $\TT^2$ and
the separatrix of the pendulum, which lies in $\{A=0\}$. Under the perturbation the manifolds
split.

We measure the splitting on the section $\varphi=\pi$. We consider the branches of the
manifolds that follow the separatrix with $I>0$, and the other separatrix is obtained by the
symmetry $(\varphi,I)\mapsto(-\varphi,-I)$. The trajectories of the unstable manifold $W^u$ of
$\mathcal T$ that stay close to the separatrix are parametrized by their angles at the moment of
the first arrival at the section. More precisely, for every $\alpha\in\TT^2$ there is exactly
one such trajectory $\varphi(t)$, $t\le0$, with $\varphi(0)=\pi$ and angles $\alpha$ at $t=0$.
We denote by $A^u(\alpha)$ the value of the actions at this point of the section. Since
$\dot A_i=-\partial H/\partial\alpha_i$ and $A\to0$ as $t\to-\infty$, we have
\begin{equation}\label{eq:Au}
A^u_i(\alpha)=\int_{-\infty}^0\partial_{\alpha_i}\bigl[\eps f(\alpha)+\epsp\,m(2\alpha_2)\bigr]
\Big|_{\alpha+\omega t}\bigl(1-\cos\varphi(t)\bigr)\,dt,\qquad i=1,2 .
\end{equation}
In the same way we define $A^s(\alpha)$ for the stable manifold $W^s$ (an integral over
$t\ge0$, with the minus sign, along the trajectory that tends to the torus as $t\to+\infty$).
Since $f$ and $m$ are even, the system is reversible with respect to the involution
$(\varphi,I,\alpha,A)\mapsto(-\varphi,I,-\alpha,A)$, which interchanges $W^u$ and $W^s$. Hence
$A^s(\alpha)=A^u(-\alpha)$, the vector $A^u(\alpha)-A^s(\alpha)$ vanishes at $\alpha=0$, and
the point of the section with $\alpha=0$ is homoclinic. We define the \emph{splitting} as the
determinant of the Jacobian matrix
\begin{equation}\label{eq:Delta}
\Delta(\eps,\epsp,\eta)=\det\partial_\alpha\bigl(A^u-A^s\bigr)\Big|_{\alpha=0}.
\end{equation}

The following question was posed in~\cite{Gallavotti98}. Let $\eps=\eta^c$ with a large fixed
$c$. Can one find the asymptotics of $\Delta$ as $\eta\to0$, or at least a positive lower bound
for $|\Delta|$ valid for $|\epsp|=O(1)$? If the amplitude $\epsp$ is small as well, the answer
is given by first-order perturbation theory~\cite{GGM99}. For generic $f$ the splitting is then
exponentially small and proportional, up to a power of $\eta$, to $\eps^2e^{-b/\sqrt\eta}$ with
$b=\pi|\bom_2|/(2g)$, and its leading term does not depend on $m$. The difficulty pointed out
in the problem statement is twofold. First, the entries of the Jacobian matrix are much larger
than its determinant, so the result comes out only after cancellations. Second, for $\epsp$ of
order one the perturbation series in $\epsp$ still converges, but its first term no longer gives
the asymptotics.

The problem statement ends with a special case, ``as difficult as the above general one'', in
which $f(\alpha)=\cos\alpha_1+\cos\alpha_2$ and $m(\alpha_2)=\cos2\alpha_2$. Since the
Hamiltonian~\eqref{eq:H} contains $m(2\alpha_2)$, the literal reading gives the strong term
$\epsp\cos(4\alpha_2)(\cos\varphi-1)$. It is possible, however, that the strong term
$\epsp\cos(2\alpha_2)(\cos\varphi-1)$ was meant. We treat both readings
(Theorems~\ref{thm:example} and~\ref{thm:cos2}).

\subsection{Main results}\label{ss:results}

We state four theorems. Theorem~\ref{thm:main} gives the asymptotics of the splitting for
general $f$ and $m$. Theorems~\ref{thm:example} and~\ref{thm:cos2} apply it to the two readings
of the example of~\cite{Gallavotti98}, and Theorem~\ref{thm:series} describes the entire
function that enters the leading coefficient.

\paragraph{Notation} We set
\begin{equation}\label{eq:mu-M}
\mu=\frac{\epsp}{Jg^2},\qquad M(\theta)=m(2\theta),
\end{equation}
so that $M$ is a real even trigonometric polynomial with zero mean and period $\pi$. We denote
Fourier coefficients by a superscript in square brackets, as in
$p(\theta)=\sum_k\fc pk e^{ik\theta}$ and
$f(\alpha)=\sum_{j,\ell}\fc f{j,\ell}e^{i(j\alpha_1+\ell\alpha_2)}$. With the perturbation $f$
we associate the number and the polynomial
\begin{equation}\label{eq:Af-pf}
\begin{gathered}
A_f=\frac{d^2}{d\alpha_1^2}\biggl[\frac1{2\pi}\int_0^{2\pi}f(\alpha_1,\alpha_2)\,d\alpha_2\biggr]_{\alpha_1=0}
=-\sum_jj^2\fc f{j,0},\\
p_f(\theta)=\frac{f(0,\theta)-f(0,\theta+\pi)}2 .
\end{gathered}
\end{equation}
The number $A_f$ is determined by the harmonics of $f$ that depend only on the slow angle. The
function $p_f$ is the sum of the odd harmonics of the function $\theta\mapsto f(0,\theta)$, and
it is a real even trigonometric polynomial.

\begin{mainthm}[general case]\label{thm:main}
Let $f(\alpha_1,\alpha_2)$ and $m(\theta)$ be real even trigonometric polynomials with zero
mean, $J,g>0$, $\bom_1\bom_2\ne0$. There is an absolute constant $c_0$ (one can take $c_0=148$)
such that for every fixed $c>c_0$ the following holds. For $\eps=\eta^c$ and $\eta$ small
enough, the manifolds near the section $\varphi=\pi$ are graphs and the
determinant~\eqref{eq:Delta} is defined. As $\eta\to0$,
\begin{equation}\label{eq:main}
\Delta(\eta^c,\epsp,\eta)=\Bigl[-\frac{8|\bom_2|}{g^3}\,A_f\,G_{M,p_f}(\mu)+o(1)\Bigr]\,
\eta^{2c-\frac12}\exp\Bigl[-\frac{\pi|\bom_2|}{2g\sqrt\eta}\Bigr]
\end{equation}
uniformly in $\epsp$ in any compact subset of $\RR$. Here $\mu$, $M$, $A_f$ and $p_f$ are
defined by~\eqref{eq:mu-M} and~\eqref{eq:Af-pf}, and the function $G_{M,p}$ by
Definition~\ref{def:G} below.
\end{mainthm}

Uniformity means that for every compact set $\mathcal K\subset\RR$ the quantity $o(1)$ tends to
zero as $\eta\to0$ uniformly in $\epsp\in\mathcal K$. This quantity depends on $J$, $g$,
$\bom_1$, $\bom_2$, $f$, $m$, $c$ and $\mathcal K$. The theorem needs no Diophantine condition
on the vector $\omega$, and $\eta$ tends to zero through all values, not along a sequence. Since
$\eps=\eta^c$, the factor multiplying the square bracket in~\eqref{eq:main} equals
$\eps^2\eta^{-1/2}\exp[-\pi|\bom_2|/(2g\sqrt\eta)]$.

\paragraph{The function $G_{M,p}$} We express the coefficient of the asymptotics through a
function constructed from the background and a polynomial $p$ as follows.

\begin{definition}\label{def:G}
Let $M$ be a real even trigonometric polynomial with zero mean and period $\pi$, and let $p$ be
a real even trigonometric polynomial containing only odd harmonics. For $n\in2\ZZ$ and $j\ge1$
define the functions $\hY^{[n]}_{(j)}(s)$ recursively:
\begin{gather}
\hY^{[n]}_{(1)}(s)=\frac{n^2\,\fc Mn}{(n-s)^4}\quad(n\ne0),\qquad \hY^{[0]}_{(1)}=0, \label{eq:Y1}\\
\hY^{[n]}_{(j)}(s)=\frac1{(n-s)^4}\int_0^s(n-t)^2\,\hQ^{[n]}_{(j)}(t)\,dt\qquad(j\ge2), \label{eq:Yj}\\
\hQ^{[n]}_{(j)}(t)=\frac12\sum_{\substack{n_1+n_2=n\\ j_1+j_2=j}}\int_0^t
\bigl[(n_1-u)-(n_2-t+u)\bigr]^2\,\hY^{[n_1]}_{(j_1)}(u)\,\hY^{[n_2]}_{(j_2)}(t-u)\,du . \label{eq:Qj}
\end{gather}
Let $*$ denote the convolution $(a*b)(s)=\int_0^sa(t)\,b(s-t)\,dt$. Set
$\hE^{[0]}_{(0)}(s)=s$, $\hE^{[n]}_{(0)}=0$ for $n\ne0$, and for $j\ge1$
\begin{equation}\label{eq:E}
\hE^{[n]}_{(j)}=\sum_{r=1}^{j}\ \sum_{\substack{n_1+\dots+n_r=n\\ j_1+\dots+j_r=j}}
\frac{s^{r+1}}{r!}*\hY^{[n_1]}_{(j_1)}*\dots*\hY^{[n_r]}_{(j_r)} ,
\end{equation}
where $n_i\in2\ZZ$, $j_i\ge1$. Finally,
\begin{equation}\label{eq:G}
G_{M,p}(\mu)=4\pi\sum_{j\ge0}\mu^j\sum_k\fc pk\,\hE^{[1-k]}_{(j)}(1).
\end{equation}
\end{definition}

The sums in~\eqref{eq:Qj} and~\eqref{eq:E} are finite, because $\hY^{[n]}_{(j)}=0$ if
$|n|>j\deg M$. For $n=0$ there is no singularity at zero in~\eqref{eq:Yj}, since
$\hQ^{[0]}_{(j)}(t)=O(t)$. The inner sum in~\eqref{eq:G} runs over the harmonics of $p$. Since
$k$ is odd there, the index $1-k$ is even. The term with $j=0$ equals $4\pi\fc p1$. By
Theorem~\ref{thm:series} below, the series~\eqref{eq:G} converges for all $\mu$.

These formulas originate as follows (Section~\ref{sec:borel}). Let $x=\xi-Y(\xi,\theta)$,
$Y=O(1/\xi)$, be the formal solution of the equation
\begin{equation}\label{eq:newton}
(\sD x)^2-x\,\sD^2x=1+zM(\theta),\qquad \sD=\partial_\xi-i\partial_\theta ,
\end{equation}
to which the inner equation of the problem with background $zM$ reduces. Then
$\hY^{[n]}_{(j)}$ is the Borel transform of the coefficient of $z^je^{in\theta}$ in the series
$Y$, and $\hE^{[n]}_{(j)}$ is the same transform for the function $x^{-2}$. Here the Borel
transform takes $\xi^{-k-1}$ to $s^k/k!$.

If $A_fG_{M,p_f}(\mu)\ne0$, formula~\eqref{eq:main} gives the asymptotics of $\Delta$ and, in
particular, a positive lower bound for $|\Delta|$. The next proposition collects what we can
say about this coefficient.

\begin{proposition}[nondegeneracy]\label{prop:nondeg}
Under the assumptions of Theorem~\ref{thm:main}:
\begin{enumerate}[label=\textup{(\roman*)}]
\item The function $G_{M,p_f}$ is entire and $G_{M,p_f}(0)=4\pi\fc{p_f}{1}$. Hence, if
$A_f\ne0$ and $\fc{p_f}{1}\ne0$, then $A_fG_{M,p_f}(\mu)\ne0$ for small $|\mu|$ and for all
real $\mu$ outside a set without finite accumulation points.
\item If $\mu\fc Mk\ge0$ and $\fc{p_f}{k}\ge0$ for all $k$, and $\fc{p_f}{1}>0$, then
$G_{M,p_f}(\mu)\ge4\pi\fc{p_f}{1}$.
\item The coefficient may vanish identically. Thus $A_f=0$ if $f$ has no harmonics depending
only on $\alpha_1$, and $G_{M,p_f}\equiv0$ for $f=\cos\alpha_1+\cos3\alpha_2$,
$m(\theta)=\cos4\theta$. In these cases \eqref{eq:main} gives only an upper bound. The leading
term then has a different form, and we do not find it in this paper.
\end{enumerate}
\end{proposition}

The conditions of item~(i) hold for generic $f$. Item~(ii) applies for both signs of $\mu$ if
the coefficients $\fc Mk$ have the appropriate sign, and it leaves no exceptional values of the
amplitude.

\paragraph{The example from the problem statement} The next two theorems correspond to the two readings
of the example.

\begin{mainthm}[strong term $\cos4\alpha_2$]\label{thm:example}
Let $f=\cos\alpha_1+\cos\alpha_2$ and $m(\theta)=\cos2\theta$, that is, the strong term is
$\epsp\cos(4\alpha_2)(\cos\varphi-1)$. Then $A_f=-1$, $p_f(\theta)=\cos\theta$,
$M(\theta)=\cos4\theta$, and for $c>c_0$, as $\eta\to0$,
\begin{equation}\label{eq:example}
\Delta(\eta^c,\epsp,\eta)=\frac{8|\bom_2|}{g^3}\;G(\mu)\;
\eta^{2c-\frac12}\exp\Bigl[-\frac{\pi|\bom_2|}{2g\sqrt\eta}\Bigr]\bigl(1+o(1)\bigr)
\end{equation}
uniformly in $\epsp$ in any compact subset of $\RR$. Here $G=G_{\cos4\theta,\cos\theta}$ is an
even entire function,
\begin{equation}\label{eq:G-example}
G(\mu)=2\pi\Bigl(1+\sum_{j\ge1}T_j\,\mu^{2j}\Bigr),\qquad T_j=\hE^{[0]}_{(2j)}(1)\ge0,\qquad
T_1=0.0053202218704\ldots
\end{equation}
Consequently, $G(\mu)\ge2\pi(1+T_1\mu^2)\ge2\pi$ for all real $\mu$, and for every compact set
$\mathcal K\subset\RR$ there is $\eta_0>0$ such that $\Delta(\eta^c,\epsp,\eta)>0$ for
$0<\eta<\eta_0$, $\epsp\in\mathcal K$.
\end{mainthm}

\begin{mainthm}[strong term $\cos2\alpha_2$]\label{thm:cos2}
Let $f=\cos\alpha_1+\cos\alpha_2$ and let the strong term be
$\epsp\cos(2\alpha_2)(\cos\varphi-1)$, that is, $m(\theta)=\cos\theta$ and
$M(\theta)=\cos2\theta$. Then for $c>c_0$, as $\eta\to0$, formula~\eqref{eq:example} holds
with the function
\begin{equation}\label{eq:GB}
G_B(\mu)=G_{\cos2\theta,\cos\theta}(\mu)
=2\pi\sum_{j\ge0}\mu^j\bigl[\hE^{[0]}_{(j)}(1)+\hE^{[2]}_{(j)}(1)\bigr]
\end{equation}
in place of $G$, uniformly in $\epsp$ in any compact subset of $\RR$. Moreover,
\begin{equation}\label{eq:GB-bounds}
\begin{gathered}
G_B(\mu)\ge2\pi\quad(\mu\ge0),\qquad
G_B(\mu)\ge\varkappa_0\cdot2\pi\,\frac{\sinh\lambda}{\lambda}\quad(\mu<0),\\
\lambda=\frac{|\mu|}{\sqrt8},\qquad\varkappa_0=\frac{6142921}{102400000}>0.0599 .
\end{gathered}
\end{equation}
Consequently, $G_B(\mu)>0$ for all real $\mu$, and for every compact set $\mathcal K\subset\RR$
there is $\eta_0>0$ such that $\Delta(\eta^c,\epsp,\eta)>0$ for $0<\eta<\eta_0$,
$\epsp\in\mathcal K$.
\end{mainthm}

For $\mu\ge0$ the bound~\eqref{eq:GB-bounds} follows from item~(ii) of
Proposition~\ref{prop:nondeg}. For $\mu<0$ the terms of the series~\eqref{eq:GB} have different
signs, and we need a separate estimate (Section~\ref{sec:negative}). Its proof includes a
computer check of a finite set of inequalities for polynomials with rational coefficients,
carried out in exact arithmetic. This is the only place where the proof relies on a computer.

\paragraph{Properties of $G_{M,p}$} The last theorem describes the coefficient of the
asymptotics.

\begin{mainthm}[series for the response of the Stokes constant]\label{thm:series}
Let $M$ be a real even trigonometric polynomial with zero mean and period $\pi$, and let $p$ be
a real even trigonometric polynomial containing only odd harmonics.
\begin{enumerate}[label=\textup{(\arabic*)}]
\item The functions $\hY^{[n]}_{(j)}$ and $\hE^{[n]}_{(j)}$ defined by
\eqref{eq:Y1}--\eqref{eq:E} are analytic in the domain
$\CC\setminus\{s\in\RR:\ |s|\ge2\}$. The series~\eqref{eq:G} converges for all $\mu\in\CC$, so
$G_{M,p}$ is an entire function.
\item For every $k\ge1$
\begin{equation}\label{eq:response}
[\mu^{k-1}]\,G_{M,p}=[t^1]\,[z^k]\,\chi\bigl[z(M+tp)\bigr],
\end{equation}
where $\chi[B]$ is the first Stokes constant of the inner equation with profile $B(\theta)$
(Section~\ref{ss:stokes}), and $[z^k]$, $[t^1]$ denote Taylor coefficients. The function
$z\mapsto\chi[zB]$ is entire, and $[z^k]\chi[z(M+tp)]$ is a polynomial in $t$ of degree at most
$k$.
\item If $\fc Mn\ge0$ for all $n$, then $\hY^{[n]}_{(j)}(s)\ge0$ and $\hE^{[n]}_{(j)}(s)\ge0$
for $0\le s<2$. Consequently, if $\mu\fc Mn\ge0$ and $\fc pk\ge0$ for all $n$ and $k$, then all
terms of the series~\eqref{eq:G} are nonnegative and $G_{M,p}(\mu)\ge4\pi\fc p1$.
\end{enumerate}
\end{mainthm}

Equality~\eqref{eq:response} means that $G_{M,p}$ is the response of the Stokes constant to
adding the profile $p$ to the background $\mu M$. Formally,
$G_{M,p}(\mu)=\partial_r\chi[\mu M+rp]|_{r=0}$. In this form the statement would require joint
analyticity of $\chi$ in the coefficients of the profile, which we neither prove nor use. We
prove the analyticity of the functions $\hY^{[n]}_{(j)}$, $\hE^{[n]}_{(j)}$ and the formula for
the Taylor coefficients of the Stokes constant in Section~\ref{sec:borel} for backgrounds all
of whose frequencies are multiples of an arbitrary integer $N\ge2$. The problem requires only
the case $N=2$.

\begin{remark}[range of the weak amplitude]\label{rem:range}
The proof gives more than Theorem~\ref{thm:main} states. The formula
\[
\Delta=\Bigl[-\frac{8|\bom_2|}{g^3}\,A_f\,G_{M,p_f}(\mu)+o(1)\Bigr]\,
\eps^2\eta^{-\frac12}\exp\Bigl[-\frac{\pi|\bom_2|}{2g\sqrt\eta}\Bigr]
\]
holds uniformly in all $\eps$ such that $e^{-\pi/(4h)}\le|\eps|/(Jg^2)\le h^{3P+2}$, where
$h=g\sqrt\eta/|\bom_2|$ and $P$ is an absolute constant, for which one can take $P=98$
(Proposition~\ref{thm:Qprime}). The condition $c>c_0$, $c_0=(3P+2)/2=148$, ensures that
$|\eps|/(Jg^2)\le h^{3P+2}$ for small $\eta$. The constant $c_0$ is not optimal. Apparently it
cannot be taken smaller than $\tfrac12$. We discuss this question in
Section~\ref{ss:exponent}, and the range $\tfrac12<c\le c_0$ remains open.
\end{remark}

\subsection{First-order perturbation theory and the ideas of the proof}\label{ss:mechanism}

For $\epsp=0$ formula~\eqref{eq:main} takes the form
\[
\Delta=\Bigl[-\frac{32\pi|\bom_2|}{g^3}\,A_f\,\fc{p_f}{1}+o(1)\Bigr]\,
\eps^2\eta^{-\frac12}\exp\Bigl[-\frac{\pi|\bom_2|}{2g\sqrt\eta}\Bigr],
\]
because $G_{M,p}(0)=4\pi\fc p1$. The coefficient here is the one given by first-order
perturbation theory, and it does not depend on $m$. By uniformity, the same formula holds if
$\epsp=\epsp(\eta)\to0$. The power of $\eta$ in it is $\eta^{-1/2}$, as in the example computed
in~\cite{GGM99}. The text of~\cite{Gallavotti98} has $\eta^{-1}$, which is apparently a misprint
(Section~\ref{ss:first-order}).

The first-order computation shows where the factors $A_f$ and $\fc{p_f}1$ come from
(Section~\ref{ss:first-order}). For $\epsp=0$ and to first order in $\eps$, the Jacobian matrix
in~\eqref{eq:Delta} is proportional to the Hessian of the Melnikov potential. Its slow entry,
the second derivative in $\alpha_1$, is proportional to $A_f$ and is not small. Its fast entry,
the second derivative in $\alpha_2$, is proportional to $\fc{p_f}1$ and is exponentially small.
The determinant is their product to leading order. For $\epsp$ of order one the Melnikov
potential is replaced by an integral of the same type along the trajectories of the pendulum
forced by the strong term alone. The slow entry is the same to leading order. The fast entry
need not be exponentially small any more, but the determinant still is. Its leading term is the
product of the slow entry and the part of the fast entry that comes from the odd harmonics of
$f$ in $\alpha_2$. In this part the number $4\pi\fc{p_f}1=G_{M,p_f}(0)$ is replaced by
$G_{M,p_f}(\mu)$, and this is the only effect of the strong term on the leading term. In the
example of Theorem~\ref{thm:example} the ratio $G(\mu)/G(0)$ equals $1+T_1\mu^2+\dots$ and is
larger than one for $\mu\ne0$.

The proof has three parts, and Section~\ref{ss:scheme} contains a detailed outline. First, we
describe the two manifolds by generating functions and express $\Delta$ through the Hessian of
their difference at the homoclinic point. This Hessian is analytic in the weak amplitude $\eps$,
and the leading term of its determinant is quadratic in $\eps$. It is the determinant of the
Hessian of a Melnikov-type integral taken along the trajectories of the pendulum forced by the
strong term alone. Second, we prove exponentially small bounds for this integral and for the
remainder. The difference of the generating functions satisfies a linear transport equation,
which we straighten on a segment of imaginary time that ends at a distance of order $\sqrt\eta$
from the singularities of the separatrix. The estimates of the manifolds near the singularities
that this requires are known for the pendulum forced by the strong term alone~\cite{BFGS12},
and we extend them to the full problem by comparison. Third, we find the limit of the
normalized fast entry and identify it with the response of the Stokes constant of the inner
equation to adding the profile $p_f$ to the strong term, in the sense of
Theorem~\ref{thm:series}. We obtain the response by interpolation in a one-parameter family of
profiles. This argument does not use the analyticity of the Stokes constant in the profile.
Finally, we construct the solutions of the inner equation as Laplace integrals, order by order
in the amplitude of the strong term, and the response becomes a residue in the Borel plane,
which gives the series~\eqref{eq:G}.

\subsection{Relation to known results}\label{ss:literature}

\paragraph{Systems with three time scales} Hamiltonians that describe a pendulum coupled to a
slow and a fast rotator were introduced in~\cite{CG94} in a study of Arnold diffusion in the
d'Alembert model of the precession of a planet. The lower bound for the splitting stated there
was incorrect~\cite{CG98}. A correct bound was obtained in~\cite{GGM99} with techniques
developed in~\cite{Gallavotti94} for systems with two time scales. Without the strong term, the
Hamiltonian~\eqref{eq:H} belongs to the class considered in~\cite{GGM99}. It is proved there
for small perturbations that the determinant is exponentially small and, generically, agrees to
leading order with the prediction of first-order perturbation theory (see
also~\cite{GGM99b}). For the example of~\cite{GGM99} one obtains
$\det D=32\pi\,\eta^{-1/2}\eps^2e^{-\frac\pi2\eta^{-1/2}}$. Applications to Arnold diffusion in
systems with several time scales are given
in~\cite{GGM00,GallavottiMPAG,BertiBolle02,Procesi03}.

A term that is not small and depends only on the fast angle, like the strong term
in~\eqref{eq:H}, arises in problems of celestial mechanics~\cite{Gallavotti98}, in particular
in the precession problem~\cite{GGM99,Mastropietro01}, and the theorem of~\cite{GGM99} does not
cover it. For such a term it is noted in~\cite[Section~8]{GGM99}, with reference
to~\cite{CG94}, that the splitting is analytic in its amplitude and, generically, can vanish
only for finitely many values of the amplitude in the domain of convergence of the perturbation
series when the other parameters are fixed. No estimate of the size of the splitting was given
(see also~\cite{GGM00,GallavottiMPAG}). The absence of a lower bound is also noted in the
review~\cite{Mastropietro01}. This is the origin of the problem~\cite{Gallavotti98}.

\paragraph{Two fast frequencies} A pendulum with a quasi-periodic forcing in which both
frequencies are fast is considered in~\cite{Simo94,DGJS97}. An averaging theorem with an
exponentially small remainder is proved in~\cite{Simo94}. For an example in which the forcing
is not small, the splitting is computed there numerically and compared with the Melnikov
function evaluated after one step of averaging. It is proved in~\cite{DGJS97} that for a
sufficiently small perturbation, under a lower bound on a sequence of its Fourier coefficients,
the splitting is exponentially small and is correctly predicted by the Melnikov function. In
the problem considered here only one frequency is fast, and the term that is not small depends
on the fast angle alone.

\paragraph{The periodic problem} If there is no weak perturbation ($\eps=0$), the forcing
depends on a single fast angle. For a pendulum with a fast periodic forcing of small amplitude
the asymptotics of the splitting is given by the Melnikov integral~\cite{DS92} (see the
survey~\cite{GL01}). When the amplitude is not small, the leading coefficient is determined by
an equation without small parameters, which is called the inner equation or the reference
system. This approach goes back to the work of Lazutkin on the standard
map~\cite{Lazutkin84}, for which a complete proof is given in~\cite{Gelfreich99}. For a
pendulum with a rapidly oscillating suspension point, asymptotic formulas of this kind are
stated in~\cite{Treschev96,Gelfreich97}, where the proofs are referred to other texts.
According to these papers, the coefficient is an entire function of the amplitude.
A proof for a pendulum with an additive forcing is given in~\cite{GOS10}. In~\cite{BFGS12} the
asymptotics is proved for a class of Hamiltonians that includes the periodic problem of this
paper. The proof covers the singular case, to which our problem belongs, for any fixed
amplitude, and the leading coefficient is an entire function of the amplitude. In this
connection the problem statement~\cite{Gallavotti98} mentions the
works~\cite{Neishtadt84,Lazutkin84,Gelfreich97}.

For high-frequency perturbations of the pendulum, the equation without small parameters is
written in~\cite[Section~4]{Gelfreich97} as the reference system $u''=(1+g(s))\,e^u$ with a
periodic function $g$. Equation~\eqref{eq:newton} is this system written for the function
$e^{-u/2}$, up to normalization (Remark~\ref{rem:newton-origin}). In
the Hamilton--Jacobi form the inner equation was studied for a wide class of systems
in~\cite{Baldoma06}, and its resurgent structure for a single-harmonic profile
in~\cite{OSS03}. Borel summation was applied earlier to the splitting of the separatrices of
the H\'enon map~\cite{GS01}. For the pendulum with the perturbation
$\mu g(\theta)(\cos q-1)$, the first Taylor coefficient of the Stokes constant is
$\chi_1=4\pi\fc g1$, the value given by the Melnikov integral~\cite{BFGS12,BMS24} (see
also~\cite[Lemma~3]{Gelfreich97}). The case $\fc g1=0$ is studied in~\cite{BMS24}. There the
first nonzero Taylor coefficient is expressed through solutions of the inner equation, an
explicit formula is given when it is the second one, and a numerical method is proposed for
the others.

\paragraph{The leading coefficient as a function of the system} The following results concern
the dependence of such coefficients on the system and their nonvanishing. The analytic
dependence of the constant of the reference system on a parameter of the function $g$ is stated
in~\cite[Section~4.3]{Gelfreich97} without proof. For the inner equation of a different
problem it is proved in~\cite{GGSZ25} that the Stokes constant is an analytic function of the
nonlinearity and does not vanish on an open dense set of nonlinearities, and a formula that
expresses it as a series is conjectured there. For the H\'enon map the nonvanishing of the
splitting constant is derived in~\cite{GS01} from the positivity of a Borel transform. For two
other problems the nonvanishing of the Stokes constant is proved with computer
assistance~\cite{BCGS23,BCGG25}. For the pendulum with an oscillating suspension point the
coefficient is written in~\cite{Treschev96} through an entire function of a variable $w$ that
is proportional to the product of two amplitudes of the forcing. Nine Taylor coefficients of
this function are computed there numerically. The computation uses a recursive system for the
coefficients of a series in the inner variable, which is integrated term by term. The
resulting numbers are close to the Taylor coefficients of $2\sinh\sqrt{2w}/\sqrt{2w}$. A
related observation is reported in~\cite[Section~2.2]{Gelfreich97}. A function of the same
form enters the lower bounds of Section~\ref{ss:growing}.

\paragraph{What is new} We use two statements on the periodic problem,
Propositions~\ref{prop:E1} and~\ref{prop:E2} of Section~\ref{ss:external}, which go back
to~\cite{BFGS12,Baldoma06}. The description of the manifolds by generating functions whose
difference satisfies a linear equation, and the straightening of the characteristic field of
this equation, go back to~\cite{Sauzin01}. This method was applied to the periodic problem
in~\cite{GOS10,BFGS12}. To the best of our knowledge, the following results are new:
\begin{itemize}
\item the asymptotics of $\Delta$ for an amplitude $\epsp$ of order one, and the expression of
its coefficient $-8A_fG_{M,p_f}(\mu)$ through the response of the Stokes constant;
\item the bound for the first variation of the splitting of the periodic problem on a complex
disc of background amplitudes, and the computation of its limit by interpolation in the
profile;
\item formula~\eqref{eq:G}, which expresses the response through the values at $s=1$ of Borel
transforms given by an explicit recursion;
\item the nonnegativity of the terms of the series~\eqref{eq:G} and the lower bounds that follow
from it.
\end{itemize}

\subsection{Structure of the paper}\label{ss:notation}

In Section~\ref{sec:setting} we reduce the problem to dimensionless form, state
Proposition~\ref{thm:Qprime} on the determinant $F$ and derive Theorem~\ref{thm:main} from it. In
Section~\ref{sec:periodic} we introduce notation and state two propositions on the periodic
problem. Section~\ref{ss:scheme} contains an outline of the proof, with the dependencies between
its steps shown in Table~\ref{tab:deps}. Sections~\ref{sec:manifolds}--\ref{sec:variation}
contain lemmas, and
Section~\ref{sec:determinant} contains the proof of Proposition~\ref{thm:Qprime}. We prove
Theorems~\ref{thm:series} and~\ref{thm:example} and Proposition~\ref{prop:nondeg} in
Section~\ref{sec:borel}, and Theorem~\ref{thm:cos2} in Section~\ref{sec:negative}.
Section~\ref{sec:discussion} contains a discussion and open questions. In
Appendix~\ref{app:bfgs} we give a self-contained proof of Proposition~\ref{prop:E1} and explain
how Proposition~\ref{prop:E2} follows from the results of~\cite{BFGS12} and~\cite{Baldoma06}, and
Appendix~\ref{app:certificate} contains the data of the computer-assisted part of
Section~\ref{sec:negative}.
\section{The dimensionless problem}\label{sec:setting}

In this section we reduce the problem to dimensionless form, state
Propositions~\ref{thm:Qprime} and~\ref{thm:Q} on the splitting determinant of the dimensionless
problem, and derive Theorem~\ref{thm:main} from the first of them.

\subsection{Dimensionless variables}\label{ss:dimensionless}

It suffices to consider the case $\bom_2>0$. Indeed, the canonical change
$(\alpha_2,A_2)\mapsto(-\alpha_2,-A_2)$ takes the Hamiltonian~\eqref{eq:H} to a Hamiltonian of
the same form with $-\bom_2$ in place of $\bom_2$ and with the function $f(\alpha_1,-\alpha_2)$
in place of $f(\alpha_1,\alpha_2)$. The function $m(2\alpha_2)$ does not change, since $m$ is
even. The Jacobian matrix in~\eqref{eq:Delta} is then multiplied on the left and on the right by
$\diag(1,-1)$, so $\Delta$ does not change. Neither do the quantities \eqref{eq:mu-M},
\eqref{eq:Af-pf}. For $M$ and $A_f$ this is obvious, and $p_f$ does not change because
$f(0,-\theta)=f(0,\theta)$. From now on we assume that $\bom_2>0$.

We pass to the dimensionless time $gt$, which we again denote by $t$, and a prime denotes the
derivative with respect to it. We set $q=\varphi$, $x=\alpha_1$, $\theta=\alpha_2$,
$\alpha=(x,\theta)$ and
\begin{equation}\label{eq:params}
\mu=\frac{\epsp}{Jg^2},\qquad\tau=\frac{\eps}{Jg^2},\qquad h=\frac{g\sqrt\eta}{\bom_2},\qquad
\nu=\frac{\bom_1\sqrt\eta}{g}.
\end{equation}
Equation~\eqref{eq:pendulum} takes the form
\begin{equation}\label{eq:A0}
q''=K(\alpha+\omega t)\,\sin q,\qquad K(x,\theta)=1+\mu M(\theta)+\tau f(x,\theta),\qquad
\omega=(\nu,\,1/h).
\end{equation}
Here and below the letter $\omega$ denotes the dimensionless frequency vector, that is, the
vector $\omega$ from~\eqref{eq:H} divided by $g$. Note that $|\nu|=\nu_0h$, where
$\nu_0=|\bom_1|\bom_2/g^2$, so the slow frequency is small together with $h$.

The dimensionless momentum and actions are $p=q'=I/(Jg)$ and $\tilde A=A/(Jg)$. From
$\dot A_i=-\partial H/\partial\alpha_i$ we get
\begin{equation}\label{eq:action-dot}
\tilde A_i'=\partial_{\alpha_i}K(\alpha+\omega t)\,\bigl(1-\cos q\bigr),\qquad i=1,2,
\end{equation}
so formula~\eqref{eq:Au} becomes
\begin{equation}\label{eq:Au-dimless}
\tilde A^u_i(\alpha)=\int_{-\infty}^0\partial_{\alpha_i}K(\alpha+\omega t)\,\bigl(1-\cos q(t)\bigr)\,dt,
\end{equation}
where $q(t)$ is the trajectory on the unstable manifold with $q(0)=\pi$ and phases $\alpha$ at
$t=0$, and similarly for $\tilde A^s$. We define the \emph{splitting determinant} of the
dimensionless problem as
\begin{equation}\label{eq:F-def}
F=\det\partial_\alpha\bigl(\tilde A^u-\tilde A^s\bigr)\Big|_{\alpha=0}.
\end{equation}
Then
\begin{equation}\label{eq:Delta-F}
\Delta=(Jg)^2F .
\end{equation}

\paragraph{The general form of the equation} The proofs deal with an equation of a somewhat more
general form:
\begin{equation}\label{eq:A}
\begin{gathered}
q''=K(\alpha+\omega t)\,\sin q,\qquad K(\alpha)=1+zM(\theta)+V(\alpha),\\
\alpha=(x,\theta),\quad\omega=(\nu,1/h).
\end{gathered}
\end{equation}
Here $0<h\le h_0$ and $|\nu|\le\nu_0h$. The function $M$ is a real even trigonometric polynomial
with zero mean. We call it the \emph{background}, and $z$ is its amplitude, which may be
complex, $|z|\le R_0$. The function $V=\sum_i\eps_if_i(\alpha)$ is the \emph{weak perturbation},
with fixed real even trigonometric polynomials $f_i$ and amplitudes $\eps_i$, possibly complex as
well. We write $\eps=(\eps_i)$ for the collection of these amplitudes and set
$|\eps|=\sum_i|\eps_i|$. The amplitudes $\eps_i$ are dimensionless, and neither they nor the
collection $\eps$ should be confused with the parameter $\eps$ of the Hamiltonian~\eqref{eq:H}.

In Gallavotti's problem the amplitude $z=\mu$ is real, the background $M(\theta)=m(2\theta)$ has
period $\pi$, and the weak perturbation consists of the single term $\tau f$. We need complex
amplitudes because the proofs use analyticity in the amplitudes and Cauchy's inequalities. In
the theorems the background has period $\pi$, while the lemmas of
Sections~\ref{sec:manifolds}--\ref{sec:stokes} do not require this. Definition~\eqref{eq:F-def}
refers to real amplitudes. In Section~\ref{ss:detdef} we write the determinant $F$ as the
Hessian determinant of a generating function, and this expression defines its analytic
continuation to complex amplitudes.

\subsection{Propositions on the splitting determinant}\label{ss:thmQ}

We set
\begin{equation}\label{eq:a-Lh}
a=\frac\pi2,\qquad L_h=\frac1h\,e^{-a/h}.
\end{equation}

\begin{proposition}[general even perturbations]\label{thm:Qprime}
Let $M(\theta)$ be a real even trigonometric polynomial with zero mean and period $\pi$, let
$f(x,\theta)$ be a real even trigonometric polynomial, and let $A_f$ and $p_f$ be defined
by~\eqref{eq:Af-pf}. Let $K=1+\mu M(\theta)+\tau f(x,\theta)$ in equation~\eqref{eq:A}, and let
the numbers $\nu_0>0$ and $R_0>0$ be fixed. There is an absolute constant $P$ (one can take
$P=98$) such that, as $h\to0$,
\begin{equation}\label{eq:Qprime}
F=\bigl(-8A_f\,G_{M,p_f}(\mu)+o(1)\bigr)\,\tau^2L_h
\end{equation}
uniformly in real $\mu\in[-R_0,R_0]$, in $\nu$ with $|\nu|\le\nu_0h$ and in $\tau$ with
\begin{equation}\label{eq:tau-range}
e^{-a/(2h)}\le|\tau|\le h^{3P+2}.
\end{equation}
Here $G_{M,p}$ is the entire function of Theorem~\ref{thm:series}. The amplitude $\tau$ may be
complex.
\end{proposition}

\begin{proposition}[two independent weak amplitudes]\label{thm:Q}
Let $M=\cos4\theta$ in equation~\eqref{eq:A}, let the amplitude $z=\mu$ be real, and let
$V=s\cos x+r\cos\theta$. Then for every fixed $\mu$, as $h\to0$,
\begin{equation}\label{eq:Q}
F=\bigl(8\,G(\mu)+o(1)\bigr)\,s\,r\,L_h
\end{equation}
uniformly in $\nu$ with $|\nu|\le\nu_0h$ and in complex amplitudes $s$, $r$ such that
\begin{equation}\label{eq:sr-range}
s\ne0,\qquad|s|,|r|\le h^{3P+2},\qquad|r|\ge e^{-2a/h}.
\end{equation}
Here $G=G_{\cos4\theta,\cos\theta}$ is the function of Theorem~\ref{thm:example}, and $P$ is the
constant of Proposition~\ref{thm:Qprime}.
\end{proposition}

We prove both propositions in Section~\ref{sec:determinant}. Both include the statement that for
small $h$ the manifolds exist and the determinant $F$ is defined, which we prove in
Sections~\ref{sec:manifolds} and~\ref{sec:generating}. Proposition~\ref{thm:Qprime} with
$M=\cos4\theta$ and $f=\cos x+\cos\theta$ gives formula~\eqref{eq:Q} for $s=r=\tau$, and
Proposition~\ref{thm:Q} shows that the amplitudes of the slow and the fast harmonic of the weak
perturbation can be varied independently.

\begin{remark}[on the order of exposition]\label{rem:G-order}
In the proof of Proposition~\ref{thm:Qprime} the function $G_{M,p}$ arises as the limit of the
normalized first variation (Section~\ref{ss:variation-limit}), which also gives
formula~\eqref{eq:response} for its Taylor coefficients. We prove in Section~\ref{sec:borel}
that this function is given by the series~\eqref{eq:G}.
\end{remark}

\begin{proof}[Deduction of Theorem~\ref{thm:main} from Proposition~\ref{thm:Qprime}]
Let $c>c_0:=(3P+2)/2$, let $\mathcal K\subset\RR$ be a compact set, and let
$R_0=\max_{\epsp\in\mathcal K}|\epsp|/(Jg^2)$. For $\eps=\eta^c$, \eqref{eq:params} gives
$\sqrt\eta=\bom_2h/g$ and
\[
\tau=\frac{\eta^c}{Jg^2}=\frac{(\bom_2/g)^{2c}}{Jg^2}\,h^{2c}.
\]
Since $2c>3P+2$, condition~\eqref{eq:tau-range} holds for $h$ small enough. Moreover,
$|\nu|=\nu_0h$. By Proposition~\ref{thm:Qprime},
$F=\bigl(-8A_fG_{M,p_f}(\mu)+o(1)\bigr)\tau^2L_h$ uniformly in $\mu\in[-R_0,R_0]$, that is, in
$\epsp\in\mathcal K$. It remains to substitute this into~\eqref{eq:Delta-F}:
\[
(Jg)^2\,\tau^2L_h=(Jg)^2\cdot\frac{\eta^{2c}}{J^2g^4}\cdot\frac{\bom_2}{g\sqrt\eta}\,
\exp\Bigl[-\frac{\pi\bom_2}{2g\sqrt\eta}\Bigr]
=\frac{\bom_2}{g^3}\,\eta^{2c-\frac12}\exp\Bigl[-\frac{\pi\bom_2}{2g\sqrt\eta}\Bigr].
\]
For $P=98$ we get $c_0=148$.
\end{proof}

The same computation with an arbitrary $\eps$ in place of $\eta^c$ gives
Remark~\ref{rem:range}. Theorems~\ref{thm:example} and~\ref{thm:cos2} follow from
Theorem~\ref{thm:main}. Indeed, for $f=\cos\alpha_1+\cos\alpha_2$ we have $A_f=-1$ and
$p_f=\cos\theta$, and we prove the lower bounds for the functions $G$ and $G_B$ in
Sections~\ref{ss:positivity} and~\ref{sec:negative}. Since these functions are bounded away
from zero on the real axis, the remainder $o(1)$ in~\eqref{eq:main} can be written as a factor
$1+o(1)$.

\begin{remark}[Proposition~\ref{thm:Q} in the original variables]\label{rem:two-amplitudes}
Suppose that in the example of Theorem~\ref{thm:example} the weak perturbation has the form
$\eps_1\cos\alpha_1+\eps_2\cos\alpha_2$ with independent amplitudes. Then $s=\eps_1/(Jg^2)$,
$r=\eps_2/(Jg^2)$, and Proposition~\ref{thm:Q} together with~\eqref{eq:Delta-F} gives
\[
\Delta=\frac{8\bom_2}{g^3}\,G(\mu)\,\eps_1\eps_2\,\eta^{-\frac12}
\exp\Bigl[-\frac{\pi\bom_2}{2g\sqrt\eta}\Bigr]\bigl(1+o(1)\bigr)
\]
for every fixed $\epsp$, provided $s$ and $r$ satisfy conditions~\eqref{eq:sr-range}.
\end{remark}

\section{Preliminaries on the periodic problem}\label{sec:periodic}

In this section we introduce notation related to the unperturbed separatrix and state
Propositions~\ref{prop:E1} and~\ref{prop:E2} on the periodically forced pendulum, on which all
subsequent constructions rely.

\subsection{The unperturbed separatrix, domains and notation}\label{ss:separatrix}

For $z=0$, $V=0$ equation~\eqref{eq:A} is the pendulum equation $q''=\sin q$. Its separatrix with
$p>0$ is given by the functions
\begin{equation}\label{eq:q0}
q^0(u)=4\arctan e^u,\qquad\psi_0(u)=\frac2{\cosh u},
\end{equation}
for which
\begin{equation}\label{eq:q0-identities}
\begin{gathered}
(q^0)'=\psi_0,\qquad\sin q^0=-\frac{2\sinh u}{\cosh^2u},\\
1-\cos q^0=\frac{\psi_0^2}2,\qquad q^0(-u)=2\pi-q^0(u).
\end{gathered}
\end{equation}
We call the variable $u$, the time along the separatrix, the \emph{label}. The value $u=0$
corresponds to $q=\pi$. The function $\psi_0$ is meromorphic, and its only poles in the strip
$|\Im u|\le a$ are the points $u=\pm ia$. The function $\cosh u/(u^2+a^2)$ is entire and has no
zeros for $|\Im u|<3a$.

For $u$ with $|\Im u|<a$ we set $\delta(u)=a-|\Im u|$. For numbers $\kappa>0$ and $\gamma>0$ we
define the set of labels
\begin{equation}\label{eq:U}
\cU^{\gamma}_\kappa=\{u:\ \delta(u)\ge\kappa h,\ \ |\Re u|\le\gamma\,\delta(u)\}.
\end{equation}
It contains the segment $\{iy:\ |y|\le a-\kappa h\}$ of the imaginary axis, whose endpoints are
at distance $\kappa h$ from the singularities (Fig.~\ref{fig:domains}).

We denote the complex neighborhoods of the circle and of the torus by
$\TT_\rho=\{\theta\in\CC/2\pi\ZZ:\ |\Im\theta|<\rho\}$ and $\TT^2_\rho$ (the same for a pair of
angles). For a trigonometric polynomial $g$ with zero mean, $P_g$ denotes its periodic primitive
with zero mean, so that $P_g'=g$. The constants $C,C_1,\dots$ do not depend on $h$, on the
amplitudes $z$, $\eps_i$ or on the point, and we indicate their dependence on the other data
where they appear.

We write Fourier coefficients with an index in square brackets ($\fc pk$, $\fc f{j,\ell}$) and
the mean over an angle as $\langle\,\cdot\,\rangle_\theta$, and $[z^k]\,\Phi$ denotes the
coefficient of $z^k$ in the Taylor series of $\Phi$ at zero. The letter $\cD$ denotes the
difference of the generating functions, $F$ the splitting determinant of the dimensionless
problem, and $\sD$ the operator $\partial_\xi-i\partial_\theta$. The letter $x$ is the slow angle
$\alpha_1$ in Sections~\ref{sec:setting}--\ref{sec:determinant} and~\ref{sec:discussion}, while
in formula~\eqref{eq:newton} and in Sections~\ref{sec:borel}, \ref{sec:negative} it is the
variable of the inner equation. The two meanings never occur together.

\subsection{The Hamilton--Jacobi equation and statements (E1), (E2)}\label{ss:external}

Let $g(\theta)$ be a real even trigonometric polynomial with zero mean, and let $z$ be a complex
amplitude. (The constant $g$ of the Hamiltonian~\eqref{eq:H} does not enter the dimensionless
equations, and up to the end of Section~\ref{sec:negative} the letter $g$ denotes only the
profile.) The \emph{periodic problem with profile $zg$} is equation~\eqref{eq:A} without the weak
perturbation and with background $zg$:
\begin{equation}\label{eq:per}
q''=\bigl(1+z\,g(\theta+t/h)\bigr)\sin q,
\end{equation}
that is, the system with Hamiltonian $p^2/2+(\cos q-1)\bigl(1+zg(\theta)\bigr)$,
$\theta'=1/h$. This is a pendulum with a fast periodic perturbation, which need not be small.
Such systems are studied in~\cite{BFGS12}.
Equation~\eqref{eq:per} is system~(7) of that paper with $H_0=y^2/2+\cos x-1$ and the
perturbation $z\,g(\tau)(\cos x-1)$, the small parameter being $h$. It belongs to the case called
singular in~\cite{BFGS12}. We give the correspondence of notation and the verification of the
hypotheses in Appendix~\ref{app:bfgs}.

We seek the invariant manifolds of the periodic orbit $q=p=0$ as graphs
$p=\partial_qS(q,\theta)$. In the variables $(u,\theta)$, where $q=q^0(u)$, the function
$T(u,\theta)=S(q^0(u),\theta)$ satisfies the Hamilton--Jacobi equation
\begin{equation}\label{eq:HJ}
\frac1h\,\partial_\theta T+\frac{(\partial_uT)^2}{2\psi_0^2}-\bigl(1+zg(\theta)\bigr)\frac{\psi_0^2}2=0
\end{equation}
(equation~(47) in~\cite{BFGS12}). For $z=0$ its solution is the function $T_0$ with
$\partial_uT_0=\psi_0^2$. We write
\begin{equation}\label{eq:E-def}
\partial_uT=\psi_0^2\,(1+E).
\end{equation}
Then the momentum on the graph equals $p=\psi_0(1+E)$, and the label of a point moving on the
manifold obeys the equations
\begin{equation}\label{eq:label-motion}
u'=1+E(u,\theta),\qquad\theta'=\frac1h .
\end{equation}

We shall use two statements about the solutions of equation~\eqref{eq:HJ}.

\begin{proposition}[E1; cf.~{\cite[Theorems~4.3 and~4.4]{BFGS12}}]\label{prop:E1}
Fix $R_0>0$ and a profile $g$. There exist $\sigma>0$ and $\varrho_1>0$ and, for every
$\varrho_2>\varrho_1$, constants $\kappa_1,h_1,C_E>0$ such that for $0<h<h_1$
equation~\eqref{eq:HJ} has a solution $T^u$ analytic in $(u,\theta,z)$ on
$(D^{\rm loc}\cup D^{\rm out})\times\TT_\sigma\times\{|z|<R_0+1\}$, where
\begin{equation}\label{eq:domains}
\begin{gathered}
D^{\rm loc}=\{u:\ \Re u<-\varrho_1\},\\
D^{\rm out}=\{u:\ |\Im u|<a-\kappa_1h-\Re u,\ \ \Re u>-\varrho_2\},
\end{gathered}
\end{equation}
and for $|z|\le R_0$ the corresponding function $E^u$ satisfies the estimates
\begin{equation}\label{eq:E1}
\begin{gathered}
|E^u(u,\theta)|\le C_E|z|h\quad(u\in D^{\rm loc}),\\
|E^u(u,\theta)|\le C_E|z|h\,\frac{|\cosh u|^2}{|u^2+a^2|^3}\quad(u\in D^{\rm out}).
\end{gathered}
\end{equation}
\end{proposition}

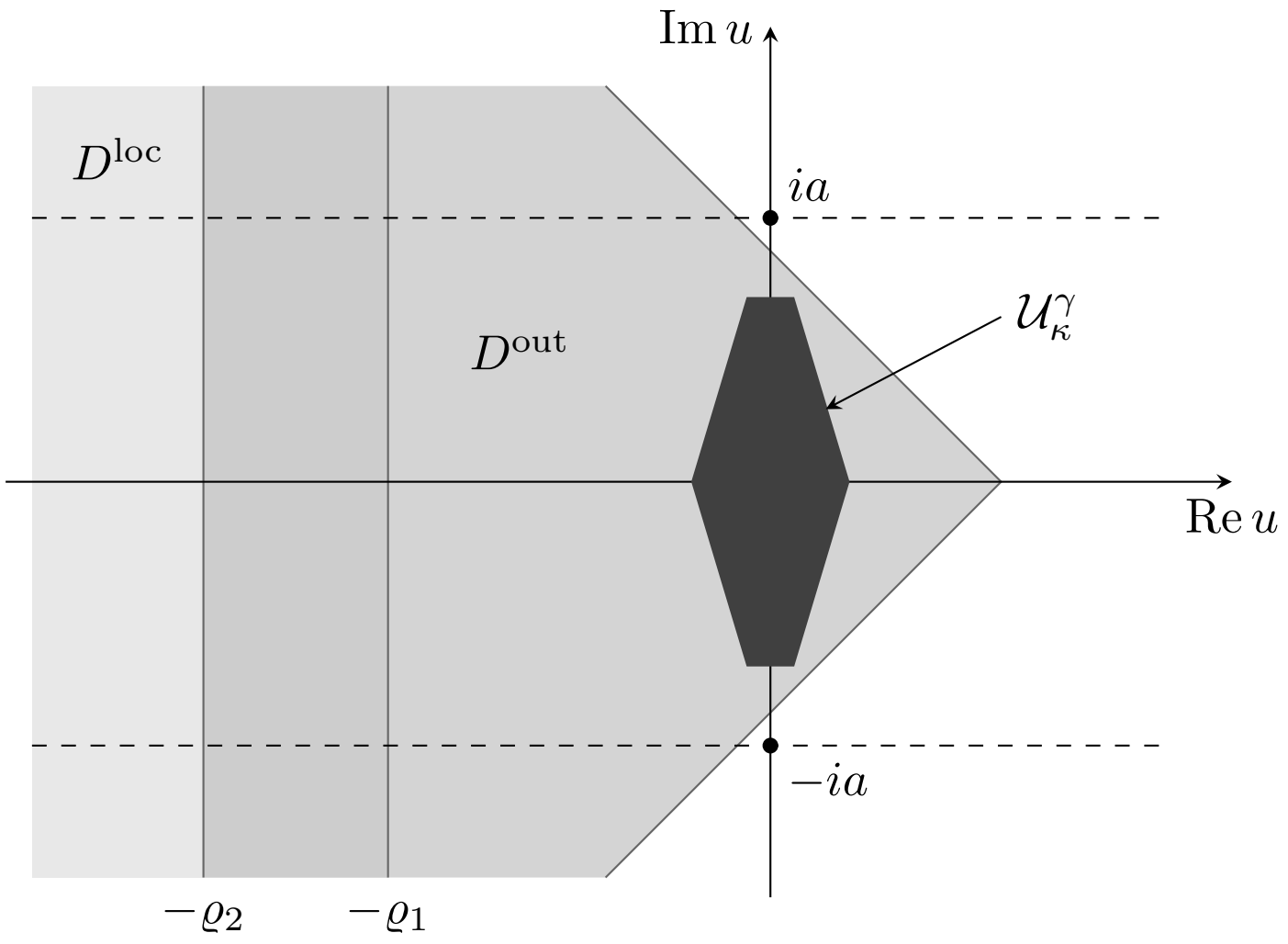
\begin{figure}[ht]
\centering
\begin{tikzpicture}[>=stealth,scale=0.95]
  \pgfmathsetmacro{\ha}{2.0} %
  \pgfmathsetmacro{\kone}{0.25} %
  \pgfmathsetmacro{\kap}{0.6} %
  \pgfmathsetmacro{\rone}{2.9} %
  \pgfmathsetmacro{\rtwo}{4.3} %
  \pgfmathsetmacro{\yt}{3.0} %
  \pgfmathsetmacro{\gam}{0.3} %
  \pgfmathsetmacro{\apex}{\ha-\kone}
  \fill[black!9] (-5.6,-\yt) rectangle (-\rone,\yt);
  \fill[black!24,opacity=0.7] (\apex,0) -- (\apex-\yt,\yt) -- (-\rtwo,\yt) -- (-\rtwo,-\yt)
        -- (\apex-\yt,-\yt) -- cycle;
  \draw[black!60] (\apex-\yt,\yt) -- (\apex,0) -- (\apex-\yt,-\yt);
  \draw[black!60] (-\rtwo,-\yt) -- (-\rtwo,\yt);
  \draw[black!60] (-\rone,-\yt) -- (-\rone,\yt);
  \draw[dashed] (-5.6,\ha) -- (3.0,\ha);
  \draw[dashed] (-5.6,-\ha) -- (3.0,-\ha);
  \draw[->] (-5.8,0) -- (3.5,0) node[below] {$\Re u$};
  \draw[->] (0,-\yt-0.15) -- (0,\yt+0.45) node[left] {$\Im u$};
  \fill[black!75] (\gam*\kap,\ha-\kap) -- (\gam*\ha,0) -- (\gam*\kap,-\ha+\kap)
        -- (-\gam*\kap,-\ha+\kap) -- (-\gam*\ha,0) -- (-\gam*\kap,\ha-\kap) -- cycle;
  \fill (0,\ha) circle (1.7pt) node[above right] {$ia$};
  \fill (0,-\ha) circle (1.7pt) node[below right] {$-ia$};
  \node at (-4.95,2.45) {$D^{\rm loc}$};
  \node at (-1.9,1.0) {$D^{\rm out}$};
  \node[below] at (-\rone,-\yt) {$-\varrho_1$};
  \node[below] at (-\rtwo,-\yt) {$-\varrho_2$};
  \draw[<-] (0.42,0.55) -- (1.75,1.25) node[right] {$\cU^{\gamma}_{\kappa}$};
\end{tikzpicture}
\caption{Domains in the plane of the complex label $u$ (schematic, not to scale). The light
half-plane is $D^{\rm loc}$ and the darker wedge is $D^{\rm out}$. They overlap for
$-\varrho_2<\Re u<-\varrho_1$. The singularities $u=\pm ia$ lie outside $D^{\rm out}$, at a
distance of order $\kappa_1h$ from its boundary. The dark rhombus is the set of labels
$\cU^{\gamma}_{\kappa}$.}\label{fig:domains}
\end{figure}

The function $T^u$ is the generating function of the unstable manifold (Lemma~\ref{lem:local},
item~4). It is defined up to an additive constant, and we use only its derivatives below. The
domain $D^{\rm out}$ approaches
the singularities $u=\pm ia$ to within a distance of order $\kappa_1h$
(Fig.~\ref{fig:domains}), and near them the second estimate~\eqref{eq:E1} takes the form
$|E^u|\le C|z|h/|u\mp ia|$. The constants $\kappa_1,h_1,C_E$ are chosen after $\varrho_2$ and do
not depend on $z$ and $h$.

Below we refer to Proposition~\ref{prop:E1} as statement~(E1). For the pendulum it corresponds
to Theorems~4.3 and~4.4 of~\cite{BFGS12}. The estimate in $D^{\rm loc}$ follows from the weighted
estimate of Proposition~6.4 of that paper, from which Theorem~4.3 is derived, and the estimate in
$D^{\rm out}$ follows from Theorem~4.4. The domain $D^{\rm out}$ is domain~(35)
of~\cite{BFGS12} with the angle $\beta_1=\pi/4$. We prove Proposition~\ref{prop:E1} in
Section~\ref{ss:A-direct} by a short self-contained argument, in which the dependence on the
complex amplitude $z$ is explicit.

\begin{proposition}[E2; {\cite[Theorem~4.16]{BFGS12}}
and {\cite[Theorem~3.1]{Baldoma06}}]\label{prop:E2}
For the same $R_0$ and $g$, fix angles $\beta_0\in(0,\pi/4)$ and $\beta_2\in(\pi/4,\pi/2)$. There
exists $\kappa_6\ge\kappa_1$ such that for all $|z|<R_0$, $\theta\in\TT_\sigma$ and every
$\zeta$ in the sector
\begin{equation}\label{eq:sector}
S^u=\{\zeta:\ \Im\zeta<-\tan\beta_2\,\Re\zeta-\kappa_6,\ \ \Im\zeta<-\tan\beta_0\,\Re\zeta-\kappa_6\}
\end{equation}
the limit
\begin{equation}\label{eq:E2}
\Phi^u_0(\zeta,\theta;z)=\lim_{h\to0}h^2\,\partial_uT^u(ia+h\zeta,\theta;z)
\end{equation}
exists.
\end{proposition}

Since $\kappa_6\ge\kappa_1$ and $\beta_0<\pi/4<\beta_2$, the sector $S^u$ lies in the half-plane
$\{\Im\zeta<-\Re\zeta-\kappa_1\}$, and the points $ia+h\zeta$ with $\zeta\in S^u$ belong to
$D^{\rm out}$ for small $h$. The limit~\eqref{eq:E2} describes the manifold in the neighborhood
of the singularity $u=ia$ stretched by the factor $1/h$. Below we refer to
Proposition~\ref{prop:E2} as statement~(E2). It follows from Theorem~4.16 of~\cite{BFGS12}
together with Theorem~3.1 of~\cite{Baldoma06}, which is quoted in~\cite{BFGS12} as
Theorem~4.12(1). There the limit is expressed through the solution $\psi^u_0$ of the inner
equation as $\Phi^u_0=-4\,\partial_\zeta\psi^u_0$. We explain in Section~\ref{ss:A-E2} how these
theorems apply to the pendulum. We do not need uniformity in the amplitude in~\eqref{eq:E2},
because we use only the pointwise limit.

Of the results of~\cite{BFGS12}, our proofs use only those on which statement~(E2) is based. In
particular, we do not need the statements of that
paper on the stable manifold, because we obtain the stable manifold from the unstable one by
reversibility (Section~\ref{ss:local}). Nor do we use the final asymptotic formulas
of~\cite{BFGS12} for the splitting of the periodic problem.

\subsection{The order in which the constants are chosen}\label{ss:constants}

We choose the constants in the lemmas of Sections~\ref{sec:manifolds}--\ref{sec:variation} in a
definite order. The order matters because the domains in (E1) depend on $\varrho_2$, while
$\varrho_2$ is determined by the length of the time interval over which the manifolds are
continued.

Let $R_0$ and the profiles be given. Statement~(E1) provides $\sigma$ and $\varrho_1$. Then we
choose the following in this order:
\begin{enumerate}[label=\textup{\arabic*)}]
\item the size $q_1$ of the neighborhood in which the local manifolds are constructed
(Lemma~\ref{lem:local});
\item a number $\varrho_1'\ge\varrho_1$ to the left of which the graph from (E1) coincides with
the local manifold (Lemma~\ref{lem:local}, item~4);
\item the time $t_*=\tfrac87(\varrho_1'+2)$ and $\varrho_2=2t_*+1$;
\item the constants $\kappa_1$, $h_1$, $C_E$ from (E1) corresponding to this $\varrho_2$;
\item the constant $C_1$ in the estimates of the background trajectories
(Lemma~\ref{lem:geometry});
\item the constant $\kappa$ that determines the distance $\kappa h$ to the singularities
(Lemmas~\ref{lem:geometry}, \ref{lem:continuation} and the conditions on $\kappa_0=4\kappa$ in
Sections~\ref{sec:straightening}, \ref{sec:stokes});
\item the remaining constants and the thresholds for $h$.
\end{enumerate}
There is no circularity here, because $q_1$ and $\varrho_1'$ do not depend on $C_E$ (the constant
$C_E$ enters the conditions of Lemma~\ref{lem:local} only through the threshold for $h$), and
$C_1$ and the straightening constants (Section~\ref{ss:straightening}) do not depend on
$\kappa$.

In Sections~\ref{sec:stokes}--\ref{sec:determinant} we consider a finite family of backgrounds
simultaneously. For it we take the constants common to all backgrounds and choose them in the
same order: a common bound $R_0$ and the smallest $\sigma$; a common $q_1$ (it is absolute); then
$\varrho_1'$ not smaller than all the $\varrho_1$; then $t_*$ and $\varrho_2$; then the largest of
the $\kappa_1$, $C_E$, $C_1$, $\kappa$ and the smallest of the thresholds for the amplitudes and
for $h$. All conditions of the lemmas of Section~\ref{sec:manifolds} are then preserved, since
they are lower bounds on $\varrho_1'$, $\kappa$ and upper bounds on $q_1$, $|\eps|$, $h$.

\section{Outline of the proof}\label{ss:scheme}

We write the manifolds as graphs over the unperturbed separatrix $q=q^0(u)=4\arctan e^u$. Let
$T^u(u,\alpha)$, $T^s(u,\alpha)$ be their generating functions, $\alpha=(x,\theta)$, and
$\cD=T^u-T^s$. The section $q=\pi$ corresponds to $u=0$, the difference of the dimensionless
actions equals $\partial_\alpha\cD(0,\alpha)$, and therefore
\[
F=\det\partial_\alpha^2\cD(0,0)
\]
(Section~\ref{sec:generating}). Theorem~\ref{thm:main} follows from formula~\eqref{eq:Qprime},
which Proposition~\ref{thm:Qprime} states uniformly in the range~\eqref{eq:tau-range}. We derived
Theorem~\ref{thm:main} from it in Section~\ref{ss:thmQ}.

\paragraph{Structure of the determinant} The matrix $\mathbf D(\tau)=\partial_\alpha^2\cD(0,0)$
is analytic in $\tau$ in a disc independent of $h$. Let
$\mathbf D(\tau)=\mathbf D_0+\tau\mathbf H+\tau^2\mathbf K_2+\dots$ be its Taylor series. For
$\tau=0$ the problem does not depend on $x$, hence $\mathbf D_0=\diag(0,d_B)$, where $d_B$ is
the splitting of the periodic problem with background $\mu M$. The background has period $\pi$
in $\theta$, and therefore $d_B=O(h^{-1}e^{-\pi/h})$, with an exponent twice as large as
in~\eqref{eq:Qprime}. This is exactly where we use that the strong term contains only even
harmonics. Consequently,
\[
F=\tau\,d_B\mathbf H_{11}+\tau^2\bigl(\det\mathbf H+d_B(\mathbf K_2)_{11}\bigr)+O(\tau^3),
\]
and the leading term turns out to be $\tau^2\det\mathbf H$. The matrix $\mathbf H$ is the second
derivative in the angles of the first variation of $\cD$ with respect to the weak amplitude,
that is, of a Melnikov-type integral taken along the trajectories of the problem with the
background. We split $f=f_E+f_O$ into the parts with even and odd harmonics in $\theta$, and
accordingly $\mathbf H=\mathbf H_E+\mathbf H_O$. Then
\[
\mathbf H_E=\diag(4A_f,0)+O(h),\qquad |\mathbf H_O|\le CL_h,\qquad
(\mathbf H_O)_{22}=-2L_h\bigl(G_{M,p_f}(\mu)+o(1)\bigr),
\]
while $\det\mathbf H_E$ is exponentially small compared with $L_h$, because the problem with the
perturbation $f_E$ again has period $\pi$ in $\theta$. Hence
$\det\mathbf H=\bigl(-8A_fG_{M,p_f}(\mu)+o(1)\bigr)L_h$. The entries of the matrices
$\mathbf D$ and $\mathbf H_E$ are not small. Only the determinant is small, as pointed out in
the problem statement.

\paragraph{The fast entry} Up to $O(hL_h)$, the quantity $(\mathbf H_O)_{22}$ equals
$\partial_\theta^2R_{p_f}(0,0)$, where $R_p$ is the first variation of the difference of the
generating functions of the periodic problem with background $\mu M$ with respect to adding the
profile $p$ to the background. We need two facts, namely the bound $R_p=O(L_h)$, uniform on a
complex disc of background amplitudes, and the limit of $R_p/L_h$ as $h\to0$.

We obtain the bound as follows. The difference $\cD$ satisfies the linear transport equation
$\bigl[(1+E)\partial_u+\omega\cdot\partial_\alpha\bigr]\cD=0$, and the first variation satisfies
the same equation with an exponentially small right-hand side. We consider the vertical segment
$u=iy$, $|y|\le\pi/2-Ch$, where $C$ is a large constant. This segment almost reaches the
singularities $u=\pm i\pi/2$ of the separatrix. On it we straighten the equation by a change of
the angles (Section~\ref{sec:straightening}). In the new angles the Fourier coefficient with
index $k=(k_1,k_2)$ is transported along the segment with the factor
$\exp\bigl[k\cdot\omega\,(y+O(h))\bigr]$, and this is what produces the
factor $e^{-|k_2|\pi/(2h)}$. The straightening requires estimates of the manifolds near the
singularities. For the periodic problem they are given by statement~(E1) of
Section~\ref{ss:external}, while we continue the manifolds of the problem with the weak
perturbation up to a distance of order $h$ from the singularities by comparison with those of
the background problem (Section~\ref{sec:manifolds}).

We express the limit through the Stokes constant. Let $g(\theta)$ be the profile of the fast
perturbation. For the periodic problem with profile $zg$ the normalized fast entry
$U_h(zg)=-\partial_\theta^2\cD(0,0)/(2L_h)$ converges as $h\to0$ to an entire function
$\chi[zg]$ (Section~\ref{sec:stokes}). Here we use statement~(E2) on the existence of the
inner limit. The function
$G_h(z)=-\partial_\theta^2R_p(0,0)/(2L_h)$ equals the derivative $\partial_rU_h(zM+rp)|_{r=0}$.
The argument below does not use the analyticity of $\chi$ in the coefficients of the profile.
For each $h$ the Taylor coefficient $[z^k]\,U_h\bigl(z(M+tp)\bigr)$ is a polynomial
in $t$ of degree at most $k$. It converges for every real $t$, hence its coefficient of $t$,
which equals $[z^{k-1}]G_h$, converges as well. Together with the uniform boundedness of $G_h$
and Vitali's theorem, this yields the convergence $G_h\to G_{M,p}$ on the whole plane and
equality~\eqref{eq:response} (Section~\ref{sec:variation}). In this argument we construct the
same periodic problem in three ways, with different decompositions of the perturbation into the
background and the weak part, and we prove in Section~\ref{ss:splitting-indep} that the
constructions agree.

\paragraph{The remainder} Near real $u$ the manifolds are analytic in $\tau$ in a disc
$|\tau|\le e_2$ independent of $h$. We can continue them up to the singularities only for
$|\tau|\le h^{P}$, and on this disc the straightening gives $|F|\le Ch^{-2}e^{-\pi/(2h)}$. By
Cauchy's inequalities $|[\tau^k]F|\le Ch^{-2}e^{-\pi/(2h)}(2/h^{P})^k$, so the terms with
$k\ge3$ are small compared with $\tau^2L_h$ for $|\tau|\le h^{3P+2}$. The linear term
$\tau d_B\mathbf H_{11}$ is small compared with $\tau^2L_h$ for $|\tau|\ge e^{-\pi/(4h)}$. The
exponent $P$ comes from a crude Gronwall estimate in the continuation of the manifolds to the
singularities, and this determines the value of $c_0$ (Section~\ref{sec:determinant}).

\paragraph{The series for $G_{M,p}$} The inner equation is the Hamilton--Jacobi equation that
arises when a neighborhood of the singularity $u=i\pi/2$ is stretched by the factor $1/h$. The
Stokes constant $\chi$ measures the exponentially small difference between two of its
solutions, corresponding to the unstable and stable manifolds. Along the characteristics of this
equation the coordinate satisfies an ordinary differential equation of second order (the inner
limit of the pendulum equation). In variables that straighten the flow of the characteristics,
it takes the form~\eqref{eq:newton}. In Section~\ref{sec:borel} we construct the solutions of
equation~\eqref{eq:newton} and of the linearized equation order by order in the background
amplitude as Laplace integrals of the functions \eqref{eq:Y1}--\eqref{eq:E}, which are analytic
outside the real rays $|s|\ge2$. The two solutions correspond to integrals along rays in the
upper and lower half-planes. Their difference in the first harmonic is determined by a simple
pole at $s=1$, and the residue there equals the sum in~\eqref{eq:G}. A uniqueness lemma shows
that the constructed functions coincide with the Taylor coefficients of the actual solutions.

\paragraph{The sign} The kernel in~\eqref{eq:Qj} is a perfect square, and the other operations
in \eqref{eq:Y1}--\eqref{eq:E} preserve nonnegativity on $[0,2)$, which gives item~(3) of
Theorem~\ref{thm:series}. For the background $\cos4\theta$ only zero modes of even orders enter
the series~\eqref{eq:G}, so $G$ is even and the bound $G\ge2\pi$ holds for both signs of $\mu$.
For the background $\cos2\theta$ with $\mu<0$ we have $G_B(\mu)/2\pi=E_0(1)-E_2(1)$, the
difference of two positive quantities built from the positive background $|\mu|\cos2\theta$.
We bound the zero mode from below by comparison with an explicitly solvable scalar equation,
which gives $E_0(1)\ge\sinh\lambda/\lambda$. We bound the second mode from above as follows.
After the factor $e^{\lambda s}$ is extracted, the system of equations for the remainders is
majorized by a system of inequalities that does not contain the amplitude, and a supersolution
of the latter is given explicitly by thirteen polynomials and verified in exact arithmetic
(Section~\ref{sec:negative}).

The dependencies between the steps of the proof are shown in Table~\ref{tab:deps}.

\begin{table}[ht]
\caption{Dependencies between the steps of the proof. (E1) and (E2) are
Propositions~\ref{prop:E1} and~\ref{prop:E2} on the periodic problem.}\label{tab:deps}
\centering\small
\renewcommand{\arraystretch}{1.15}
\begin{tabular}{@{}>{\raggedright\arraybackslash}p{0.49\linewidth}%
  >{\raggedright\arraybackslash}p{0.13\linewidth}>{\raggedright\arraybackslash}p{0.30\linewidth}@{}}
\toprule
Step & Section & Uses \\
\midrule
Geometry of the background trajectories near the singularity & \ref{ss:geometry} & (E1) \\
Local manifolds, agreement with the graph from (E1), reversibility & \ref{ss:local} & (E1) \\
Continuation of the perturbed manifolds to the singularities & \ref{ss:continuation} & \ref{ss:geometry}, \ref{ss:local} \\
Generating functions, transport equation, determinant & \ref{ss:genfun}, \ref{ss:detdef} & \ref{ss:continuation} \\
Independence of the decomposition of the perturbation into background and weak part & \ref{ss:splitting-indep} & \ref{sec:manifolds}, \ref{ss:genfun} \\
Straightening on the vertical segment & \ref{ss:norms}, \ref{ss:straightening} & self-contained \\
Bound for the determinant on the disc of weak amplitudes & \ref{ss:polydisc} & \ref{sec:generating}, \ref{ss:straightening} \\
Structure of the difference in the periodic problem & \ref{ss:difference} & (E1), \ref{ss:straightening} \\
Stokes constant as a limit & \ref{ss:stokes} & (E1), (E2), \ref{ss:difference} \\
Bound for the first variation on a complex disc of background amplitudes & \ref{ss:variation-bound} & \ref{ss:geometry}, \ref{ss:straightening}, \ref{ss:difference} \\
Limit of the first variation, the function $G_{M,p}$ & \ref{ss:variation-limit} & \ref{ss:continuation}, \ref{ss:splitting-indep}, \ref{ss:stokes}, \ref{ss:variation-bound} \\
Asymptotics of the determinant (Propositions~\ref{thm:Qprime}, \ref{thm:Q} and Theorem~\ref{thm:main}) & \ref{sec:determinant} & \ref{sec:generating}, \ref{ss:polydisc}, \ref{ss:variation-limit} \\
Series for $G_{M,p}$ (Theorem~\ref{thm:series}, items (1), (2)) & \ref{ss:inner}--\ref{ss:jump} & (E1), (E2), \ref{ss:stokes}, \ref{ss:variation-limit} \\
Nonnegativity (Theorem~\ref{thm:series}, item (3); Theorem~\ref{thm:example}) & \ref{ss:positivity} & explicit integrals \\
Background $\cos2\theta$, $\mu<0$ (Theorem~\ref{thm:cos2}) & \ref{sec:negative} & \ref{sec:borel}, computer check \\
\bottomrule
\end{tabular}
\end{table}
\section{Invariant manifolds near the singularities of the separatrix}\label{sec:manifolds}

In this section we construct the unstable and stable manifolds of equation~\eqref{eq:A} and
estimate their deviation from the manifolds of the periodic problem with the same background.
First we study the label trajectories of the periodic problem near the singularities $u=\pm ia$
(Lemma~\ref{lem:geometry}). Then we construct the manifolds near the torus
(Lemma~\ref{lem:local}) and continue them along these trajectories
(Lemma~\ref{lem:continuation}). Of the two statements of Section~\ref{ss:external}, we use only
(E1) in this section.

\subsection{Geometry of the background trajectories}\label{ss:geometry}

Let $g$ be the profile of the periodic problem and $E^u$ the function from (E1). Fix a number
$t_*\ge1$ (it will be chosen in Section~\ref{ss:local}) and take $\varrho_2=2t_*+1$ in (E1). For a
label $u_0$ with $|\Im u_0|<a$ and $t\in[-t_*,0]$ set
\begin{equation}\label{eq:d}
d(t)=\delta(u_0)+|t| .
\end{equation}
By~\eqref{eq:label-motion}, the label of a point moving on the unstable manifold of the periodic
problem satisfies the equation $U'=1+E^u(U,\theta+t/h)$. The following lemma describes such
trajectories issued backward in time from points of the set $\cU^{1/8}_\kappa$.

\begin{lemma}\label{lem:geometry}
There exist constants $\kappa_*$, $h_2$, $C_1$ depending on $R_0$, $g$, $t_*$, and absolute
constants $C_\psi$, $C_S$ such that for $\kappa\ge\kappa_*$, $h<h_2$, $|z|<R_0$,
$\theta\in\TT_\sigma$ and $u_0\in\cU^{1/8}_\kappa$ the solution of the problem
\begin{equation}\label{eq:label-ode}
U'=1+E^u(U,\theta+t/h),\qquad U(0)=u_0,
\end{equation}
exists on $[-t_*,0]$, lies in $D^{\rm out}$ and has the following properties.
\begin{enumerate}[label=\textup{\arabic*.}]
\item $|U(t)-u_0-t|\le C_1|z|h\log\bigl(1+|t|/\delta(u_0)\bigr)\le\min\{d(t)/8,\ 1/10\}$.
\item The distance from $U(t)$ to the nearer of the singularities $\pm ia$ is at least
$0.49\,d(t)$, and to the farther one at least $\tfrac5{24}\,d(t)$.
\item $|E^u(U(t),\theta+t/h)|\le C_1|z|h/d(t)$.
\item If $d(t)\le1$, then $|\psi_0(U(t))|\le C_\psi/d(t)$ and
$|\sin q^0(U(t))|+|\cos q^0(U(t))|\le C_S/d(t)^2$, and one can take $C_\psi=6$, $C_S=46$. If
$d(t)>1$, the same estimates hold with constants depending on $t_*$.
\end{enumerate}
\end{lemma}

\begin{proof}
We may assume that $\Im u_0\ge0$ (otherwise all sets are reflected in the real axis). Let
$\delta=\delta(u_0)$ and
\[
u=u_0+t+\xi,\qquad t\in[-t_*,0],\quad|\xi|\le d(t)/8 .
\]

\emph{Step 1: the point $u$ lies in $D^{\rm out}$.} Since $|\Re u_0|\le\delta/8$, we have
$|\Im u|\le a-\delta+|\xi|$ and $-\Re u\ge|t|-\delta/8-|\xi|$. Therefore the inequality
$|\Im u|<a-\kappa_1h-\Re u$ follows from $2|\xi|<\tfrac78\delta+|t|-\kappa_1h$. The latter
holds if $\kappa\ge2\kappa_1$, since then $\kappa_1h\le\delta/2$, while
$2|\xi|\le(\delta+|t|)/4$. Moreover,
$\Re u\ge-\tfrac\delta8-t_*-\tfrac{\delta+t_*}8\ge-\tfrac a4-\tfrac98t_*>-\varrho_2$.

\emph{Step 2: distances to the singularities.} For the nearer singularity $ia$
\[
|u_0+t-ia|=\sqrt{(\Re u_0+t)^2+\delta^2}\ \ge\ \frac{\tfrac78\delta+|t|}{\sqrt2}\ \ge\ \frac7{8\sqrt2}\,d(t),
\]
whence $|u-ia|\ge\bigl(\tfrac7{8\sqrt2}-\tfrac18\bigr)d(t)\ge0.49\,d(t)$. For the farther
singularity $|u_0+t+ia|\ge\sqrt{a^2+(|t|-\delta/8)_+^2}$. If $|t|\le2a$, then $d(t)\le3a$ and
$|u+ia|\ge a-\tfrac{d(t)}8\ge\tfrac58a\ge\tfrac5{24}d(t)$. If $|t|>2a$, then
$|t|-\tfrac\delta8\ge\tfrac{15}{16}|t|$, $d(t)\le\tfrac32|t|$ and
$|u+ia|\ge\tfrac{15}{16}|t|-\tfrac{d(t)}8\ge\tfrac34|t|\ge\tfrac12d(t)$. Thus
\begin{equation}\label{eq:pole-distances}
|u-ia|\ge0.49\,d(t),\qquad|u+ia|\ge\tfrac5{24}\,d(t),
\end{equation}
and the sum of these distances is at least $2a$.

\emph{Step 3: estimate of $E^u$.} The function $\Pi(u)=\cosh u/(u^2+a^2)$ is entire, and its zeros
$\pm3ia,\pm5ia,\dots$ lie outside the compact set
$\Sigma=\{|\Re u|\le\varrho_2,\ |\Im u|\le a+\tfrac1{10}\}$. Hence $c_\Pi\le|\Pi(u)|\le C_\Pi$
on $\Sigma$, with positive constants depending on $t_*$. By~\eqref{eq:pole-distances} the
smaller of the two distances is at least $\tfrac5{24}d(t)$, and the larger one is at least $a$.
Consequently, for $u\in D^{\rm out}\cap\Sigma$ and any $\vartheta\in\TT_\sigma$ the second
estimate~\eqref{eq:E1} gives
\begin{equation}\label{eq:E-bound-traj}
|E^u(u,\vartheta)|\le C_E|z|h\,\frac{|\Pi(u)|^2}{|u-ia|\,|u+ia|}
\le\frac{24\,C_E\,C_\Pi^2}{5a}\cdot\frac{|z|h}{d(t)}=:C_1\,\frac{|z|h}{d(t)} .
\end{equation}

\emph{Step 4: closing the argument.} Let $t_1$ be the largest number in $[0,t_*]$ such that the
solution of~\eqref{eq:label-ode} exists on $[-t_1,0]$ and satisfies
$|U(t)-u_0-t|\le\min\{d(t)/8,\,1/10\}$ on this interval. Then for $-t_1\le t\le0$ the point
$U(t)$ has the form $u_0+t+\xi$ with $|\xi|\le d(t)/8$, lies in $D^{\rm out}\cap\Sigma$, and
integration of~\eqref{eq:E-bound-traj} gives
\[
|U(t)-u_0-t|\le C_1|z|h\int_0^{|t|}\frac{ds}{\delta+s}=C_1|z|h\log\Bigl(1+\frac{|t|}\delta\Bigr)
\le\frac{C_1R_0}{\kappa}\,|t| .
\]
If $\kappa\ge16C_1R_0$ and $C_1R_0h\log\bigl(1+t_*/(\kappa h)\bigr)<1/20$, the right-hand side is
strictly smaller than both bounds $d(t)/8$ and $1/10$, and by continuity $t_1=t_*$. These
conditions and the condition from Step~1 are satisfied by $\kappa_*=\max\{2\kappa_1,16C_1R_0\}$ and
a sufficiently small $h_2$. This proves items 1--3.

\emph{Step 5: item 4.} Let $d(t)\le1$ and $\varsigma=U(t)-ia$. By Step~2,
$|\varsigma|\ge0.49\,d(t)$.
On the other hand, $|u_0+t-ia|^2\le(|t|+\delta/8)^2+\delta^2\le1.02\,d(t)^2$, so
$|\varsigma|\le\bigl(1.01+\tfrac18\bigr)d(t)\le1.14\,d(t)\le1.25$. Next, $\cosh U=i\sinh\varsigma$
and $\sinh U=i\cosh\varsigma$. The product
$\sinh\varsigma/\varsigma=\prod_{k\ge1}\bigl(1+\varsigma^2/\pi^2k^2\bigr)$ implies
$|\sinh\varsigma|\ge\sin|\varsigma|\ge0.75\,|\varsigma|$ for $|\varsigma|\le1.25$. Hence
$|\cosh U|\ge0.36\,d(t)$, and
\begin{gather*}
|\psi_0(U)|\le\frac{5.6}{d(t)},\qquad
|\cos q^0(U)|=\Bigl|1-\frac2{\cosh^2U}\Bigr|\le1+\frac{15.5}{d(t)^2},\\
|\sin q^0(U)|=\frac{2|\cosh\varsigma|}{|\cosh U|^2}\le\frac{29.5}{d(t)^2} .
\end{gather*}
Therefore one can take $C_\psi=6$ and $C_S=46$. For $d(t)>1$ we use
$|\cosh U|=|\Pi(U)|\,|U-ia|\,|U+ia|\ge c_\Pi\cdot0.49\cdot\tfrac5{24}$.
\end{proof}

\subsection{Local manifolds and reversibility}\label{ss:local}

We now consider equation~\eqref{eq:A} with background $zM$ and weak perturbation $V$. Let $P_M$
be the periodic primitive of the background with zero mean (Section~\ref{ss:separatrix}),
\begin{equation}\label{eq:varpi}
\varpi(t)=hz\,P_M(\theta+t/h),\qquad v=p-\varpi\sin q .
\end{equation}
In the variables $(q,v)$ equation~\eqref{eq:A} takes the form
\begin{equation}\label{eq:qv}
q'=v+\varpi\sin q,\qquad v'=(1+V)\sin q-\varpi\cos q\,\bigl(v+\varpi\sin q\bigr),
\end{equation}
where $V=V(\alpha+\omega t)$. Indeed,
$\frac{d}{dt}[\varpi\sin q]=zM\sin q+\varpi\cos q\cdot q'$, and the term $zM\sin q$ cancels. The
point of this change of variables is that the leading term $\sin q$ enters~\eqref{eq:qv} with the
coefficient $1$ rather than $1+zM$. Set
\[
\varpi_0=hR_0\sup_{\TT_\sigma}|P_M|,\qquad C_V=\sum_i\sup_{\TT^2_\sigma}|f_i| .
\]
Then $|\varpi|\le\varpi_0$ and $|V|\le C_V|\eps|$ for complex phases in the strip and real $t$.

\begin{lemma}\label{lem:local}
There exist $q_1\in(0,\tfrac12]$, $e_1>0$, $h_3>0$ and $C_2$ such that for $h<h_3$,
$|z|\le R_0$, $|\eps|\le e_1$, $\alpha\in\TT^2_\sigma$ and complex $q_*$ with $|q_*|\le q_1$ the
following holds.
\begin{enumerate}[label=\textup{\arabic*.}]
\item Let $w_\pm=(q\pm v)/2$, $\mathbf w=(w_+,w_-)$ and
$\|\mathbf w\|_*=\sup_{t\le0}e^{-t/2}\bigl(|w_+(t)|+|w_-(t)|\bigr)$. System~\eqref{eq:qv} has a
unique solution on $(-\infty,0]$ with $q(0)=q_*$ and $\|\mathbf w\|_*\le2q_1$. For this solution
$\|\mathbf w\|_*\le\tfrac32|q_*|$ and, consequently, $|q(t)|+|v(t)|\le3|q_*|e^{t/2}$.
\item The function $\Gamma(q_*,\alpha;z,\eps):=v(0)$ is analytic in all arguments, and
$|\Gamma(q_*,\alpha;z,\eps)-\Gamma(q_*,\alpha;z,0)|\le C_2\,|\eps|\,|q_*|$.
\item \textup{(Backward invariance.)} If $|q_*|\le q_1/2$, then
$v(t)=\Gamma\bigl(q(t),\alpha+\omega t;z,\eps\bigr)$ for all $t\le0$.
\item \textup{(Agreement with the graph from (E1).)} There exists $\varrho_1'\ge\varrho_1$,
depending only on $q_1$ and $\varrho_1$, such that for $\eps=0$, $\Re u\le-\varrho_1'$,
$\theta\in\TT_\sigma$ and $h<h_3':=\min\{h_3,\,h_1,\,1/(10C_ER_0)\}$ one has $|q^0(u)|\le q_1/4$ and
\begin{equation}\label{eq:coincidence}
\psi_0(u)\bigl(1+E^u(u,\theta)\bigr)=\Gamma\bigl(q^0(u),\theta;z,0\bigr)+hz\,P_M(\theta)\sin q^0(u).
\end{equation}
\end{enumerate}
Here $\eps$ denotes the collection of weak amplitudes $(\eps_i)$. For $\eps=0$ the function
$\Gamma$ depends, of the angles, only on $\theta$.
\end{lemma}

The threshold $h_3'$ depends on the constants $h_1$, $C_E$ from (E1) and is therefore chosen after
$\varrho_2$ (see Section~\ref{ss:constants}).

\begin{proof}
In the variables $w_\pm$ system~\eqref{eq:qv} takes the form
\[
w_\pm'=\pm w_\pm+N_\pm,\qquad N_\pm=\tfrac12(N_q\pm N_v),
\]
where $N_q=\varpi\sin q$ and $N_v=(\sin q-q)+V\sin q-\varpi\cos q\,(v+\varpi\sin q)$. Here
$q=w_++w_-$, $v=w_+-w_-$, so $|q|,|v|\le|w_+|+|w_-|\le|q|+|v|$. A solution on $(-\infty,0]$
with $q(0)=q_*$ and finite norm $\|\mathbf w\|_*$ is the same thing as a fixed point of the
operator $\mathbf w\mapsto\tilde{\mathbf w}$,
\begin{equation}\label{eq:local-operator}
\tilde w_-(t)=\int_{-\infty}^te^{-(t-s)}N_-(s)\,ds,\qquad
\tilde w_+(t)=e^t\bigl(q_*-\tilde w_-(0)\bigr)-\int_t^0e^{t-s}N_+(s)\,ds,
\end{equation}
where $N_\pm$ are computed from $\mathbf w$. If $|N_\pm(s)|\le\bar Ne^{s/2}$, then
$|\tilde w_-(t)|\le\tfrac23\bar Ne^{t/2}$ and
$|\tilde w_+(t)|\le\bigl(|q_*|+\tfrac23\bar N+2\bar N\bigr)e^{t/2}$, so the norm of the image
does not exceed $|q_*|+\tfrac{10}3\bar N$.

In the ball $\|\mathbf w\|_*\le2q_1\le1$ we have $|q|,|v|\le1$, and by the addition formulas the
functions $N_\pm$ are Lipschitz in $\mathbf w$ with respect to $|w_+|+|w_-|$, with constant
$L\le C\,(q_1^2+\varpi_0+C_V|\eps|)$. Here the term $\sin q-q$ contributes $q_1^2$, and the other
terms contain the factor $\varpi$ or $V$. Moreover, $N_\pm=0$ for $\mathbf w=0$. Choose $q_1$,
$h_3$, $e_1$ so that $\tfrac{10}3L\le\tfrac13$. The condition on $h_3$ also includes
$\varpi_0\le\tfrac12$. Then the operator~\eqref{eq:local-operator} maps the ball of radius $2q_1$
into itself, since $|q_*|+\tfrac13\cdot2q_1\le2q_1$, and is a contraction on it with
coefficient~$\tfrac13$. For the fixed point $\|\mathbf w\|_*\le|q_*|+\tfrac13\|\mathbf w\|_*$,
that is, $\|\mathbf w\|_*\le\tfrac32|q_*|$. This is item~1.

The iterates are analytic in $(q_*,\alpha,z,\eps)$ and converge uniformly, hence $\Gamma$ is
analytic. Only the term $V\sin q$ in $N_v$ depends on $\eps$, and $|V\sin q|\le1.2\,C_V|\eps|\,|q|$
for $|q|\le1$. Therefore the difference between the operators with the given $\eps$ and with
$\eps=0$, evaluated at the fixed point, does not exceed
$\tfrac{10}3\cdot0.6\,C_V|\eps|\,\|\mathbf w\|_*=2C_V|\eps|\,\|\mathbf w\|_*$ in norm. Hence
$\|\mathbf w_\eps-\mathbf w_0\|_*\le\tfrac32\cdot2C_V|\eps|\cdot\tfrac32|q_*|$ and
\[
|\Gamma(q_*,\alpha;z,\eps)-\Gamma(q_*,\alpha;z,0)|\le\|\mathbf w_\eps-\mathbf w_0\|_*
\le\tfrac92\,C_V|\eps|\,|q_*| .
\]
This is item~2 with $C_2=\tfrac92C_V$.

For item~3, note that the shifted function $s\mapsto\mathbf w(t+s)$, $s\le0$, is a solution with
phases $\alpha+\omega t$ and value $q(t)$ at zero. Its norm is at most
$e^{t/2}\|\mathbf w\|_*\le\tfrac34q_1$, and $|q(t)|\le\tfrac32|q_*|\le q_1$. The uniqueness in
item~1 gives the claim.

Item~4. Let $\Re u\le-\varrho_1'$, where $\varrho_1'\ge\max\{\varrho_1,2\}$. The solution of the
equation $U'=1+E^u(U,\theta+t/h)$ with $U(0)=u$ stays in $D^{\rm loc}$ for $t\le0$. Indeed,
there $|E^u|\le C_ER_0h<\tfrac1{10}$ by the first estimate~\eqref{eq:E1}, so
$\Re U(t)\le\Re u-0.9\,|t|$. The functions $q(t)=q^0(U(t))$ and
$p(t)=\psi_0(U)\bigl(1+E^u\bigr)$ form a solution of equation~\eqref{eq:A} with $\eps=0$, since
the graph of the gradient of a solution of~\eqref{eq:HJ} is invariant. For $\Re U\le-2$ we have
$|q^0(U)|\le4.1\,e^{\Re U}$ and $|\psi_0(U)|\le4.1\,e^{\Re U}$, hence
\[
|q(t)|+|p(t)|\le12\,e^{\Re U(t)}\le12\,e^{-\varrho_1'}e^{-0.9|t|}.
\]
For $v=p-\varpi\sin q$ we get $|w_+|+|w_-|\le|q|+|v|\le(1+2\varpi_0)\bigl(|q|+|p|\bigr)$. Choose
$\varrho_1'$ so that $12\,e^{-\varrho_1'}(1+2\varpi_0)\le q_1/4$. Since $\varpi_0\le\tfrac12$,
such a $\varrho_1'$ depends only on $q_1$ and $\varrho_1$. Then $q_*=q^0(u)$ satisfies
$|q_*|\le q_1/4$, and the norm of the constructed solution is at most $q_1/4$, so it lies in the
uniqueness ball of item~1. Consequently, $v(0)=\Gamma(q^0(u),\theta;z,0)$, and this equality
coincides with~\eqref{eq:coincidence}, since $\varpi(0)=hz\,P_M(\theta)$.
\end{proof}

\paragraph{The choice of $t_*$} We now fix
\begin{equation}\label{eq:T}
t_*=\tfrac87\,(\varrho_1'+2).
\end{equation}
Then for the trajectories of Lemma~\ref{lem:geometry}
\[
\Re U(-t_*)\le\tfrac{\delta(u_0)}8-t_*+\tfrac{d(-t_*)}8\le\tfrac a4-\tfrac78t_*\le-\varrho_1'-1 .
\]
In particular, $|q^0(U(-t_*))|\le q_1/4$.

\paragraph{Reversibility} Since $K(-\alpha)=K(\alpha)$, the map $q(t)\mapsto2\pi-q(-t)$ takes
solutions of~\eqref{eq:A} with phases $\alpha$ to solutions with phases $-\alpha$. Unstable
half-trajectories, which tend to $q=0$ as $t\to-\infty$, are thereby taken to stable ones, which
tend to $q=2\pi$ as $t\to+\infty$. Since $q^0(-u)=2\pi-q^0(u)$, the stable graph over the label
$u$ is expressed through the unstable one:
\begin{equation}\label{eq:reversibility}
p^s(u,\alpha)=p^u(-u,-\alpha),\qquad E^s(u,\alpha)=E^u(-u,-\alpha).
\end{equation}
We take this equality as the definition of the stable graph everywhere below. All estimates for
it are obtained from the estimates for the unstable graph by reflection, since the sets
$\cU^\gamma_\kappa$ and the segment of the imaginary axis are symmetric with respect to zero.

\subsection{Continuation of the perturbed manifolds}\label{ss:continuation}

In this section the background is $zM$ with $|z|\le R_0$, that is, $g=M$ in (E1), and
$V=\sum_i\eps_if_i$ is the weak perturbation. Set
\begin{equation}\label{eq:Lambda}
\Lambda=2C_S+1,\qquad\vartheta=\frac{R_0\sup_{\TT_\sigma}|P_M|}{\kappa},
\end{equation}
where $C_S$ is the absolute constant of Lemma~\ref{lem:geometry}, and assume that $\kappa$ is so
large that
\begin{equation}\label{eq:kappa-cond}
\vartheta\,\bigl(2C_S+16C_S+7C_S^2\bigr)\le1 .
\end{equation}

For $u_0\in\cU^{1/8}_\kappa$ and $\alpha\in\TT^2_\sigma$ let $U(t)$ be the trajectory of
Lemma~\ref{lem:geometry}. The functions
\begin{equation}\label{eq:background-solution}
q_B(t)=q^0(U(t)),\qquad
v_B(t)=\psi_0(U)\bigl(1+E^u(U,\theta+t/h)\bigr)-\varpi(t)\sin q^0(U)
\end{equation}
form a solution of system~\eqref{eq:qv} with $V=0$. We call it the \emph{background solution},
and we define the \emph{perturbed solution} $(q,v)(t)$, $t\in[-t_*,0]$, as the solution of
system~\eqref{eq:qv} with the initial data
\begin{equation}\label{eq:perturbed-initial}
q(-t_*)=q_B(-t_*),\qquad v(-t_*)=\Gamma\bigl(q_B(-t_*),\alpha-\omega t_*;z,\eps\bigr).
\end{equation}
By the choice of $t_*$ and $\varrho_1'$ we have $|q_B(-t_*)|\le q_1/4$, and by
Lemma~\ref{lem:local}(4) $v_B(-t_*)=\Gamma\bigl(q_B(-t_*),\theta-t_*/h;z,0\bigr)$. For $t\le-t_*$
the perturbed solution is continued by the solution of Lemma~\ref{lem:local} and tends to zero. Thus
the point $(q,p)(0)$, $p=v+\varpi\sin q$, lies on the unstable manifold of the torus $q=p=0$ of
equation~\eqref{eq:A} with phases $\alpha$.

\begin{lemma}\label{lem:continuation}
There exist constants $C_3,C_4,c_3>0$ and $h_4>0$ depending on $R_0$, $M$, $f_i$, $t_*$, $\kappa$
such that for $h<h_4$ the following holds.
\begin{enumerate}[label=\textup{\arabic*.}]
\item If $u_0\in\cU^{1/8}_\kappa$ and $|\eps|\le c_3\,\delta(u_0)^{\Lambda}$, then for
$t\in[-t_*,0]$
\begin{equation}\label{eq:deviation}
|q(t)-q_B(t)|\le C_3\,|\eps|\,\delta(u_0)^{-\Lambda},\qquad
|v(t)-v_B(t)|\le\frac{C_3\,|\eps|\,\delta(u_0)^{-\Lambda}}{d(t)} .
\end{equation}
\item If $|\eps|\le c_3(\kappa h)^{\Lambda+2}$, then on $\cU^{1/32}_{4\kappa}\times\TT^2_\sigma$
there is defined a function $p^u(u,\alpha;z,\eps)$, analytic in all arguments, such that the point
$\bigl(q^0(u),p^u(u,\alpha)\bigr)$ is the endpoint (at $t=0$) of one of the perturbed solutions
constructed above with phases $\alpha$, and
\begin{equation}\label{eq:pu-bound}
|p^u(u,\alpha;z,\eps)-p^u_B(u,\theta;z)|\le C_4\,|\eps|\,\delta(u)^{-\Lambda-2},\qquad
p^u_B=\psi_0\,(1+E^u).
\end{equation}
\item For real labels and labels close to them there is no loss. Namely, there exists $e_2>0$
independent of $h$ such that for $|u|\le a/64$ and $|\eps|\le e_2$ the function $p^u$ is defined and
analytic, and $|p^u-p^u_B|\le C_4\,|\eps|$.
\end{enumerate}
\end{lemma}

\begin{proof}
1. Let $\delta q=q-q_B$, $\delta v=v-v_B$, $S_B=|\sin q_B|+|\cos q_B|$ and $d=d(t)$. As long as
$|\delta q|\le1$, the addition formulas give
\begin{gather*}
|\sin q-\sin q_B|\le2S_B|\delta q|,\qquad|\cos q-\cos q_B|\le2S_B|\delta q|,\\
|\sin q|,\ |\cos q|\le3S_B,\qquad
|\sin q\cos q-\sin q_B\cos q_B|\le7S_B^2|\delta q| .
\end{gather*}
Subtracting the equations~\eqref{eq:qv} for the background solution from those for the perturbed
one, we obtain
\begin{equation}\label{eq:delta-system}
\begin{aligned}
|\delta q'|&\le|\delta v|+2\varpi_0S_B|\delta q|,\\
|\delta v'|&\le2S_B|\delta q|+3C_V|\eps|\,S_B
+\varpi_0\bigl(3S_B|\delta v|+2S_B|v_B|\,|\delta q|\bigr)+7\varpi_0^2S_B^2|\delta q| .
\end{aligned}
\end{equation}
Suppose first that $d(t)\le1$. By Lemma~\ref{lem:geometry}
\[
S_B\le\frac{C_S}{d^2},\qquad|v_B|\le\frac{1.1\,C_\psi+\vartheta C_S}{d}\le\frac8d,\qquad
\frac{\varpi_0}{d}\le\vartheta .
\]
Here we used that $|E^u|\le C_1R_0/\kappa\le\tfrac1{16}$ and
$\varpi_0/d\le\varpi_0/(\kappa h)=\vartheta$. Consider the function $W=|\delta q|+d\,|\delta v|$.
Since $d'=-1$, for its upper right Dini derivative $D^+W$ we obtain from~\eqref{eq:delta-system}
\begin{align*}
D^+W&\le|\delta q'|-|\delta v|+d\,|\delta v'|\\
&\le\frac1d\Bigl[\bigl(2C_S+\vartheta(2C_S+16C_S+7C_S^2\vartheta)\bigr)|\delta q|
+3C_S\vartheta\,d\,|\delta v|\Bigr]+\frac{3C_SC_V|\eps|}{d}\\
&\le\frac{\Lambda}{d}\,\bigl(W+2C_V|\eps|\bigr),
\end{align*}
where the last inequality follows from~\eqref{eq:kappa-cond}. For $d(t)>1$ the same computation
gives $D^+W\le\Lambda'\bigl(W+2C_V|\eps|\bigr)/d$ with a constant $\Lambda'=\Lambda'(t_*,R_0)$. We
integrate from $-t_*$ to $t$. On the part where $d\le1$ we have
$\int\Lambda\,ds/d\le\Lambda\log(1/\delta(u_0))$, therefore
\begin{equation}\label{eq:gronwall}
\begin{gathered}
W(t)+2C_V|\eps|\le C(t_*,R_0)\,\delta(u_0)^{-\Lambda}\bigl(W(-t_*)+2C_V|\eps|\bigr),\\
W(-t_*)=d(-t_*)\,|\delta v(-t_*)|\le(a+t_*)\,C_2\,q_1\,|\eps| .
\end{gathered}
\end{equation}
The bound for $W(-t_*)$ follows from~\eqref{eq:perturbed-initial} and Lemma~\ref{lem:local}(2).
The resulting estimate is exactly~\eqref{eq:deviation}. The condition
$|\eps|\le c_3\delta(u_0)^\Lambda$ with $c_3$ small enough ensures $|\delta q|\le1$ and closes the
argument by continuity.

2. Denote $(q_{\rm end},p_{\rm end})(u_0,\alpha)=(q,\,v+\varpi\sin q)(0)$. These functions are
analytic in $(u_0,\alpha,z,\eps)$, and by item~1
\begin{equation}\label{eq:end-bounds}
|q_{\rm end}-q^0(u_0)|\le C_3\,|\eps|\,\delta_0^{-\Lambda},\qquad
|p_{\rm end}-p^u_B(u_0,\theta)|\le C\,|\eps|\,\delta_0^{-\Lambda-1},\qquad\delta_0=\delta(u_0).
\end{equation}
The second estimate takes into account the term $\varpi(0)(\sin q-\sin q_B)(0)$, which for
$\delta_0\le1$ does not exceed $2C_S\vartheta\,C_3\,|\eps|\,\delta_0^{-\Lambda-1}$.

\emph{The label of the endpoint.} In the disc $|\xi|\le\delta_0/4$ we have
$q^0(u_0+\xi)-q^0(u_0)=\psi_0(u_0)\,\xi+R(\xi)$, where $|R(\xi)|\le C|\xi|^2/\delta_0^2$, since
$(q^0)''=\sin q^0$. Moreover, $|\psi_0(u_0)|\ge2/\cosh(a/8)>1.9$. If
$C_3\,|\eps|\,\delta_0^{-\Lambda}\le c\,\delta_0^2$ with $c$ small enough, then the equation
$q^0(u_0+\xi)=q_{\rm end}$ has a unique solution $\xi(u_0,\alpha)$ in this disc. This solution is
analytic and satisfies $|\xi|\le C_3\,|\eps|\,\delta_0^{-\Lambda}$. Existence is given by the
contraction $\xi\mapsto\bigl(q_{\rm end}-q^0(u_0)-R(\xi)\bigr)/\psi_0(u_0)$ in the disc of radius
$C_3\,|\eps|\,\delta_0^{-\Lambda}$. Uniqueness in the whole disc $|\xi|\le\delta_0/4$ follows from
the injectivity of $q^0$ in the strip $|\Im u|<a$. Indeed, there $\Re e^u>0$, and the function
$w\mapsto4\arctan w$ is injective in the right half-plane.

\emph{Inversion.} The set
$\cU^{1/16}_{2\kappa}=\{16|\Re u|+|\Im u|\le a,\ |\Im u|\le a-2\kappa h\}$ is convex, and the
disc of radius $\delta_0/32$ centered at its point $u_0$ lies in $\cU^{1/8}_\kappa$. By Cauchy's
inequality $|\partial_{u_0}\xi|\le C\,|\eps|\,\delta_0^{-\Lambda-1}\le\tfrac12$. For
$u\in\cU^{1/32}_{4\kappa}$ the disc of radius $\delta(u)/64$ centered at $u$ lies in
$\cU^{1/16}_{2\kappa}$, and the map $u_0\mapsto u-\xi(u_0,\alpha)$ is a contraction of this disc
into itself. The fixed point $u_0(u,\alpha)$ of this map is analytic,
$|u_0-u|\le C\,|\eps|\,\delta(u)^{-\Lambda}$, and $q_{\rm end}(u_0,\alpha)=q^0(u)$ for it. Set
$p^u(u,\alpha)=p_{\rm end}(u_0(u,\alpha),\alpha)$. Then
\[
|p^u(u,\alpha)-p^u_B(u,\theta)|\le|p_{\rm end}(u_0)-p^u_B(u_0)|+\sup|\partial_up^u_B|\cdot|u_0-u|
\le C_4\,|\eps|\,\delta(u)^{-\Lambda-2},
\]
since $|p^u_B|\le1.1\,C_\psi/\delta$ and, by Cauchy's inequality, $|\partial_up^u_B|\le C/\delta^2$.
All the smallness conditions used follow from $|\eps|\le c_3(\kappa h)^{\Lambda+2}$ and
$\delta\ge\kappa h$.

3. For $|u_0|\le a/32$ we have $\delta_0\ge a/2$, and all the estimates of items 1--2 hold with
constants independent of $h$.
\end{proof}

The stable graph $p^s(u,\alpha)=p^u(-u,-\alpha)$ is defined on the same set and satisfies
estimate~\eqref{eq:pu-bound} with $p^s_B$ in place of $p^u_B$.
\section{Generating functions and the splitting determinant}\label{sec:generating}

In this section we describe the manifolds constructed in Section~\ref{sec:manifolds} by
generating functions. For the difference of the generating functions we derive a linear
transport equation, on which the exponentially small estimates of
Section~\ref{sec:straightening} are based. Then we write the splitting determinant as the
Hessian determinant of this difference and prove that the construction does not depend on how
the function $K-1$ is decomposed into the background and the weak perturbation.

\subsection{Generating functions}\label{ss:genfun}

Set $\mathcal N_\kappa=\cU^{1/32}_{4\kappa}$, and let $|\eps|\le c_3(\kappa h)^{\Lambda+2}$. For
$(u,\alpha)\in\mathcal N_\kappa\times\TT^2_\sigma$ let $(q,p)(t)$, $t\le0$, be the perturbed
solution with endpoint $\bigl(q^0(u),p^u(u,\alpha)\bigr)$ and phases $\alpha$ constructed in
Lemma~\ref{lem:continuation}. Here $p=v+\varpi\sin q$, and for $t\le-t_*$ the solution is
continued by Lemma~\ref{lem:local}. Set
\begin{equation}\label{eq:Tu-def}
\begin{gathered}
T^u(u,\alpha)=\int_{-\infty}^0\Bigl[\tfrac12\,p(t)^2+K(\alpha+\omega t)\bigl(1-\cos q(t)\bigr)\Bigr]dt,\\
T^s(u,\alpha)=-T^u(-u,-\alpha),
\end{gathered}
\end{equation}
and
\begin{equation}\label{eq:D-def}
\cD=T^u-T^s .
\end{equation}
The function $T^u$ is the action along the half-trajectory, and it also depends on the
amplitudes $z$ and $\eps=(\eps_i)$.

\begin{lemma}\label{lem:genfun}
Under the above assumptions the following holds.
\begin{enumerate}[label=\textup{\arabic*.}]
\item The function $T^u$ is analytic in $(u,\alpha,z,\eps)$, and $|T^u|\le C_5/\delta(u)$. In
particular, $|\cD|\le2C_5/(4\kappa h)$ on $\mathcal N_\kappa\times\TT^2_\sigma$.
\item $\partial_uT^u=\psi_0\,p^u$ and
\begin{equation}\label{eq:first-variation}
\begin{gathered}
\partial_{\alpha_i}T^u=\int_{-\infty}^0\partial_{\alpha_i}K(\alpha+\omega t)\bigl(1-\cos q(t)\bigr)\,dt,\\
\partial_{\eps_i}T^u=\int_{-\infty}^0f_i(\alpha+\omega t)\bigl(1-\cos q(t)\bigr)\,dt .
\end{gathered}
\end{equation}
The second formula is the first variation with respect to the amplitude of the weak term
$\eps_if_i$.
\item The function $T^u$ satisfies the Hamilton--Jacobi equation
\begin{equation}\label{eq:HJ-qp}
\omega\cdot\partial_\alpha T^u+\frac{(\partial_uT^u)^2}{2\psi_0^2}-K(\alpha)\,\frac{\psi_0^2}2=0 .
\end{equation}
The function $T^s$ satisfies the same equation, and $\partial_uT^s=\psi_0\,p^s$.
\item The difference $\cD$ satisfies the transport equation
\begin{equation}\label{eq:transport}
\bigl[(1+E)\,\partial_u+\omega\cdot\partial_\alpha\bigr]\cD=0,\qquad 1+E=\frac{p^u+p^s}{2\psi_0}.
\end{equation}
\item $\cD(0,\alpha)=T^u(0,\alpha)+T^u(0,-\alpha)$, and this function is even in $\alpha$. The
vector $\partial_\alpha\cD(0,\alpha)$ is the difference of the dimensionless actions of the
unstable and stable manifolds on the section $q=\pi$, that is,
$\partial_\alpha\cD(0,\alpha)=\tilde A^u(\alpha)-\tilde A^s(\alpha)$.
\end{enumerate}
\end{lemma}

\begin{proof}
1. Let $u_0=u_0(u,\alpha)$ be the label of the perturbed solution and $d(t)=\delta(u_0)+|t|$. By
the proof of Lemma~\ref{lem:continuation}, $\delta(u_0)\ge\tfrac{63}{64}\delta(u)$. On
$(-\infty,-t_*]$, by Lemma~\ref{lem:local}, $|q|+|p|\le Ce^{(t+t_*)/2}$, and the integrand
in~\eqref{eq:Tu-def} does not exceed $Ce^{t+t_*}$. On $[-t_*,0]$, by Lemmas~\ref{lem:geometry}
and~\ref{lem:continuation}, $|p|\le C/d(t)$ and $|1-\cos q|\le1+3S_B\le C/d(t)^2$, so the
integral over this interval does not exceed $C/\delta(u_0)$. Analyticity follows from the
uniform convergence of the integral.

2. Let $\partial$ denote the derivative with respect to one of the parameters $u$, $\alpha_i$,
$\eps_i$. The solution $(q,p)(t)$ is analytic in these parameters, and the bound
$|q(t)|+|p(t)|\le Ce^{(t+t_*)/2}$ for $t\le-t_*$ is uniform in them. By Cauchy's inequality the
same bound holds for $\partial q$ and $\partial p$. Therefore we may differentiate under the
integral sign. Integrating by parts and using the equations $q'=p$, $p'=K\sin q$, we obtain
\begin{equation}\label{eq:variation-formula}
\partial T^u=\bigl[p\,\partial q\bigr]_{t=0}
+\int_{-\infty}^0\partial K\,(1-\cos q)\,dt ,
\end{equation}
where $\partial K$ is the derivative of the function $K(\alpha+\omega t)$ with respect to the
parameter entering it explicitly. For the derivative in $u$ the endpoint of the trajectory is
$q^0(u)$, and $K$ does not depend on $u$, hence $\partial_uT^u=p^u\psi_0$. For the derivatives
in $\alpha_i$ and $\eps_i$ the endpoint is fixed, and only the integral remains.

3. Let $\mathcal H(t)=\tfrac12p(t)^2+K(\alpha+\omega t)\bigl(\cos q(t)-1\bigr)$ be the energy
along the solution. Since $q'=p$ and $p'=K\sin q$,
\[
\mathcal H'(t)=(\omega\cdot\partial_\alpha K)(\alpha+\omega t)\,\bigl(\cos q(t)-1\bigr),
\]
and $\mathcal H(t)\to0$ as $t\to-\infty$. Integrating and using item~2 and the identity
$1-\cos q^0=\psi_0^2/2$, we obtain
\[
\tfrac12\,(p^u)^2-K(\alpha)\,\tfrac12\,\psi_0^2=\mathcal H(0)
=-\int_{-\infty}^0\omega\cdot\partial_\alpha K\,(1-\cos q)\,dt=-\omega\cdot\partial_\alpha T^u .
\]
The substitution $p^u=\partial_uT^u/\psi_0$ gives~\eqref{eq:HJ-qp}. Further,
by~\eqref{eq:reversibility}
\[
\partial_uT^s(u,\alpha)=(\partial_uT^u)(-u,-\alpha)=\psi_0(u)\,p^u(-u,-\alpha)=\psi_0(u)\,p^s(u,\alpha),
\]
and the equation for $T^s$ is obtained from the equation for $T^u$ by the substitution
$(u,\alpha)\mapsto(-u,-\alpha)$, since $K$ and $\psi_0$ are even.

4. Subtract the equations of item~3 for $T^u$ and $T^s$ and use the identity
$(\partial_uT^u)^2-(\partial_uT^s)^2=\psi_0\,(p^u+p^s)\,\partial_u\cD$.

5. The first claim follows from the definition of $T^s$. By item~2,
$\partial_{\alpha_i}T^u(0,\alpha)$ is the integral~\eqref{eq:Au-dimless}, that is,
$\tilde A^u_i(\alpha)$. It remains to check that $\partial_\alpha T^s(0,\alpha)=\tilde A^s(\alpha)$.
We have $\partial_\alpha T^s(0,\alpha)=(\partial_\alpha T^u)(0,-\alpha)=\tilde A^u(-\alpha)$. On
the other hand, the stable trajectory with phases $\alpha$ at $t=0$ has the form
$q^s(t)=2\pi-q(-t)$, where $q$ is the unstable trajectory with phases $-\alpha$
(Section~\ref{ss:local}), and
\begin{align*}
\tilde A^s_i(\alpha)&=-\int_0^{\infty}(\partial_{\alpha_i}K)(\alpha+\omega t)\bigl(1-\cos q(-t)\bigr)\,dt\\
&=\int_{-\infty}^0(\partial_{\alpha_i}K)(-\alpha+\omega t)\bigl(1-\cos q(t)\bigr)\,dt
=\tilde A^u_i(-\alpha),
\end{align*}
since the derivative of the even function $K$ is odd.
\end{proof}

\begin{remark}\label{rem:periodic-T}
If there is no weak perturbation ($V=0$), the perturbed solution coincides with the background
one, and $p^u=p^u_B=\psi_0(1+E^u)$. Then the derivative $\partial_uT^u$ of the
function~\eqref{eq:Tu-def} coincides with the derivative of the function $T^u$ from (E1), and by
equations~\eqref{eq:HJ} and~\eqref{eq:HJ-qp} the derivatives in $\theta$ coincide as well. Hence
these two functions differ by a constant, and using the same notation for them causes no
confusion, because only derivatives are used throughout.
\end{remark}

\subsection{The splitting determinant}\label{ss:detdef}

\begin{definition}\label{def:F}
The \emph{splitting determinant} of equation~\eqref{eq:A} is
\begin{equation}\label{eq:F-hessian}
F=\det\partial_\alpha^2\cD(0,0).
\end{equation}
\end{definition}

By item~5 of Lemma~\ref{lem:genfun}, for real amplitudes this is the same determinant as
in~\eqref{eq:F-def}, so $\Delta=(Jg)^2F$. Since the function $\cD(0,\alpha)$ is even,
$\partial_\alpha\cD(0,0)=0$, that is, the point of the section with $\alpha=0$ is homoclinic.

For labels close to real ones, the weak perturbation need not be small compared with $h$. By
Lemma~\ref{lem:continuation}(3), the function $p^u$, and with it $T^u$ and $\cD$, are defined and
analytic for
\[
|u|\le a/64,\qquad\alpha\in\TT^2_\sigma,\qquad|z|<R_0,\qquad|\eps|\le e_2,
\]
where $e_2$ does not depend on $h$. In this domain all statements of Lemma~\ref{lem:genfun}
remain valid, and the estimate of item~1 gives $|\cD|\le C$ with a constant independent of $h$.
Consequently, the matrix $\partial_\alpha^2\cD(0,0)$ and the determinant $F$ are analytic in the
amplitudes $(z,\eps)$ for $|z|<R_0$, $|\eps|\le e_2$, and by Cauchy's inequality they are bounded
there uniformly in $h$. For complex amplitudes formula~\eqref{eq:F-hessian} serves as the
definition of $F$.

\subsection{Independence of the decomposition of the perturbation}\label{ss:splitting-indep}

The construction of Sections~\ref{sec:manifolds} and~\ref{ss:genfun} depends on which part of the
function $K-1$ is called the background and which the weak perturbation, since the variable $v$,
the background trajectories and the label $u_0$ depend on this choice. In
Sections~\ref{sec:stokes}--\ref{sec:determinant} we decompose the same function $K$ in different
ways, so we need the following lemma.

\begin{lemma}\label{lem:splitting}
Let
\[
K=1+z^{(1)}g^{(1)}(\theta)+V^{(1)}(\alpha)=1+z^{(2)}g^{(2)}(\theta)+V^{(2)}(\alpha)
\]
be two decompositions of the same function, where $g^{(i)}$ are real even trigonometric
polynomials with zero mean (a zero background is allowed), $|z^{(i)}|\le R_0$, and $V^{(i)}$ are
the weak parts. We write $\eps^{(i)}$ for the collection of the amplitudes in $V^{(i)}$, so that
$|\eps^{(i)}|$ is the sum of their absolute values. Let the constants of
Lemmas~\ref{lem:geometry}--\ref{lem:continuation} be chosen
common to both decompositions, as described in Section~\ref{ss:constants}: common $R_0$,
$\sigma$, $q_1$, $\varpi_0$, $\varrho_1'$, $t_*$, $\kappa$, and $e_1$, $e_2$ the smaller of the
two. There exist $e_3\in(0,e_2]$ and $\tilde h_4>0$ (not larger than the thresholds of
Lemmas~\ref{lem:geometry}--\ref{lem:continuation} for both decompositions) such that for
$|\eps^{(1)}|,|\eps^{(2)}|\le e_3$, $h<\tilde h_4$, $|u|\le a/64$ and $\alpha\in\TT^2_\sigma$ the
graphs of the two constructions coincide:
\begin{equation}\label{eq:splitting-indep}
p^u_{(1)}(u,\alpha)=p^u_{(2)}(u,\alpha).
\end{equation}
Consequently, for $|u|\le a/64$ the functions $T^u(u,\alpha)$ and $\cD(u,\alpha)$ of the two
constructions coincide, and hence so does the determinant $F$.
\end{lemma}

\begin{proof}
Let $\varpi_i(t)=hz^{(i)}P_{g^{(i)}}(\theta+t/h)$ and $v_i=p-\varpi_i\sin q$. Let
$\mathbf w_i=(w_{i,+},w_{i,-})$ be the variables of Lemma~\ref{lem:local} for the $i$th
decomposition, and $\Gamma_i$ the corresponding local graphs. By assumption
$|\varpi_i|\le\varpi_0\le\tfrac12$.

\emph{Step 1: the local graphs coincide in the coordinates $(q,p)$.} For $|q_*|\le q_1/2$ and
$\alpha\in\TT^2_\sigma$ we have
\begin{equation}\label{eq:local-graphs-coincide}
\Gamma_1(q_*,\alpha)+\varpi_1(0)\sin q_*=\Gamma_2(q_*,\alpha)+\varpi_2(0)\sin q_* .
\end{equation}
Indeed, let $(q,v_1)$ be the solution of Lemma~\ref{lem:local}(1) for the first decomposition.
Then $\|\mathbf w_1\|_*\le\tfrac32|q_*|\le\tfrac34q_1$. The pair $(q,p)$ with
$p=v_1+\varpi_1\sin q$ is a solution of equation~\eqref{eq:A}, hence the pair $(q,v_2)$ with
$v_2=v_1+(\varpi_1-\varpi_2)\sin q$ is a solution of system~\eqref{eq:qv} of the second
decomposition with $q(0)=q_*$. Since
\[
|w_{2,\pm}-w_{1,\pm}|\le\varpi_0\,|\sin q|\le1.2\,\varpi_0\,\bigl(|w_{1,+}|+|w_{1,-}|\bigr),
\]
we have $\|\mathbf w_2\|_*\le(1+2.4\,\varpi_0)\cdot\tfrac34q_1\le2q_1$. By the uniqueness in
Lemma~\ref{lem:local}(1), $v_2(0)=\Gamma_2(q_*,\alpha)$, which is exactly
equality~\eqref{eq:local-graphs-coincide}.

\emph{Step 2: the same trajectories.} Let $|u|\le a/64$. Denote by $\gamma$ the perturbed solution
of the first construction with endpoint $\bigl(q^0(u),p^u_{(1)}(u,\alpha)\bigr)$, and by
$u_0=u_0^{(1)}(u,\alpha)$ its label. By the estimate of the label in the proof of
Lemma~\ref{lem:continuation} (the step ``Inversion''), $|u_0-u|\le C\,|\eps^{(1)}|$. Let
$U_i(t;u',\theta)$ be the background trajectory of the $i$th decomposition with $U_i(0)=u'$. At
$t=-t_*$ we have $q(-t_*)=q^0(u_*)$, where $u_*=U_1(-t_*;u_0,\theta)$, and $|q(-t_*)|\le q_1/4$. The
point $\gamma(-t_*)$ lies on the local graph of the first decomposition with phases
$\alpha-\omega t_*$, and by Step~1 also on that of the second. By Lemma~\ref{lem:geometry}(1)
\[
|u_*-(u_0-t_*)|\le\epsilon_h:=C_1R_0h\log(1+2t_*/a).
\]
The function $\Upsilon(u')=U_2(-t_*;u',\theta)$ is analytic in the disc $|u'-u|\le a/16$, because
this disc lies in $\cU^{1/8}_\kappa$ and $\delta(u')\ge a/2$ in it. In this disc
$|\Upsilon(u')-(u'-t_*)|\le\epsilon_h$. By Cauchy's inequality the map
\[
u'\mapsto u_*+t_*-\bigl[\Upsilon(u')-u'+t_*\bigr]
\]
has derivative at most $32\,\epsilon_h/a\le\tfrac12$ in the disc $|u'-u|\le a/32$ and maps this
disc into the disc $|u'-u|\le2\epsilon_h+C\,|\eps^{(1)}|$. Its fixed point $u_0'$ satisfies
$\Upsilon(u_0')=u_*$ and $|u_0'-u|\le2\epsilon_h+C\,|\eps^{(1)}|$.

At $t=-t_*$ the perturbed solution of the second construction with label $u_0'$ and phases
$\alpha$ has $q=q^0(\Upsilon(u_0'))=q(-t_*)$ and lies on the local graph of the second
decomposition, that is, it passes through the point $\gamma(-t_*)$. By the uniqueness of the
solution of the Cauchy problem for equation~\eqref{eq:A}, it coincides with $\gamma$ on
$[-t_*,0]$, and its endpoint is the point $\bigl(q^0(u),p^u_{(1)}(u,\alpha)\bigr)$. Hence
$\xi=u-u_0'$ solves the equation $q^0(u_0'+\xi)=q^{(2)}_{\rm end}(u_0',\alpha)$, and
$|\xi|<\delta(u_0')/4$. By the uniqueness in the step ``The label of the endpoint'' of the proof
of Lemma~\ref{lem:continuation}, $\xi=\xi^{(2)}(u_0',\alpha)$. Consequently, $u_0'$ is a fixed
point of the map $u_0\mapsto u-\xi^{(2)}(u_0,\alpha)$ lying in the disc of radius $\delta(u)/64$
centered at $u$, provided $2\epsilon_h+Ce_3<a/128$, which is a condition on $e_3$ and
$\tilde h_4$. By the uniqueness of the fixed point in the step ``Inversion'',
$u_0'=u_0^{(2)}(u,\alpha)$, whence
\[
p^u_{(2)}(u,\alpha)=p^{(2)}_{\rm end}(u_0',\alpha)=p^u_{(1)}(u,\alpha).
\]

\emph{Step 3: conclusion.} The function $T^u(u,\alpha)$ is the integral of
$\tfrac12p^2+K(1-\cos q)$ along the solution of equation~\eqref{eq:A} on $(-\infty,0]$ with
endpoint
$\bigl(q^0(u),p^u(u,\alpha)\bigr)$. It is determined by the function $K$ and the graph, not by
the decomposition.
\end{proof}

\begin{remark}\label{rem:splitting-domain}
The domains of applicability of different decompositions overlap only for small amplitudes,
since one needs $|\eps^{(i)}|\le e_3$. Below, we use Lemma~\ref{lem:splitting} only to compare
Taylor coefficients in the amplitudes at zero, and this suffices.
\end{remark}
\section{Straightening and exponentially small estimates}\label{sec:straightening}

In this section we consider the transport equation~\eqref{eq:transport} for the difference of the
generating functions on the segment $u=iy$ of the imaginary axis whose endpoints are at a distance
of order $h$ from the singularities $u=\pm ia$. A change of the angles reduces it to an equation
with coefficients independent of the angles (Lemma~\ref{lem:straightening}). After this, each
Fourier coefficient of the difference is transported along the segment with an explicit
exponential factor, which yields exponentially small estimates (Lemma~\ref{lem:polydisc}).
Lemmas~\ref{lem:G-bounds} and~\ref{lem:straightening} do not depend on the other sections.

\subsection{Norms and the inversion of the transport operator}\label{ss:norms}

Let $\kappa_0>0$ and
\[
y_*=a-\kappa_0h,\qquad\delta(y)=a-|y|\quad(|y|\le y_*).
\]
We consider functions $f(y,\alpha)=\sum_{k\in\ZZ^2}f_k(y)e^{ik\cdot\alpha}$ on
$[-y_*,y_*]\times\TT^2_\rho$, continuous in $y$ and analytic in $\alpha$, with the norms
\begin{equation}\label{eq:norms}
\begin{gathered}
\|f\|_\rho=\sum_ke^{\rho|k|}\sup_y|f_k(y)|,\\
\|f\|_{\rho,j}=\sum_ke^{\rho|k|}\sup_y\Bigl(\frac{\delta(y)}h\Bigr)^{j}|f_k(y)|\quad(j=1,2),
\end{gathered}
\end{equation}
where $|k|=|k_1|+|k_2|$. The following properties are obvious:
\begin{itemize}
\item $\sup|f|\le\|f\|_\rho$, $\|fg\|_\rho\le\|f\|_\rho\|g\|_\rho$ and
$\|fg\|_{\rho,i+j}\le\|f\|_{\rho,i}\|g\|_{\rho,j}$, where the index $0$ corresponds to the norm
$\|\cdot\|_\rho$;
\item $\|f\|_\rho\le\kappa_0^{-j}\|f\|_{\rho,j}$, since $\delta(y)\ge\kappa_0h$;
\item if $\sup_{\TT^2_{\rho'}}|f(y,\cdot)|\le A(y)$ for some $\rho'>\rho$, then
$|f_k(y)|\le e^{-\rho'|k|}A(y)$ and $\|f\|_\rho\le c(\rho'-\rho)\sup_yA$, where
$c(r)=\sum_{k\in\ZZ^2}e^{-r|k|}$, and similarly for the weighted norms.
\end{itemize}

On the vertical line $u=iy$ we have $\partial_u=-i\,d/dy$. Set $\lambda_k=k\cdot\omega$ and define
the operator $\cG$ by
\begin{equation}\label{eq:G-operator}
(\cG f)_k(y)=i\int_{y_k}^ye^{\lambda_k(y-v)}f_k(v)\,dv,\qquad
y_k=\begin{cases}y_*,&\lambda_k>0,\\-y_*,&\lambda_k<0,\\0,&\lambda_k=0.\end{cases}
\end{equation}
It is a right inverse of the operator $\partial_u+\omega\cdot\partial_\alpha$, that is,
$-i(\cG f)_k'+i\lambda_k(\cG f)_k=f_k$. We choose the lower limit so that the exponential under
the integral does not exceed one.

\begin{lemma}[estimates of the operator $\cG$]\label{lem:G-bounds}
\begin{enumerate}[label=\textup{(\roman*)}]
\item $\|\cG f\|_\rho\le2y_*\|f\|_\rho$ and $\|\partial_u\cG f\|_\rho\le2\|f\|_\rho$.
\item Let the functions depend only on one angle $\theta$ (then $\lambda_k=k/h$, $k\in\ZZ$), and
let $\kappa_0\ge2$. Then for $j=1,2$
\[
\|\cG f-(\cG f)_0\|_{\rho,j}\le9h\,\|f\|_{\rho,j},\qquad
\|\partial_u\cG f\|_{\rho,j}\le10\,\|f\|_{\rho,j}.
\]
\end{enumerate}
\end{lemma}

\begin{proof}
(i) If $\lambda_k\ne0$, then $|e^{\lambda_k(y-v)}|\le1$ on the interval of integration and
$\int|\lambda_k|e^{-|\lambda_k||y-v|}dv\le1$. Moreover,
$(\partial_u\cG f)_k=f_k-i\lambda_k(\cG f)_k$. If $\lambda_k=0$, then $(\partial_u\cG f)_k=f_k$.

(ii) Let $k\ge1$ and $|f_k(v)|\le N\,(h/\delta(v))^j$. Then
$\tfrac kh\,|(\cG f)_k(y)|\le N\,I(y)$, where
\[
I(y)=\int_y^{y_*}\frac kh\,e^{-k(v-y)/h}\Bigl(\frac h{\delta(v)}\Bigr)^{j}dv .
\]
Suppose first that $y\ge0$, and set $D=\delta(y)$ and $v=y+s$. For $s\le D/2$ the weight is at
most $(2h/D)^j$, and the mass of the kernel is at most one. Values $s\ge D/2$ are possible only if
$D\ge2\kappa_0h$. For them the weight is at most $\kappa_0^{-j}$, and the mass of the kernel is at
most $e^{-D/(2h)}$. Set $r=D/h\ge2\kappa_0\ge4$. The function $r^je^{-r/2}$ decreases for
$r\ge2j$, therefore
\[
\kappa_0^{-j}e^{-D/(2h)}=\Bigl(\frac hD\Bigr)^{j}\Bigl(\frac r{\kappa_0}\Bigr)^{j}e^{-r/2}
\le2^j\Bigl(\frac hD\Bigr)^{j}.
\]
Thus $I(y)\le2^{j+1}(h/\delta(y))^j$ for $y\ge0$. For $y<0$ the part $v\in[y,|y|]$ contributes at
most $(h/\delta(y))^j$, since $\delta(v)\ge\delta(y)$ there, and the rest does not exceed
$I(|y|)$. Consequently, $I(y)\le(1+2^{j+1})(h/\delta(y))^j\le9\,(h/\delta(y))^j$ for all $y$.
Hence
\[
\sup_y\Bigl(\frac{\delta(y)}h\Bigr)^{j}|(\cG f)_k(y)|\le\frac{9h}{k}\,N,\qquad
\sup_y\Bigl(\frac{\delta(y)}h\Bigr)^{j}|(\partial_u\cG f)_k(y)|\le10\,N,
\]
since $(\partial_u\cG f)_k=f_k-i\tfrac kh(\cG f)_k$. The case $k\le-1$ is symmetric, and
$(\partial_u\cG f)_0=f_0$.
\end{proof}

\subsection{The straightening lemma}\label{ss:straightening}

\begin{lemma}[straightening]\label{lem:straightening}
Let $E=E_B(y,\theta)+E_1(y,\alpha)$, where
\begin{equation}\label{eq:straightening-hyp}
\|E_B\|_{\rho,1}\le A_E,\qquad|\langle E_B\rangle_\theta(y)|\le A_0\Bigl(\frac h{\delta(y)}\Bigr)^{2},
\qquad\|E_1\|_\rho\le\epsilon_E\le h^2 .
\end{equation}
Fix $\varsigma\in(0,\rho/4)$ and set $\rho_1=\rho-\varsigma$, $\rho_2=\rho-3\varsigma$. Let
\[
\kappa_0\ge\kappa_0^*:=\max\{2,\ 40A_E,\ 576A_E/\varsigma\},\qquad
h<h_5:=\min\{\tfrac14,\ \varsigma/(128a),\ \nu_0^{-1/2}\}
\]
(the last condition gives $|\nu|\le\nu_0h\le1/h$, that is, $|\omega|\le2/h$). Then the following
holds.
\begin{enumerate}[label=\textup{\arabic*.}]
\item There exist a function $W(y,\alpha)$ of class $C^1$ in $y$ and analytic in
$\alpha\in\TT^2_\rho$, and $H=\partial_uW$, such that
\begin{equation}\label{eq:straightening-eq}
(1+E)(1+H)+\omega\cdot\partial_\alpha W=1 .
\end{equation}
\item $H=H_B+H_1$, where $H_B$ depends only on $(y,\theta)$, $\|H_B\|_{\rho,1}\le20A_E$ and
$\|H_1\|_\rho\le4\epsilon_E$.
\item Let $W_0(y)$ be the zero mode of the function $W(y,\cdot)$. Then $W_0(0)=0$ and
\begin{equation}\label{eq:m-bound}
\begin{gathered}
|W_0(y)|\le A_W\frac{h^2}{\delta(y)}+4a\,\epsilon_E,\qquad
|W(y,\alpha)-W_0(y)|\le18A_E\frac{h^2}{\delta(y)}+4a\,\epsilon_E,\\
A_W=A_0+20A_E^2 .
\end{gathered}
\end{equation}
\item For every $y$ the map
\begin{equation}\label{eq:b-map}
\mathfrak b_y:\ \alpha\mapsto\beta=\alpha-\omega\,(W-W_0)(y,\alpha)
\end{equation}
is injective on $\TT^2_{\rho_1}$, its image contains $\TT^2_{\rho_2}$, the inverse map
$\mathfrak a_y$ is analytic in $\beta$ and of class $C^1$ in $y$, and
$|\mathfrak b_y(\alpha)-\alpha|\le36A_Eh/\delta(y)+8ah$.
\item Let a function $\cD(y,\alpha)$ be of class $C^1$ in $y$, analytic in $\alpha\in\TT^2_\rho$,
and satisfy
\[
(1+E)\,\partial_u\cD+\omega\cdot\partial_\alpha\cD=S .
\]
Then the Fourier coefficients of the function $X(y,\beta)=\cD(y,\mathfrak a_y(\beta))$,
$\beta\in\TT^2_{\rho_2}$, satisfy the equation
\begin{equation}\label{eq:Xk-ode}
X_k'=\lambda_k\bigl(1-iW_0'(y)\bigr)X_k+i\tilde S_k,\qquad
\tilde S(y,\beta)=\frac S{1+E}\bigl(y,\mathfrak a_y(\beta)\bigr),
\end{equation}
and, consequently,
\begin{equation}\label{eq:Xk-transport}
X_k(0)=X_k(y)\,e^{-\lambda_k(y-iW_0(y))}-i\int_0^ye^{-\lambda_k(v-iW_0(v))}\tilde S_k(v)\,dv .
\end{equation}
\end{enumerate}
\end{lemma}

In the applications $E_B$ is the contribution of the background and $E_1$ that of the weak
perturbation. The bound on the mean in~\eqref{eq:straightening-hyp} ensures that the zero mode
$W_0$, which enters the exponent in~\eqref{eq:Xk-transport}, is small.

\begin{proof}
The conditions on $\kappa_0$ and $h$ imply the following five inequalities, which we use below:
\begin{equation}\label{eq:five-inequalities}
\begin{gathered}
\frac{20A_E}{\kappa_0}\le\frac12,\qquad\epsilon_E\le h^2\le\frac1{16},\qquad
\frac{A_E}{\kappa_0}+\epsilon_E\le\frac18,\\
\frac{24A_E}{\kappa_0}+4\epsilon_E\le1,\qquad
\frac{36A_E}{\kappa_0}+8ah\le\frac\varsigma8 .
\end{gathered}
\end{equation}

\emph{The background.} The map $H\mapsto\partial_u\cG\bigl[-E_B(1+H)\bigr]$ acts in the space of
functions of $(y,\theta)$. By Lemma~\ref{lem:G-bounds}(ii) it is Lipschitz in the norm
$\|\cdot\|_{\rho,1}$ with constant $10A_E/\kappa_0\le\tfrac12$ and takes zero to an element of
norm at most $10A_E$. Its fixed point $H_B$ satisfies $\|H_B\|_{\rho,1}\le20A_E$. Set
$W_B=\cG\bigl[-E_B(1+H_B)\bigr]$ and denote by $W_{B,0}$ the zero mode of $W_B$. Then
$\partial_uW_B=H_B$, and for the nonzero modes
\[
\|W_B-W_{B,0}\|_{\rho,1}\le9h\cdot A_E\Bigl(1+\frac{20A_E}{\kappa_0}\Bigr)\le18A_Eh .
\]
For the zero mode $W_{B,0}'=-i\langle E_B(1+H_B)\rangle_\theta$, and pointwise in $y$
\begin{gather*}
|W_{B,0}'(y)|\le|\langle E_B\rangle_\theta|+\sup_\theta|E_B|\,\sup_\theta|H_B|
\le\bigl(A_0+20A_E^2\bigr)\Bigl(\frac h{\delta(y)}\Bigr)^{2},\\
|W_{B,0}(y)|\le A_Wh^2\int_0^{|y|}\frac{dv}{(a-v)^2}\le\frac{A_Wh^2}{\delta(y)} .
\end{gather*}

\emph{The full problem.} We seek $H=H_B+H_1$. The equation $H=\partial_u\cG\bigl[-E(1+H)\bigr]$
is equivalent to the equation
\[
H_1=\partial_u\cG\bigl[-E_BH_1-E_1(1+H_B+H_1)\bigr].
\]
By Lemma~\ref{lem:G-bounds}(i) the right-hand side is a contraction in the norm $\|\cdot\|_\rho$
with constant $2(A_E/\kappa_0+\epsilon_E)\le\tfrac14$. Since
$\|H_B\|_\rho\le20A_E/\kappa_0\le\tfrac12$, we get
$\|H_1\|_\rho\le\tfrac43\cdot2\epsilon_E\bigl(1+\|H_B\|_\rho\bigr)\le4\epsilon_E$. Set
$W=\cG\bigl[-E(1+H)\bigr]$. Then $\partial_uW=H$ and
$(\partial_u+\omega\cdot\partial_\alpha)W=-E(1+H)$, which is equivalent
to~\eqref{eq:straightening-eq}. Further,
\[
\|W-W_B\|_\rho\le2y_*\bigl(\|E_B\|_\rho\|H_1\|_\rho+\epsilon_E(1+\|H\|_\rho)\bigr)
\le2a\,\epsilon_E\Bigl(1+\frac{24A_E}{\kappa_0}+4\epsilon_E\Bigr)\le4a\,\epsilon_E .
\]
This gives items 2 and 3. The equality $W_0(0)=0$ follows from the definition of $\cG$ on the zero
mode.

\emph{Item 4.} Let $g=\omega\,(W-W_0)$. By item~3 and the inequality $|\omega|\le2/h$
\[
\sup_{\TT^2_\rho}|g|\le\frac{36A_Eh}{\delta(y)}+8ah\le\epsilon_\beta:=\frac{36A_E}{\kappa_0}+8ah
\le\frac\varsigma8,
\]
and by Cauchy's inequality $|\partial_\alpha g|\le2\epsilon_\beta/\varsigma\le\tfrac14$ on
$\TT^2_{\rho_1}$. For $\beta\in\TT^2_{\rho_2}$ the map $\alpha\mapsto\beta+g(y,\alpha)$ is a
contraction of the closed ball of radius $2\epsilon_\beta\le\varsigma$ centered at $\beta$ into
itself, and its fixed point is $\mathfrak a_y(\beta)$. To prove injectivity, note that if
$\mathfrak b_y(\alpha)=\mathfrak b_y(\alpha')$ (in the covering space), then
$|\alpha-\alpha'|\le\tfrac14|\alpha-\alpha'|$.

\emph{Item 5.} Let $\mathsf L=(1+E)\,\partial_u+\omega\cdot\partial_\alpha$. Item~1 implies
$\mathsf LW=-E$, and $\mathsf LW_0=-i(1+E)W_0'$. Therefore
\[
\mathsf L\beta=\omega-\omega\,\mathsf L(W-W_0)=\omega\,(1+E)(1-iW_0').
\]
Differentiating the equality $\cD(y,\alpha)=X(y,\mathfrak b_y(\alpha))$, we obtain
\[
S=\mathsf L\cD=(1+E)\bigl[-i\,\partial_yX+(1-iW_0')\,\omega\cdot\partial_\beta X\bigr],
\]
that is, $-i\,\partial_yX+(1-iW_0')\,\omega\cdot\partial_\beta X=\tilde S$. The division is
legitimate, since $|E|\le A_E/\kappa_0+\epsilon_E\le\tfrac12$. Passing to Fourier coefficients in
$\beta$ gives~\eqref{eq:Xk-ode}, and the integrating factor $e^{-\lambda_k(y-iW_0(y))}$ gives
formula~\eqref{eq:Xk-transport}.
\end{proof}

\subsection{Bound for the determinant on the polydisc of weak amplitudes}\label{ss:polydisc}

We return to equation~\eqref{eq:A} with background $zM$ and weak perturbation $V$. Set
\begin{equation}\label{eq:P-def}
\kappa_0=4\kappa,\qquad\rho=\frac\sigma2,\qquad\varsigma=\frac\sigma{12},\qquad P=\Lambda+5 .
\end{equation}
Then the segment $\{iy:\ |y|\le y_*\}$ lies in $\mathcal N_\kappa$, and $\rho_2=\sigma/4$. We
denote by $E^u_B$ and $E^s_B$ the functions \eqref{eq:E-def}, \eqref{eq:reversibility} of the
periodic problem with background $zM$. They depend only on $(u,\theta)$.

\begin{lemma}\label{lem:polydisc}
Let $\kappa$ be large enough, $h<h_6$, $|z|\le R_0$ and $|\eps|\le h^{P}$. Then the function $E$
from~\eqref{eq:transport} on the segment $u=iy$, $|y|\le y_*$, satisfies the assumptions of
Lemma~\ref{lem:straightening} with
\[
E_B=\tfrac12\bigl(E^u_B+E^s_B\bigr),\qquad A_E=c(\sigma/2)\,C_1R_0,\qquad A_0=\tfrac12(C_1R_0)^2,
\]
and the following holds.
\begin{enumerate}[label=\textup{\arabic*.}]
\item For the function $X(y,\beta)=\cD(y,\mathfrak a_y(\beta))$ and all $k\in\ZZ^2$
\begin{equation}\label{eq:Xk-bound}
|X_k(0)|\le\frac{C_6}h\,e^{-\rho_2|k|}\exp\bigl[-|\lambda_k|\,\bigl(a-(\kappa_0+1)h\bigr)\bigr].
\end{equation}
\item $F=\bigl(\det\partial_\alpha\mathfrak b_0(0)\bigr)^2\,\det\partial_\beta^2X(0,\beta_0)$,
where $\beta_0=\mathfrak b_0(0)$. Moreover, $\partial_\alpha\mathfrak b_0=\mathrm{Id}+O(h)$ and
$|\Im\beta_0|\le Ch$.
\item If $K(x,\theta+2\pi/N)\equiv K(x,\theta)$ for some integer $N\ge1$, then
\begin{equation}\label{eq:F-bound}
\begin{gathered}
|\partial_{\beta_1}\partial_{\beta_2}X(0,\beta_0)|+|\partial_{\beta_2}^2X(0,\beta_0)|
\le\frac{C_6}h\,e^{-Na/h},\qquad
|\partial_{\beta_1}^2X(0,\beta_0)|\le\frac{C_6}h,\\
|F|\le\frac{C_6}{h^{2}}\,e^{-Na/h}.
\end{gathered}
\end{equation}
\end{enumerate}
\end{lemma}

\begin{proof}
\emph{The assumptions of Lemma~\ref{lem:straightening}.} By Lemma~\ref{lem:geometry}(3) at $t=0$
and by~\eqref{eq:reversibility}
\[
|E^{u}_B(iy,\theta)|\le\frac{C_1R_0h}{\delta(y)},\qquad
|E^{s}_B(iy,\theta)|\le\frac{C_1R_0h}{\delta(y)}\qquad(\theta\in\TT_\sigma),
\]
whence $\|E_B\|_{\rho,1}\le c(\sigma/2)\,C_1R_0=A_E$. Average equation~\eqref{eq:HJ} with $g=M$
over real $\theta$ for fixed $u$. Since $T^u$ is periodic in $\theta$ and the mean of $M$ is zero,
we obtain the exact identity
\begin{equation}\label{eq:mean-identity}
\bigl\langle(1+E^u_B)^2\bigr\rangle_\theta=1,\qquad\text{that is,}\qquad
\langle E^u_B\rangle_\theta=-\tfrac12\bigl\langle(E^u_B)^2\bigr\rangle_\theta,
\end{equation}
and the same for $E^s_B$. Hence $|\langle E_B\rangle_\theta|\le\tfrac12(C_1R_0)^2(h/\delta)^2$.
Further, $E_1=E-E_B=(p^u-p^u_B+p^s-p^s_B)/(2\psi_0)$, and by~\eqref{eq:pu-bound}
\[
\sup_{\TT^2_\sigma}|E_1|\le C_4\,|\eps|\,(\kappa_0h)^{-\Lambda-2}\le Ch^3,
\]
so $\|E_1\|_\rho\le h^2$ for small $h$. The condition $\kappa_0\ge\kappa_0^*$ is a condition on
$\kappa$, since the constants $A_E$, $A_0$ do not depend on $\kappa$.

1. By Lemma~\ref{lem:genfun}(4), $S=0$ in the transport equation, and
formula~\eqref{eq:Xk-transport} with $y=y_*\operatorname{sgn}\lambda_k$ gives
\[
|X_k(0)|\le|X_k(y)|\,\exp\bigl[-|\lambda_k|\,\bigl(y_*-|W_0(y)|\bigr)\bigr].
\]
Here $|X_k(y)|\le e^{-\rho_2|k|}\sup|\cD|\le e^{-\rho_2|k|}\cdot2C_5/(\kappa_0h)$ by
Lemma~\ref{lem:genfun}(1), and by~\eqref{eq:m-bound} $|W_0|\le A_Wh/\kappa_0+4ah^2\le h$ if
$\kappa_0\ge2A_W$ and $h\le1/(8a)$.

2. By Lemma~\ref{lem:genfun}(5), $\partial_\alpha\cD(0,0)=0$. Since
$\cD(0,\alpha)=X(0,\mathfrak b_0(\alpha))$ and the matrix $\partial_\alpha\mathfrak b_0$ is
invertible, $\partial_\beta X(0,\beta_0)=0$, and in the formula for the second derivative of the
composite function only the term
$(\partial_\alpha\mathfrak b_0)^{\mathsf T}\,\partial_\beta^2X\,(\partial_\alpha\mathfrak b_0)$
remains. The estimates of the map $\mathfrak b_0$ follow from items 3 and 4 of
Lemma~\ref{lem:straightening} at $y=0$ and Cauchy's inequality.

3. All the constructions (the local manifolds, the functions $p^{u}$, $p^s$, $T^{u}$, $T^s$, the
fixed points of Lemma~\ref{lem:straightening}, the maps $\mathfrak b_y$) commute with the shift
$\theta\mapsto\theta+2\pi/N$ by uniqueness. Therefore $X_k=0$ if $k_2$ is not divisible by $N$.
The derivatives with index $2$ involve only the modes with $k_2\ne0$, that is, with
$|k_2|\ge N$. For them $|\lambda_k|\ge|k_2|/h-\nu_0h|k_1|$, and by~\eqref{eq:Xk-bound}
\[
\sum_{k_1}\sum_{|k_2|\ge N}|k|^2\,|X_k(0)|\,e^{|k||\Im\beta_0|}
\le\frac Ch\sum_{|k_2|\ge N}e^{-|k_2|(a/h-\kappa_0-1)}\le\frac{C'}h\,e^{-Na/h}.
\]
The estimate of $\partial_{\beta_1}^2X$ follows from~\eqref{eq:Xk-bound} without the exponential
factor. Finally,
$|\det\partial_\beta^2X|\le|\partial_{\beta_1}^2X|\,|\partial_{\beta_2}^2X|
+|\partial_{\beta_1}\partial_{\beta_2}X|^2$.
\end{proof}

\begin{remark}\label{rem:not-small}
It is the entries of the Hessian in the coordinates $\beta$ and the determinant that are
exponentially small, not the entries of the matrix $\partial_\alpha^2\cD(0,0)$ itself. Indeed, if
$\partial_{\beta_1}^2X\ne0$, the entries of this matrix contain corrections of order
$h^3\partial_{\beta_1}^2X$ coming from the off-diagonal entry
$\partial_\theta\mathfrak b_{0,1}=O(\nu h^2)$ of the matrix $\partial_\alpha\mathfrak b_0$.
\end{remark}
\section{The periodic problem and the Stokes constant}\label{sec:stokes}

In this section there is no weak perturbation. We consider the periodic problem~\eqref{eq:per}
with profile $zg$, where $g$ is a real even trigonometric polynomial with zero mean and
$|z|<R_0$. We use statement~(E1) with this $g$. First we describe the structure of the
difference of the generating functions on the vertical segment (Lemma~\ref{lem:difference}).
Then we introduce the Stokes constant $\chi[zg]$ and prove that the normalized fast entry of the
splitting matrix converges to it (Lemma~\ref{lem:stokes}). This is the first place where we use
statement~(E2).

\subsection{Structure of the difference of the generating functions}\label{ss:difference}

Let $T^u$ be the function from (E1). Set
\begin{equation}\label{eq:D-per}
T^s(u,\theta)=-T^u(-u,-\theta),\qquad\cD=T^u-T^s .
\end{equation}
By Remark~\ref{rem:periodic-T} this function $\cD$ differs from the difference~\eqref{eq:D-def}
constructed for the background $zg$ and $V=0$ by a constant. Unlike the latter, it is defined on
the whole intersection of the domains of definition of $T^u(u,\theta)$ and $T^u(-u,-\theta)$.
Since $g$ and $\psi_0$ are even, the function $T^s$ satisfies the same equation~\eqref{eq:HJ},
and $\partial_uT^s=\psi_0^2(1+E^s)$ with $E^s(u,\theta)=E^u(-u,-\theta)$. Subtracting the two
equations, we obtain
\begin{equation}\label{eq:transport-per}
\Bigl[(1+E)\,\partial_u+\frac1h\,\partial_\theta\Bigr]\cD=0,\qquad E=\tfrac12\bigl(E^u+E^s\bigr),
\qquad\partial_u\cD=\psi_0^2\,\bigl(E^u-E^s\bigr).
\end{equation}

Let $\kappa_0$, $\rho$, $\varsigma$ be chosen as in~\eqref{eq:P-def}, and $y_*=a-\kappa_0h$. As in
the proof of Lemma~\ref{lem:polydisc}, the function $E$ on the segment $u=iy$, $|y|\le y_*$,
satisfies the assumptions of Lemma~\ref{lem:straightening} with $E_1=0$ and a single angle
$\theta$. Let $W$, $H$, $W_0$, $\mathfrak b_y$, $\mathfrak a_y$ be the objects of that lemma, and
\begin{equation}\label{eq:Xi-def}
X(y,\beta)=\cD\bigl(iy,\mathfrak a_y(\beta)\bigr),\qquad
\Xi(y,\beta)=-\frac1h\,\partial_\beta X(y,\beta)=\sum_{k\ne0}\Xi_k(y)\,e^{ik\beta}.
\end{equation}

\begin{lemma}\label{lem:difference}
Let $\kappa_0\ge\kappa_0^*$ and, in addition, $2A_W/\kappa_0\le\rho_2/4$. There exist $C_7$ and
$h_7>0$ such that for $h<h_7$ and $|z|<R_0$ the following holds.
\begin{enumerate}[label=\textup{\arabic*.}]
\item $\partial_u\cD(iy,\theta)=\bigl(1+H(y,\theta)\bigr)\,\Xi\bigl(y,\mathfrak b_y(\theta)\bigr)$.
\item $\Xi_k(y)=\Xi_k(0)\exp\bigl[\tfrac kh\,\bigl(y-iW_0(y)\bigr)\bigr]$.
\item $|\Xi_k(y)|\le C_7\,h^{-2}\,e^{-3\rho_2|k|/4}
\exp\bigl[-\tfrac{|k|}h\,\bigl(y_*-y\operatorname{sgn}k\bigr)\bigr]$.
\item If $g$ has period $2\pi/N$, then $\Xi_k=0$ for $k$ not divisible by $N$, and for
$\theta\in\TT_{\rho_2/2}$
\begin{equation}\label{eq:EuEs}
|E^u-E^s|(iy,\theta)\le C_7\,\frac{\delta(y)^2}{h^{2}}\,e^{N(\kappa_0+1)}\,e^{-N\delta(y)/h}.
\end{equation}
\item $\partial_u\cD(-u,-\theta)=-\partial_u\cD(u,\theta)$ and
$\partial_\theta^2\cD(0,0)=-h\,\bigl(1+E(0,0)\bigr)\,\partial_\theta\partial_u\cD(0,0)$.
\end{enumerate}
\end{lemma}

\begin{proof}
1. We have
$\partial_u\cD=-i\,\partial_y\bigl[X(y,\mathfrak b_y(\theta))\bigr]
=-i\bigl[\partial_yX+\partial_\beta X\cdot\partial_y\mathfrak b_y\bigr]$. By
Lemma~\ref{lem:straightening}(5) with $S=0$, $-i\,\partial_yX=-(1-iW_0')\,\partial_\beta X/h$, and
$\partial_y\mathfrak b_y=-h^{-1}\partial_y(W-W_0)=-h^{-1}(iH-W_0')$. Adding up, we obtain
$-(1+H)\,\partial_\beta X/h$.

2. This follows from the equation $X_k'=\tfrac kh(1-iW_0')X_k$ and the equality
$\Xi_k=-\tfrac{ik}hX_k$.

3. By Lemma~\ref{lem:geometry}(3)
\[
|\partial_u\cD(iy,\theta)|\le|\psi_0(iy)|^2\cdot\frac{2C_1R_0h}{\delta(y)}\le\frac{2\pi^2C_1R_0h}{\delta(y)^3},
\]
since $|\psi_0(iy)|=2/\cos y\le\pi/\delta(y)$. At $y=\pm y_*$ this is at most
$Ch^{-2}\kappa_0^{-3}$, and $|1+H|\ge\tfrac12$. Hence $|\Xi_k(\pm y_*)|\le Ch^{-2}e^{-\rho_2|k|}$.
The transport by item~2 from the point $y=y_*\operatorname{sgn}k$ gives a factor not exceeding
\[
\exp\Bigl[-\frac{|k|}h\,\bigl(y_*-y\operatorname{sgn}k\bigr)
+\frac{|k|}h\,\bigl(|W_0(y)|+|W_0(y_*)|\bigr)\Bigr],
\]
and $\tfrac{|k|}h\cdot2\sup|W_0|\le2|k|A_W/\kappa_0\le\rho_2|k|/4$.

4. By uniqueness, all the constructions commute with the shift $\theta\mapsto\theta+2\pi/N$,
which gives the first claim. To prove the second claim, note that by items 1 and 3 the sum over
$k\ge N$ does not exceed $Ch^{-2}e^{N(\kappa_0+1)}e^{-N\delta(y)/h}$ for $y\ge0$, and the terms
with $k\le-N$ are smaller. The case $y<0$ is symmetric. It remains to divide by
$|\psi_0(iy)|^2=4/\cos^2y\ge4/\delta(y)^2$.

5. The function $\cD(u,\theta)=T^u(u,\theta)+T^u(-u,-\theta)$ is even with respect to the change
$(u,\theta)\mapsto(-u,-\theta)$. This gives the first equality and, in particular,
$\partial_u\cD(0,0)=0$. Differentiating equation~\eqref{eq:transport-per} with respect to
$\theta$ at the point $(0,0)$ and using this, we obtain the second equality.
\end{proof}

\subsection{The Stokes constant}\label{ss:stokes}

We keep the notation of Section~\ref{ss:difference}. Since the coefficients of $g$ are real, the
function $E^u$ has the conjugation symmetry
\begin{equation}\label{eq:conjugation}
\overline{E^u(\bar u,\bar\theta;\bar z)}=E^u(u,\theta;z).
\end{equation}
Indeed, conjugation takes solutions of~\eqref{eq:per} with parameters $(z,\theta)$ to solutions
with parameters $(\bar z,\bar\theta)$. Therefore, by the uniqueness in Lemma~\ref{lem:local}, the
function $\Gamma$ has this symmetry. By formula~\eqref{eq:coincidence}, so does the function
$E^u$ for $\Re u\le-\varrho_1'$, and by analyticity, on the whole connected symmetric domain
$D^{\rm loc}\cup D^{\rm out}$.

From~\eqref{eq:D-per} we get $\partial_uT^s(u,\theta)=(\partial_uT^u)(-u,-\theta)$. If
$u=i(a-xh)$, then the point conjugate to $-u$ equals $ia+h\zeta$ with $\zeta=-ix$, and for
$x>\kappa_6$ the point $-ix$ lies in the sector $S^u$. Therefore~\eqref{eq:conjugation} and
statement~(E2) imply that for $x>\kappa_6$ and $\theta\in\TT_\sigma$ the limit
\begin{equation}\label{eq:J-def}
\cJ(x,\theta;z)=\lim_{h\to0}h^2\,\partial_u\cD\bigl(i(a-xh),\theta;z\bigr)
=\Phi_0^u(-ix,\theta;z)-\overline{\Phi_0^u(-ix,-\bar\theta;\bar z)}
\end{equation}
exists. Set
\begin{equation}\label{eq:Xh-Uh}
\mathcal X_h(z)=i\,h^2e^{a/h}\,\Xi_1(0;z,h),\qquad
U_h(zg)=-\frac{\partial_\theta^2\cD(0,0)}{2L_h}.
\end{equation}
The quantity $U_h(zg)$ is the normalized fast entry of the splitting matrix of the periodic
problem with profile $zg$.

\begin{lemma}\label{lem:stokes}
There exist $C_8$ and a threshold for $h$ such that the following holds.
\begin{enumerate}[label=\textup{\arabic*.}]
\item $\bigl|h^2\,\partial_u\cD\bigl(i(a-xh),\theta\bigr)\bigr|\le C_8\,e^{-x}$ for
$\kappa_0+2\le x\le a/h$ and $\theta\in\TT_{\rho_2/4}$, uniformly in $h$ and $|z|<R_0$.
Consequently, $|\cJ(x,\theta;z)|\le C_8\,e^{-x}$, and the function $\cJ$ is analytic in
$(\theta,z)$.
\item The limit
\begin{equation}\label{eq:chi-def}
\chi[zg]=\lim_{x\to+\infty}i\,e^{x}\,\fc{\cJ}{1}(x;z)
\end{equation}
exists, where $\fc{\cJ}{1}$ is the first Fourier coefficient in $\theta$. The function
$z\mapsto\chi[zg]$ is analytic for $|z|<R_0$, and $\chi[0]=0$.
\item The function $\mathcal X_h$ is analytic in $z$, $|\mathcal X_h(z)|\le C_8$, and
$\mathcal X_h(z)\to\chi[zg]$ as $h\to0$ locally uniformly in the disc $|z|<R_0$.
\item $|U_h(zg)-\mathcal X_h(z)|\le C_8\,h$. In particular, $U_h(zg)\to\chi[zg]$ locally
uniformly, and all Taylor coefficients of the function $U_h(zg)$ in $z$ at zero converge to the
Taylor coefficients of the function $\chi[zg]$.
\end{enumerate}
\end{lemma}

Since $R_0$ is arbitrary, we obtain the following.

\begin{corollary}\label{cor:chi-entire}
The function $z\mapsto\chi[zg]$ is entire.
\end{corollary}

\begin{definition}\label{def:stokes}
The function $\chi[zg]$ defined by~\eqref{eq:chi-def} is called the \emph{first Stokes constant}
of the inner equation with profile $zg$.
\end{definition}

The name is explained as follows. In Section~\ref{ss:inner} we show that the function
$\varphi^{u}=\Phi_0^{u}+4\zeta^{-2}$ and its stable counterpart $\varphi^s$ satisfy the equation
obtained by differentiating with respect to $\zeta$ the inner equation
$(\partial_\theta+\partial_\zeta)\psi=\tfrac18\zeta^2(\partial_\zeta\psi)^2-2zg(\theta)/\zeta^2$,
written in the normalization of~\cite{BMS24}. Moreover, $\cJ=\varphi^u-\varphi^s$.
Formula~\eqref{eq:chi-def} means that $\fc{\cJ}{1}(x;z)=-i\chi[zg]\,e^{-x}(1+o(1))$ as
$x\to+\infty$. In other words, the constant $\chi$ measures the exponentially small difference
between two solutions of the inner equation.

\begin{proof}[Proof of Lemma~\ref{lem:stokes}]
1. Let $y=a-xh$. By item~3 of Lemma~\ref{lem:difference}
\begin{align*}
|\Xi_k(y)|&\le C_7h^{-2}e^{-3\rho_2k/4}\,e^{-k(x-\kappa_0)}&&(k\ge1),\\
|\Xi_k(y)|&\le C_7h^{-2}e^{-3\rho_2|k|/4}\,e^{-|k|y_*/h}&&(k\le-1).
\end{align*}
For $\theta\in\TT_{\rho_2/4}$ the point $\mathfrak b_y(\theta)$ lies in $\TT_{\rho_2/2}$, so the
series $\sum_k|\Xi_k(y)|e^{\rho_2|k|/2}$ converges and $h^2|\Xi(y,\mathfrak b_y(\theta))|\le Ce^{-x}$.
By item~1 of the same lemma this gives the estimate for $h^2\partial_u\cD$. The limit of uniformly
bounded analytic functions is analytic (Vitali's theorem).

2--3. The analyticity of $\mathcal X_h$ in $z$ follows from the analyticity of $T^u$ and of the
fixed points of Lemma~\ref{lem:straightening}, and the boundedness from item~3 of
Lemma~\ref{lem:difference} with $k=1$, $y=0$. By item~2 of the same lemma, for $y=a-xh$
\begin{equation}\label{eq:Xi1-transport}
h^2\,\Xi_1(y)=-i\,\mathcal X_h\,e^{-x}\,e^{-iW_0(y)/h},\qquad\frac{|W_0(y)|}h\le\frac{A_W}x .
\end{equation}
On the other hand, $\Xi_1(y)=\frac1{2\pi}\oint\Xi(y,\beta)\,e^{-i\beta}\,d\beta$. We replace the
contour by the image of the real circle under the map $\mathfrak b_y$ (it is homotopic to the
circle in the strip of analyticity) and use item~1 of Lemma~\ref{lem:difference}:
\[
\Xi_1(y)=\frac1{2\pi}\int_0^{2\pi}\frac{\partial_u\cD(iy,\theta)}{1+H(y,\theta)}\,
e^{-i\mathfrak b_y(\theta)}\,\mathfrak b_y'(\theta)\,d\theta .
\]
By Lemma~\ref{lem:straightening} the quantities $|H|$, $|\mathfrak b_y(\theta)-\theta|$ and
$|\mathfrak b_y'-1|$ do not exceed $C/x$. Together with item~1 this gives
$h^2\Xi_1(y)=\bigl(h^2\partial_u\cD(iy,\cdot)\bigr)^{[1]}+O(e^{-x}/x)$, that is,
by~\eqref{eq:Xi1-transport},
\begin{equation}\label{eq:Xh-via-D}
\mathcal X_h(z)=i\,e^{x}\,\bigl(h^2\,\partial_u\cD(i(a-xh),\cdot\,;z)\bigr)^{[1]}+O(1/x)
\end{equation}
uniformly in $h\le a/x$ and $|z|<R_0$. Set $\chi_x(z)=i\,e^x\,\fc{\cJ}{1}(x;z)$. Passing to the
limit $h\to0$ in~\eqref{eq:Xh-via-D} for fixed $x$ (dominated convergence under the integral over
$\theta$), we obtain
\[
\limsup_{h\to0}|\mathcal X_h(z)-\chi_x(z)|\le\frac Cx .
\]
Hence $|\chi_x-\chi_{x'}|\le C/x+C/x'$, so the limit $\chi[zg]=\lim_{x\to\infty}\chi_x(z)$
exists, and $\mathcal X_h(z)\to\chi[zg]$ pointwise. The family $\mathcal X_h$ is uniformly
bounded and analytic, so the convergence is locally uniform and the limit is analytic. For $z=0$
we have $\partial_u\cD\equiv0$, whence $\chi[0]=0$.

4. By item~5 of Lemma~\ref{lem:difference}, $\partial_u\cD(0,0)=0$. By item~1 of the same lemma
at $y=0$ we have
$\partial_u\cD(0,\theta)=\bigl(1+H(0,\theta)\bigr)\,\Xi\bigl(0,\mathfrak b_0(\theta)\bigr)$, and
since $|H|\le\tfrac12$, we get $\Xi(0,\beta_0)=0$, where $\beta_0=\mathfrak b_0(0)$.
Differentiating with respect to $\theta$ at zero,
\[
\partial_\theta\partial_u\cD(0,0)=\bigl(1+H(0,0)\bigr)\,\mathfrak b_0'(0)\,\partial_\beta\Xi(0,\beta_0).
\]
By Lemma~\ref{lem:straightening} at $y=0$ we have $|H(0,0)|\le Ch$, $|\beta_0|\le Ch$ and, by
Cauchy's inequality, $|\mathfrak b_0'(0)-1|\le Ch$. By item~3 of Lemma~\ref{lem:difference}
\[
\sum_{|k|\ge2}|k|\,|\Xi_k(0)|\,e^{|k||\Im\beta_0|}\le Ch^{-2}e^{-2a/h}.
\]
Therefore the equality $\Xi(0,\beta_0)=0$ gives
$\Xi_{-1}(0)\,e^{-i\beta_0}=-\Xi_1(0)\,e^{i\beta_0}+O(h^{-2}e^{-2a/h})$, and
\begin{align*}
\partial_\beta\Xi(0,\beta_0)&=i\,\Xi_1(0)\,e^{i\beta_0}-i\,\Xi_{-1}(0)\,e^{-i\beta_0}+O(h^{-2}e^{-2a/h})\\
&=2i\,\Xi_1(0)\,\bigl(1+O(h)\bigr)+O(h^{-2}e^{-2a/h}).
\end{align*}
Since $|\Xi_1(0)|\le C_8h^{-2}e^{-a/h}$ by item~3, we get
$\partial_\theta\partial_u\cD(0,0)=2i\,\Xi_1(0)+O(h^{-1}e^{-a/h})$. By item~5 of
Lemma~\ref{lem:difference} and the estimate $|E(0,0)|\le Ch$
\[
\partial_\theta^2\cD(0,0)=-h\,\bigl(1+O(h)\bigr)\bigl[2i\,\Xi_1(0)+O(h^{-1}e^{-a/h})\bigr]
=-2L_h\bigl[\mathcal X_h+O(h)\bigr].
\]
The convergence of the Taylor coefficients follows from the locally uniform convergence of
analytic functions.
\end{proof}

\begin{remark}\label{rem:straightened-phase}
The Fourier coefficients of the function $\partial_u\cD(0,\theta)$ in the angle $\theta$ itself do
not have this property, because the factor $(1+H)\,e^{i(\mathfrak b_0(\theta)-\theta)}$ is not
constant, and the first mode $\Xi_1$ contributes terms of order $h\,|\Xi_1|$ to all harmonics in
$\theta$. This is why in item~4 we carry out the argument in the straightened phase $\beta$.
\end{remark}

\begin{proposition}[normalization]\label{prop:normalization}
$[z^1]\,\chi[zg]=4\pi\,\fc g1$.
\end{proposition}

\begin{proof}
For small $|z|$ the same periodic problem can be constructed with zero background and the weak
term $zg$. By Lemma~\ref{lem:splitting}, applied to the decompositions $(zg,0)$ and $(0,zg)$, one
obtains the same function $\partial_\theta^2\cD(0,0)$. For zero background the perturbed solution
at $z=0$ is the separatrix itself, and formula~\eqref{eq:first-variation} for the first variation
with respect to the amplitude, together with the identity $1-\cos q^0(t)=2\cosh^{-2}t$, gives
\begin{gather*}
\partial_z\cD(0,\theta)\big|_{z=0}=\int_{-\infty}^{\infty}g(\theta+t/h)\,\frac{2\,dt}{\cosh^2t}
=\sum_\ell\fc g\ell\,e^{i\ell\theta}\,\hat K(\ell/h),\\
\hat K(s)=\int_{-\infty}^{\infty}\frac{2\cos(st)}{\cosh^2t}\,dt=\frac{2\pi s}{\sinh(\pi s/2)} .
\end{gather*}
Consequently,
$\partial_\theta^2\cD(0,0)=-z\sum_\ell\ell^2\,\fc g\ell\,\hat K(\ell/h)+O(z^2)$. The contribution
of the harmonics $\ell=\pm1$ equals $-2z\,\fc g1\cdot4\pi L_h\bigl(1+O(e^{-2a/h})\bigr)$, and the
others are exponentially smaller. Therefore $[z^1]\,U_h(zg)\to4\pi\,\fc g1$ as $h\to0$, and the
claim follows from item~4 of Lemma~\ref{lem:stokes}.
\end{proof}
\section{The first variation with respect to the perturbation profile}\label{sec:variation}

Here we estimate the derivative of the difference of the generating functions with respect to the
amplitude of the weak term at a fixed background (Lemma~\ref{lem:variation-bound}) and find the
limit of this derivative as $h\to0$ (Lemma~\ref{lem:variation-limit}). This is how the function
$G_{M,p}$ arises. Throughout this section the background $zM$ has period $\pi$ in $\theta$, and
$|z|<R_0$.

\subsection{Bound for the first variation}\label{ss:variation-bound}

Let $\phi(\alpha)=e^{i(jx+\ell\theta)}$, where $\ell$ is odd. Denote by $q^u_B(t;u,\theta)$,
$t\le0$, and $q^s_B(t;u,\theta)$, $t\ge0$, the unstable and stable half-trajectories of the
periodic problem with background $zM$ for which $q(0)=q^0(u)$ and the phase at $t=0$ equals
$\theta$. Set
\begin{equation}\label{eq:P-us}
\begin{gathered}
\mathcal P^u(u,\alpha)=\int_{-\infty}^0\phi(\alpha+\omega t)\bigl(1-\cos q^u_B(t)\bigr)\,dt,\\
\mathcal P^s(u,\alpha)=-\int_0^{\infty}\phi(\alpha+\omega t)\bigl(1-\cos q^s_B(t)\bigr)\,dt,
\end{gathered}
\end{equation}
and
\begin{equation}\label{eq:R-phi}
R_\phi=\mathcal P^u-\mathcal P^s .
\end{equation}
By formula~\eqref{eq:first-variation}, for a real even trigonometric polynomial
$f=\sum c_\phi\,\phi$ the sum $\sum c_\phi R_\phi$ is the derivative of the difference $\cD$ with
respect to the amplitude of the weak term $\eps f$ at $\eps=0$. In this section the frequency
$\nu$ may be complex with $|\nu|\le\nu_1$, where $|j|\,\nu_1\le\tfrac14$. Then
$\phi(\alpha+\omega t)=e^{ij(x+\nu t)}e^{i\ell(\theta+t/h)}$.

\begin{lemma}\label{lem:variation-bound}
There exist $C_9$ and $h_9>0$ such that for $h<h_9$, $|z|<R_0$ and $|\nu|\le\nu_1$
\begin{equation}\label{eq:variation-bound}
\sup_{\alpha\in\TT^2_{\rho_2/8}}|R_\phi(0,\alpha)|\le C_9\,L_h .
\end{equation}
\end{lemma}

The bound is of the right order. Each of the functions $\mathcal P^u$, $\mathcal P^s$ separately
is bounded but not small, while their difference is exponentially small.

\begin{proof}
\emph{Step 1: each branch.} Let $U^b(s)$, $b=u,s$, be the label trajectories of the background
half-trajectories. Along them
$\frac{d}{ds}\mathcal P^b\bigl(U^b(s),\alpha+\omega s\bigr)=\phi(\alpha+\omega s)\bigl(1-\cos q^0(U^b(s))\bigr)$,
that is,
\begin{equation}\label{eq:P-equation}
\bigl[(1+E^b)\,\partial_u+\omega\cdot\partial_\alpha\bigr]\mathcal P^b=\tfrac12\,\psi_0^2\,\phi,\qquad b=u,s.
\end{equation}
Let $u\in\cU^{1/8}_\kappa$ and $\alpha\in\TT^2_\sigma$. On the interval $[-t_*,0]$, by
Lemma~\ref{lem:geometry}, $|1-\cos q^0(U)|=\tfrac12|\psi_0(U)|^2\le C/d(t)^2$ and $|\phi|\le C$.
For $t\le-t_*$, by Lemma~\ref{lem:local}, $|1-\cos q|\le Ce^{t+t_*}$, while $|\phi|\le Ce^{|t|/4}$.
Therefore $|\mathcal P^b(u,\alpha)|\le C/\delta(u)$ and, by Cauchy's inequality,
$|\partial_u\mathcal P^b|\le C/\delta(u)^2$ on $\cU^{1/16}_{2\kappa}$.

\emph{Step 2: the equation for the difference.} Subtracting equations~\eqref{eq:P-equation} for
$b=u$ and $b=s$, we obtain
\begin{equation}\label{eq:R-equation}
\bigl[(1+E)\,\partial_u+\omega\cdot\partial_\alpha\bigr]R_\phi=S,\qquad
S=-\tfrac12\,(E^u-E^s)\,\partial_u\bigl(\mathcal P^u+\mathcal P^s\bigr),
\end{equation}
where $E=\tfrac12(E^u+E^s)$. By estimate~\eqref{eq:EuEs} with $N=2$ and by Step~1, on the
segment $u=iv$ for $\theta\in\TT_{\rho_2/2}$
\begin{equation}\label{eq:S-bound}
|S(iv,\alpha)|\le C\,h^{-2}\,e^{-2\delta(v)/h}.
\end{equation}
The estimate is obtained in the strip $|\Im\theta|<\rho_2/2$. Accordingly, the function
$\tilde S(v,\beta)=\bigl[S/(1+E)\bigr]\bigl(v,\mathfrak a_v(\beta)\bigr)$ is controlled for
$|\Im\beta_2|<\rho_2/3$, and its Fourier coefficients satisfy
\[
|\tilde S_{(j,n)}(v)|\le C\,e^{-\rho_2|n|/3}\,h^{-2}\,e^{-2\delta(v)/h}.
\]

\emph{Step 3: straightening.} We apply Lemma~\ref{lem:straightening} with $E=E_B$ and $E_1=0$. The
function $W$ depends only on $(y,\theta)$, so its construction does not use $\nu$, and the change
of variables $\beta_1=x-\nu\,(W-W_0)$, $\beta_2=\theta-(W-W_0)/h$ and the derivation of
formula~\eqref{eq:Xk-transport} remain valid for complex $\nu$ as well, since
$R_\phi=e^{ijx}\tilde R(u,\theta)$ is an entire function of $x$. The expansion of
$X(y,\beta)=R_\phi\bigl(iy,\mathfrak a_y(\beta)\bigr)$ contains only the modes $k=(j,n)$ with odd
$n$, because the shift $\theta\mapsto\theta+\pi$ preserves the background and $W$ and changes the
sign of $\phi$. Here $\lambda_k=j\nu+n/h$ and
\[
\bigl|e^{-\lambda_k(y-iW_0)}\bigr|\le\exp\bigl[-\Re\lambda_k\,y+|\lambda_k|\,|W_0|\bigr],\qquad
|\lambda_k|\,|W_0|\le\frac{\rho_2|n|}8+1,
\]
since $|W_0|\le A_Wh/\kappa_0$ and $2A_W/\kappa_0\le\rho_2/4$.

\emph{Step 4: estimate of the modes.} Let $n\ge1$, and take $y=y_*$ in~\eqref{eq:Xk-transport}. In
the boundary term we have $|X_k(y_*)|\le e^{-\rho_2|k|}\cdot2C/(\kappa_0h)$, and the transport
factor is at most $C\,e^{-ny_*/h}\,e^{\rho_2n/8}$. The integral term does not exceed
\[
C\,e^{-(\frac13-\frac18)\rho_2|n|}\,h^{-2}\,e^{-2a/h}\int_0^{y_*}e^{(2-n)v/h}\,dv
\le C\,e^{-5\rho_2|n|/24}\cdot
\begin{cases}h^{-1}e^{-a/h}\,e^{-\kappa_0},&n=1,\\[2pt] h^{-1}e^{-2a/h},&n\ge3.\end{cases}
\]
For $n\le-1$ we take $y=-y_*$ and obtain the same. Thus
\begin{equation}\label{eq:X-jn-bound}
|X_{(j,n)}(0)|\le\frac Ch\,e^{-5\rho_2|n|/24}\Bigl(e^{-|n|(a/h-\kappa_0)}+e^{-a/h}\Bigr).
\end{equation}

\emph{Step 5: the sum.} We have
$R_\phi(0,\alpha)=\sum_nX_{(j,n)}(0)\,e^{i(j\beta_1+n\beta_2)}$, where
$\beta=\mathfrak b_0(\alpha)$. For $\alpha\in\TT^2_{\rho_2/8}$ we have
$|\Im\beta|\le\rho_2/8+Ch$, and since $\tfrac5{24}>\tfrac18$, estimate~\eqref{eq:X-jn-bound}
implies that the series is majorized by $C_9\,h^{-1}e^{-a/h}$.
\end{proof}

\subsection{The limit of the first variation and the function \texorpdfstring{$G_{M,p}$}{G}}\label{ss:variation-limit}

Let $p$ be a real even trigonometric polynomial containing only odd harmonics. Set
$R_p=\sum_\ell\fc p\ell\,R_{e^{i\ell\theta}}$ (here $j=0$) and
\begin{equation}\label{eq:Gh-def}
G_h(z)=-\frac{\partial_\theta^2R_p(0,\theta;z)\big|_{\theta=0}}{2L_h},\qquad|z|<R_0 .
\end{equation}

For a profile $B$, an even trigonometric polynomial with harmonics from the union of the spectra
of $M$ and $p$, denote by
\begin{equation}\label{eq:Uh-B}
U_h(B)=-\frac{\partial_\theta^2\cD^{\rm per}_B(0,0)}{2L_h}
\end{equation}
the normalized fast entry of the splitting matrix of the periodic problem with profile $B$. We
construct the function $\cD^{\rm per}_B$ as in Sections~\ref{sec:manifolds}
and~\ref{sec:generating} for some decomposition of the profile $B$ into the background and the
weak part. We use three decompositions:
\begin{itemize}
\item[(I)] background $zg_t$, where $g_t=M+tp$ with real $t$, no weak part, $|z|<R_0$ (this is
the function $\cD$ of Section~\ref{sec:stokes}, constructed from the graph of (E1));
\item[(II)] background $zM$, weak part $rp$, $|z|<R_0$, $|r|\le e_2$;
\item[(III)] zero background, weak part $B$ with small coefficients (the sum of their absolute
values is at most $e_2$).
\end{itemize}
By Lemma~\ref{lem:splitting}, wherever two decompositions are applicable and the amplitudes of
the weak parts do not exceed $e_3$, the values of $U_h(B)$ coincide. We take the constants of the
lemmas of Section~\ref{sec:manifolds} common to the backgrounds involved
(Section~\ref{ss:constants}). For each $k$ only finitely many values of $t$ are involved below.

\begin{lemma}\label{lem:variation-limit}
\begin{enumerate}[label=\textup{\arabic*.}]
\item The function $G_h$ is analytic, and $|G_h(z)|\le C_{10}$ for $|z|<R_0$, $h<h_9$.
\item $G_h(z)=\partial_rU_h(zM+rp)\big|_{r=0}$, where $U_h(zM+rp)$ is constructed by
decomposition~(II).
\item As $h\to0$ the functions $G_h$ converge locally uniformly on $\CC$ to an entire function,
which we denote by $G_{M,p}$. Its Taylor coefficients are
\begin{equation}\label{eq:G-taylor}
[z^{k-1}]\,G_{M,p}=[t^1]\,[z^k]\,\chi\bigl[z(M+tp)\bigr],\qquad k\ge1,
\end{equation}
and $[z^k]\,\chi\bigl[z(M+tp)\bigr]$ is a polynomial in $t$ of degree at most $k$.
\end{enumerate}
\end{lemma}

\begin{proof}
1. This follows from Lemma~\ref{lem:variation-bound} with $j=0$ and Cauchy's inequality in
$\theta$.

2. By Lemma~\ref{lem:continuation}(3), applied to the background $zM$ and the weak term
$rp(\theta)$ (decomposition~(II)), the function $\cD^{\rm per}_{zM+rp}(0,\theta)$ is analytic in
$r$ in the disc $|r|\le e_2$ independent of $h$, and its derivative at $r=0$ equals
$R_p(0,\theta)$ by~\eqref{eq:first-variation}.

3. Fix $k$, real nodes $t_0=0,t_1,\dots,t_k$ and a sufficiently small $h$. We take the constants
of the lemmas of Section~\ref{sec:manifolds} common to the zero background and the backgrounds
$g_{t_0}=M,g_{t_1},\dots,g_{t_k}$. We apply the construction of Sections~\ref{sec:manifolds},
\ref{sec:generating} with zero background (then $E^u\equiv0$ and $U(t)=u_0+t$) and the weak term
$V=B$, that is, decomposition~(III). By Lemma~\ref{lem:continuation}(3) the resulting function
$U^{\rm III}_h(B)$ is analytic in the coefficients of $B$ in a neighborhood of zero, and
$U^{\rm III}_h(0)=0$, since for $K\equiv1$ the difference $\cD$ is constant. We expand it in
homogeneous polynomials as $U^{\rm III}_h(B)=\sum_{k\ge1}\mathcal Q_{k,h}(B)$. Then
\[
[z^k]\,U^{\rm III}_h\bigl(z(M+tp)\bigr)=\mathcal Q_{k,h}(M+tp)
\]
is a polynomial in $t$ of degree at most $k$.

\emph{Comparison with (II).} By Lemma~\ref{lem:splitting}, applied to the decompositions
$(zM,rp)$ and $(0,zM+rp)$, for small $|z|$, $|r|$ we have $U_h(zM+rp)=U^{\rm III}_h(zM+rp)$.
Differentiating with respect to $r$ at $r=0$ and using item~2, we obtain
$G_h(z)=\sum_kz^{k-1}D\mathcal Q_{k,h}(M)[p]$ for small $|z|$, that is,
\begin{equation}\label{eq:Gh-coeff}
[z^{k-1}]\,G_h=[t^1]\,\mathcal Q_{k,h}(M+tp).
\end{equation}

\emph{Comparison with (I).} For every real $t$ the profile $g_t=M+tp$ satisfies the assumptions
of Section~\ref{sec:stokes}. The function $U_h(zg_t)$ from~\eqref{eq:Xh-Uh} is defined through
$\cD$ and does not depend on the choice of the constants of the lemmas of
Section~\ref{sec:manifolds}, because the graph $E^u$ is uniquely determined for
$\Re u\le-\varrho_1'$ by the uniqueness in Lemma~\ref{lem:local} (the local manifold does not
change when $q_1$ is decreased) and is continued analytically. The function $U_h(zg_t)$ is
analytic in $z$ for $|z|<R_0$ and coincides with the function $U_h$ constructed by the
decomposition $(zg_t,0)$ (Remark~\ref{rem:periodic-T}). At the nodes, by
Lemma~\ref{lem:splitting} applied to the decompositions $(zg_{t_i},0)$ and $(0,zg_{t_i})$, it
equals $U^{\rm III}_h(zg_{t_i})$ for small $|z|$. Therefore
\[
\mathcal Q_{k,h}(M+t_ip)=[z^k]\,U_h(zg_{t_i}),\qquad i=0,\dots,k .
\]

The same equality holds for every real $t'$. Indeed, add the node $t'$ and construct, with new
common constants, the function of decomposition~(III) and its polynomials $\mathcal Q'_{k,h}$.
Then $\mathcal Q'_{k,h}(M+tp)=[z^k]\,U_h(zg_t)$ for $t=t_0,\dots,t_k,t'$. The polynomials
$\mathcal Q_{k,h}(M+tp)$ and $\mathcal Q'_{k,h}(M+tp)$ in $t$ have degree at most $k$ and coincide
at the $k+1$ points $t_0,\dots,t_k$. Hence they are equal, and
$\mathcal Q_{k,h}(M+t'p)=[z^k]\,U_h(zg_{t'})$.

Thus, by Lemma~\ref{lem:stokes}(4), for every real $t$
\[
\mathcal Q_{k,h}(M+tp)=[z^k]\,U_h(zg_t)\ \longrightarrow\ [z^k]\,\chi[zg_t]\qquad(h\to0)
\]
(the threshold for $h$ depends on $t$). Polynomials of degree at most $k$ that converge at $k+1$
points converge coefficientwise. Their limit is a polynomial in $t$ of degree at most $k$ equal to
$[z^k]\,\chi[zg_t]$ for all real $t$. Therefore~\eqref{eq:Gh-coeff} implies
$[z^{k-1}]\,G_h\to[t^1]\,[z^k]\,\chi\bigl[z(M+tp)\bigr]$ for every $k$. By item~1 the family
$\{G_h\}$ is normal in every disc. All its subsequential limits have the same Taylor coefficients
at zero, hence the limit is unique. The radius $R_0$ is arbitrary.
\end{proof}

\begin{remark}\label{rem:joint-analyticity}
If the function $\chi$ were jointly analytic in the coefficients of the profile, the right-hand
side of~\eqref{eq:G-taylor} would be a Taylor coefficient of the directional derivative, and we
would have $G_{M,p}(\mu)=\partial_r\chi[\mu M+rp]\big|_{r=0}$. The proof does not use this. We
apply only the one-parameter periodic theory for the profiles $z(M+tp)$ with real $t$.
\end{remark}

\begin{corollary}\label{cor:G-zero}
$G_{M,p}(0)=4\pi\,\fc p1$.
\end{corollary}

\begin{proof}
By~\eqref{eq:G-taylor} and Proposition~\ref{prop:normalization}
$G_{M,p}(0)=[t^1]\,[z^1]\,\chi\bigl[z(M+tp)\bigr]=[t^1]\,4\pi\,\bigl(\fc M1+t\,\fc p1\bigr)=4\pi\,\fc p1$.
\end{proof}

Equality~\eqref{eq:G-taylor} coincides with equality~\eqref{eq:response} of
Theorem~\ref{thm:series}. Until Section~\ref{sec:borel}, $G_{M,p}$ means the limit of
Lemma~\ref{lem:variation-limit}. We prove in Section~\ref{sec:borel} that this function is given
by the series~\eqref{eq:G}.
\section{Asymptotics of the splitting determinant}\label{sec:determinant}

In this section we prove Propositions~\ref{thm:Qprime} and~\ref{thm:Q}. Theorem~\ref{thm:main} was
derived from Proposition~\ref{thm:Qprime} in Section~\ref{ss:thmQ}.

We begin with an auxiliary estimate. Let $q^u_B(t)$, $t\le0$, and $q^s_B(t)$, $t\ge0$, be the
background half-trajectories of the periodic problem with background $zM$, $|z|\le R_0$, for
which $q(0)=\pi$ and the phase at $t=0$ equals $\theta\in\TT_\sigma$.

\begin{lemma}\label{lem:real-label}
Let $\Psi(t)$ be a continuous function on $\RR$ with $|\Psi(t)|\le C_\Psi$. Then
\begin{equation}\label{eq:real-label}
\Bigl|\int_{-\infty}^0\Psi(t)\bigl(1-\cos q^u_B(t)\bigr)\,dt
+\int_0^{\infty}\Psi(t)\bigl(1-\cos q^s_B(t)\bigr)\,dt
-\int_{-\infty}^{\infty}\Psi(t)\,\frac{2\,dt}{\cosh^2t}\Bigr|\le C\,C_\Psi\,h
\end{equation}
with a constant $C$ independent of $h$, $z$, $\theta$ and $\Psi$.
\end{lemma}

\begin{proof}
We have $1-\cos q^u_B(t)=2\cosh^{-2}U(t)$, where $U$ is the background label trajectory with
$U(0)=0$. By Lemma~\ref{lem:geometry}(1), $|U(t)-t|\le C_1R_0h\log(1+t_*/a)$ for $-t_*\le t\le0$.
For $t\le-t_*$ the trajectory lies in $D^{\rm loc}$, where $|E^u|\le C_ER_0h$ by~\eqref{eq:E1}.
Hence $|U(t)-t|\le Ch(1+|t|)$ for all $t\le0$, whence
\[
\bigl|2\cosh^{-2}U(t)-2\cosh^{-2}t\bigr|\le Ch\,(1+|t|)\,e^{-1.8|t|},
\]
and the difference of the integrals over $t\le0$ does not exceed $C\,C_\Psi\,h$. (The tails
$t\le-t_*$ themselves are not small, since $t_*$ is fixed.) The stable branch is reduced to the
unstable one by reversibility~\eqref{eq:reversibility}.
\end{proof}

\subsection{General even perturbations}\label{ss:general-f}

We prove Proposition~\ref{thm:Qprime}. Let $V=\tau f(x,\theta)$, and let $f=f_E+f_O$ be the
decomposition into the parts with even and odd Fourier indices in $\theta$. Then
\[
A_f=\partial_x^2\langle f\rangle_\theta\big|_{x=0},\qquad p_f(\theta)=f_O(0,\theta).
\]
First we fix a real $\mu$ and take $R_0>|\mu|$ in the lemmas. We discuss uniformity in $\mu$ in
Section~\ref{ss:proof-main}. Let
\[
\mathbf D(\tau)=\partial_\alpha^2\cD(0,0),\qquad F(\tau)=\det\mathbf D(\tau).
\]
By Section~\ref{ss:detdef}, the matrix $\mathbf D(\tau)$ is analytic for $|\tau|\le e_2$ and
bounded there uniformly in $h$. Write
\[
\mathbf D(\tau)=\mathbf D_0+\tau\,\mathbf H+\tau^2\,\mathbf K_2+\dots,
\]
where $|\mathbf H|,|\mathbf K_2|\le C$ by Cauchy's inequality.

\paragraph{Step 1: zeroth order} For $\tau=0$ the function $K$ does not depend on $x$, hence
\[
\mathbf D_0=\diag(0,d_B),\qquad|d_B|\le Ch^{-1}e^{-2a/h}.
\]
The estimate follows from Lemma~\ref{lem:polydisc} with $N=2$ and zero weak perturbation, since
the background has period $\pi$. Consequently,
\begin{equation}\label{eq:F-coeffs}
[\tau^1]\,F=d_B\,\mathbf H_{11},\qquad[\tau^2]\,F=\det\mathbf H+d_B\,(\mathbf K_2)_{11}.
\end{equation}

\paragraph{Step 2: the odd part} By formula~\eqref{eq:first-variation}, the matrix $\mathbf H$ is
the second derivative in $\alpha$ at zero of the first variation
$\sum_{j,\ell}\fc f{j,\ell}R_{e^{i(jx+\ell\theta)}}(0,\alpha)$. Let
$\mathbf H=\mathbf H_E+\mathbf H_O$, where $\mathbf H_E$ and $\mathbf H_O$ correspond to the sums
over even and odd $\ell$. By Lemma~\ref{lem:variation-bound} (real $\nu$, finite spectrum) and
Cauchy's inequality
\begin{equation}\label{eq:HO-bound}
|\mathbf H_O|\le C\,L_h .
\end{equation}
Consider the entry with indices $22$. For $\nu=0$ the sum
$\sum_{j,\ell}\fc f{j,\ell}R_{e^{i(jx+\ell\theta)}}(0,0,\theta)$ over odd $\ell$ equals
$R_{p_f}(0,\theta)$. By Lemma~\ref{lem:variation-bound} this sum is analytic in complex $\nu$ in
the disc $|\nu|\le\nu_1$ and bounded there by $CL_h$. By Cauchy's inequality its values at $\nu$
and at $0$ differ by at most $C|\nu|L_h$, together with the derivatives in $\theta$. Since
$|\nu|\le\nu_0h$, by Lemma~\ref{lem:variation-limit}
\begin{equation}\label{eq:HO22}
(\mathbf H_O)_{22}=\partial_\theta^2R_{p_f}(0,0)+O(hL_h)=-2L_h\bigl(G_{M,p_f}(\mu)+o(1)\bigr).
\end{equation}

\paragraph{Step 3: the even part} We show that
\begin{equation}\label{eq:HE}
\mathbf H_E=\diag(4A_f,0)+O(h).
\end{equation}
The first variation for $f_E$ on the section $u=0$ equals the sum of the integrals over the two
background half-trajectories with integrand $f_E(\alpha+\omega t)$. For
$\alpha\in\TT^2_{\sigma/2}$ this function is bounded, and by Lemma~\ref{lem:real-label} the first
variation equals
\[
\int_{-\infty}^{\infty}f_E(\alpha+\omega t)\,\frac{2\,dt}{\cosh^2t}+O(h)
=\sum_{j}\sum_{\ell\ \text{even}}\fc f{j,\ell}\,e^{i(jx+\ell\theta)}\,\hat K(j\nu+\ell/h)+O(h)
\]
uniformly on $\TT^2_{\sigma/2}$, where $\hat K(s)=2\pi s/\sinh(\pi s/2)$. Since
$0<\hat K(s)\le4.2\,\pi|s|e^{-\pi|s|/2}$ for $|s|\ge1$ and $\hat K(s)=4+O(s^2)$ for small $s$,
the sum equals $4\langle f\rangle_\theta(x)+O(h^2)$. Cauchy's inequality transfers this to the
second derivatives at zero. Note that the entries $(\mathbf H_E)_{12}$ and $(\mathbf H_E)_{22}$
are only polynomially small.

\paragraph{Step 4: the determinant of the even part} The determinant $\det\mathbf H_E$ is
exponentially small. Indeed, apply~\eqref{eq:F-coeffs} to the problem with $f_E$ in place of $f$,
in which the function $K$ has period $\pi$ in $\theta$. We obtain
$\det\mathbf H_E=[\tau^2]\,F_E-d_B\,(\mathbf K_2^E)_{11}$. By Lemma~\ref{lem:polydisc} with $N=2$
and Cauchy's inequality on the disc $|\tau|\le h^P$
\begin{equation}\label{eq:detHE}
|\det\mathbf H_E|\le Ch^{-2-2P}e^{-2a/h}.
\end{equation}

\paragraph{Step 5: the quadratic coefficient} For symmetric matrices of order two
\[
\det(\mathbf H_E+\mathbf H_O)=\det\mathbf H_E+\det\mathbf H_O
+(\mathbf H_E)_{11}(\mathbf H_O)_{22}+(\mathbf H_E)_{22}(\mathbf H_O)_{11}
-2(\mathbf H_E)_{12}(\mathbf H_O)_{12}.
\]
Substituting \eqref{eq:HO-bound}--\eqref{eq:detHE} and using~\eqref{eq:F-coeffs}, we obtain
\begin{equation}\label{eq:tau2}
[\tau^2]\,F=\bigl(-8A_f\,G_{M,p_f}(\mu)+o(1)\bigr)\,L_h .
\end{equation}

\paragraph{Step 6: the remainder} By Lemma~\ref{lem:polydisc} with $N=1$,
$|F(\tau)|\le C_6h^{-2}e^{-a/h}$ for $|\tau|\le h^P$, and by Cauchy's inequality
$|[\tau^k]\,F|\le C_6h^{-2}e^{-a/h}(2/h^P)^k$. Therefore for $|\tau|\le h^P/4$ the sum of the
terms with $k\ge3$ does not exceed $Ch^{-2-3P}e^{-a/h}|\tau|^3$, and the linear term,
by~\eqref{eq:F-coeffs}, does not exceed $C|\tau|h^{-1}e^{-2a/h}$. Consequently,
\[
\frac{F(\tau)}{\tau^2L_h}=-8A_f\,G_{M,p_f}(\mu)+o(1)+O\bigl(h^{-1-3P}|\tau|\bigr)
+O\Bigl(\frac{e^{-a/h}}{|\tau|}\Bigr),
\]
and both remainders tend to zero uniformly for $e^{-a/(2h)}\le|\tau|\le h^{3P+2}$. This is
exactly~\eqref{eq:Qprime}.

The lower bound in~\eqref{eq:tau-range} is needed, because the term
$\tau\,d_B\,\mathbf H_{11}$, linear in $\tau$, is not small compared with $\tau^2L_h$ for
exponentially small $\tau$.

\subsection{Two independent weak amplitudes}\label{ss:two-amplitudes}

We prove Proposition~\ref{thm:Q}. Let $M=\cos4\theta$, $V=s\cos x+r\cos\theta$, let $\mu$ be real,
$R_0=|\mu|+1$, $G=G_{\cos4\theta,\cos\theta}$, $\mathbf D=\partial_\alpha^2\cD(0,0)$ and
$F=\det\mathbf D$.

\paragraph{Step 1: structure at $s=0$} For $s=0$ the function $K$ does not depend on $x$. By
uniqueness, all the constructions do not depend on $x$, hence $\mathbf D_{11}=\mathbf D_{12}=0$
and $F(0,r)=0$. The function $F$ is analytic in the polydisc $|s|+|r|\le e_2$. Differentiating
the determinant, we obtain the exact equality
\begin{equation}\label{eq:dsF}
\partial_sF(0,r)=\partial_s\mathbf D_{11}(0,r)\cdot\mathbf D_{22}(0,r).
\end{equation}
By formula~\eqref{eq:first-variation}
\[
\partial_s\cD(0,\alpha)\big|_{s=0}=\int_{-\infty}^0\cos(x+\nu t)\bigl(1-\cos q^u_r\bigr)\,dt
+\int_0^{\infty}\cos(x+\nu t)\bigl(1-\cos q^s_r\bigr)\,dt,
\]
where $q^{u}_r$, $q^s_r$ are the half-trajectories of the periodic problem with profile
$\mu M+r\cos\theta$ and $q(0)=\pi$. They do not depend on $x$. Consequently,
\begin{equation}\label{eq:Kh}
\begin{gathered}
\partial_s\mathbf D_{11}(0,r)=-\mathcal K_h(r),\\
\mathcal K_h(r)=\int_{-\infty}^0\cos(\nu t)\bigl(1-\cos q^u_r\bigr)\,dt
+\int_0^{\infty}\cos(\nu t)\bigl(1-\cos q^s_r\bigr)\,dt,
\end{gathered}
\end{equation}
and $\mathbf D_{22}(0,r)=\partial_\theta^2\cD^{\rm per}_{\mu M+r\cos\theta}(0,0)=:D^{\rm per}_{22}(r)$.

\paragraph{Step 2: the mixed derivative} From~\eqref{eq:dsF} and~\eqref{eq:Kh}
\[
\partial_r\partial_sF(0,0)=-\mathcal K_h'(0)\,D^{\rm per}_{22}(0)-\mathcal K_h(0)\,\partial_rD^{\rm per}_{22}(0).
\]
Here:
\begin{itemize}
\item $D^{\rm per}_{22}(0)$ is the fast entry for the background $\mu\cos4\theta$ alone, and by
Lemma~\ref{lem:polydisc} with $N=4$ we have $|D^{\rm per}_{22}(0)|\le Ch^{-1}e^{-4a/h}$;
\item $|\mathcal K_h'(0)|\le C$, since the matrix $\mathbf D$ is analytic and bounded uniformly in
$h$ in the polydisc $|s|+|r|\le e_2$;
\item $\mathcal K_h(0)=4+O(h)$ by Lemma~\ref{lem:real-label} with $\Psi(t)=\cos\nu t$, since
$\cos\nu t=1+O(h^2t^2)$ and $\int_{-\infty}^{\infty}2\cosh^{-2}t\,dt=4$;
\item $\partial_rD^{\rm per}_{22}(0)=\partial_\theta^2R_{\cos\theta}(0,0)=-2L_hG_h(\mu)$ by item~2
of Lemma~\ref{lem:variation-limit}, and $G_h(\mu)\to G(\mu)$ by item~3 of the same lemma.
\end{itemize}
Thus
\begin{equation}\label{eq:mixed}
\partial_r\partial_sF(0,0)=\bigl(8\,G(\mu)+o(1)\bigr)\,L_h .
\end{equation}

\paragraph{Step 3: estimates on the polydisc} By Lemma~\ref{lem:polydisc} with $N=1$,
$|F(s,r)|\le C_6h^{-2}e^{-a/h}$ for $|s|+|r|\le h^P$. For $r=0$ the function $K$ has period
$\pi/2$ in $\theta$, and by the same lemma with $N=4$, $|F(s,0)|\le C_6h^{-2}e^{-4a/h}$ for
$|s|\le h^P$.

\paragraph{Step 4: the remainder} In the bidisc $|s|,|r|\le h^P/2$ expand
$F=\sum_{i\ge1,\,j\ge0}c_{ij}\,s^ir^j$. There are no terms with $i=0$ by Step~1. By Cauchy's
inequality $|c_{ij}|\le C_6h^{-2}e^{-a/h}(2/h^P)^{i+j}$. For $|s|,|r|\le h^P/8$
\[
\bigl|F(s,r)-F(s,0)-c_{11}\,sr\bigr|\le\sum_{\substack{i,j\ge1\\ i+j\ge3}}|c_{ij}|\,|s|^i|r|^j
\le C\,h^{-2-3P}e^{-a/h}\,|sr|\,(|s|+|r|),
\]
and by Schwarz's lemma, since $F(0,0)=0$,
$|F(s,0)|\le C_6h^{-2}e^{-4a/h}\cdot2|s|/h^P$. Since $c_{11}=\partial_r\partial_sF(0,0)$,
\eqref{eq:mixed} gives
\[
\frac{F(s,r)}{s\,r\,L_h}=8\,G(\mu)+o(1)+O\bigl(h^{-1-3P}(|s|+|r|)\bigr)
+O\Bigl(\frac{h^{-1-P}e^{-3a/h}}{|r|}\Bigr).
\]
Both remainders tend to zero for $|s|,|r|\le h^{3P+2}$, $|r|\ge e^{-2a/h}$, $s\ne0$. This is
exactly~\eqref{eq:Q}. The equality $G(0)=2\pi$ follows from Corollary~\ref{cor:G-zero}.

\subsection{Uniformity in the background amplitude}\label{ss:proof-main}

In Section~\ref{ss:general-f} the amplitude $\mu$ was fixed. We show that the remainder $o(1)$
in~\eqref{eq:Qprime} is uniform in $\mu\in[-R_0,R_0]$. We apply the lemmas of
Sections~\ref{sec:manifolds}--\ref{sec:variation} with $R_0+1$ in place of $R_0$. All estimates of
these lemmas are uniform in complex $z$ with $|z|\le R_0+1$, that is, the constants depend on
$R_0$ but not on $z$. Therefore the quantities $O(\cdot)$ in Steps 1--6 of
Section~\ref{ss:general-f} are uniform in $\mu$.
The only passage to the limit is the convergence $G_h(\mu)\to G_{M,p_f}(\mu)$
in~\eqref{eq:HO22}. By item~3 of Lemma~\ref{lem:variation-limit} it is locally uniform, hence
uniform on the interval $[-R_0,R_0]$. Uniformity in $\nu$ with $|\nu|\le\nu_0h$ follows from the
fact that the constants of the lemmas depend only on $\nu_0$. Proposition~\ref{thm:Qprime} is
proved, and with it Theorem~\ref{thm:main}.

\begin{remark}\label{rem:P-value}
The exponent $P=\Lambda+5$ depends neither on $\mu$ nor on the form of the weak perturbation,
since $\Lambda=2C_S+1$ is determined only by the geometry of Lemma~\ref{lem:geometry}. The
constants and the threshold for $h$ depend on $\mu$. For $C_S=46$ we get $\Lambda=93$ and
$P=98$. The value of $P$ is not optimal, because it comes from the crude Gronwall estimate in
Lemma~\ref{lem:continuation}.
\end{remark}

\begin{remark}\label{rem:complex-mu}
The proofs of formulas~\eqref{eq:Qprime} and~\eqref{eq:Q} do not use the reality of $\mu$. The
lemmas are stated for complex $|z|<R_0$, so both formulas hold also for a fixed complex $\mu$.
Only the parameter $t$ in Lemma~\ref{lem:variation-limit} must be real, because the profile in
statements (E1) and (E2) has real coefficients.
\end{remark}
\section{A series for the response of the Stokes constant}\label{sec:borel}

In this section we prove that the function $G_{M,p}$ of Lemma~\ref{lem:variation-limit} is given
by the series~\eqref{eq:G}. The main step is the following proposition on the Taylor
coefficients of the Stokes constant. In it the background may have any frequency lattice step
$N\ge2$.

Let $M$ be a real even trigonometric polynomial with zero mean all of whose frequencies are
multiples of an integer $N\ge2$. Define the functions $\hY^{[n]}_{(j)}$, $\hQ^{[n]}_{(j)}$,
$\hE^{[n]}_{(j)}$ by \eqref{eq:Y1}--\eqref{eq:E}, in which the indices $n$, $n_i$ range over
$N\ZZ$, and set $\hE^{[n]}_{(j)}=0$ for $n\notin N\ZZ$. Let
\[
\Omega_N=\CC\setminus\{s\in\RR:\ |s|\ge N\}.
\]

\begin{proposition}\label{thm:S}
Let $M$ be as above, and let $p$ be a real even trigonometric polynomial with zero mean. Then all
the functions $\hY^{[n]}_{(j)}$ and $\hE^{[n]}_{(j)}$ are analytic in $\Omega_N$, and for every
$k\ge1$
\begin{equation}\label{eq:S1}
[t^1]\,[z^k]\,\chi\bigl[z(M+tp)\bigr]=4\pi\sum_m\fc pm\,\hE^{[1-m]}_{(k-1)}(1).
\end{equation}
\end{proposition}

The proof occupies Sections~\ref{ss:inner}--\ref{ss:jump}. Its plan is as follows. The limit
functions from (E2) satisfy the inner equation and are uniquely determined by their decay at
infinity (Lemmas~\ref{lem:S1}, \ref{lem:S2}). Therefore we can construct their Taylor
coefficients in the amplitude explicitly, as Laplace integrals of the functions
$\hY^{[n]}_{(j)}$ along rays in the upper or lower half-plane. The difference of the two branches
in the first harmonic equals the residue at the point $s=1$.

\subsection{Inner solutions and their Taylor coefficients}\label{ss:inner}

Let $g$ be a real even trigonometric polynomial with zero mean, and $|z|<R_0$. By statement~(E2)
and formula~\eqref{eq:J-def}, the functions
\begin{gather*}
\Phi_0^u(\zeta,\theta;z)=\lim_{h\to0}h^2\,\partial_uT^u(ia+h\zeta,\theta;z)\quad(\zeta\in S^u),\\
\Phi_0^s(\zeta,\theta;z)=\overline{\Phi_0^u(-\bar\zeta,-\bar\theta;\bar z)}\quad(\zeta\in S^s),
\end{gather*}
are defined, where $S^s=\{-\bar\zeta:\ \zeta\in S^u\}$. Moreover,
$\cJ(x,\theta;z)=(\Phi_0^u-\Phi_0^s)(-ix,\theta;z)$. Write
\begin{equation}\label{eq:phi-def}
\Phi^b_0=-\frac4{\zeta^{2}}+\varphi^b,\qquad b=u,s.
\end{equation}

\begin{lemma}\label{lem:S1}
The functions $\varphi^b$ are analytic in $(\zeta,\theta,z)$ on
$S^b\times\TT_\sigma\times\{|z|<R_0\}$ and satisfy the estimate
\begin{equation}\label{eq:phi-bound}
|\varphi^b(\zeta,\theta;z)|\le C_0\,|z|\,|\zeta|^{-3},\qquad C_0=\frac{4C_E}{\pi^3},
\end{equation}
and the equation
\begin{equation}\label{eq:inner}
(\partial_\theta+\partial_\zeta)\varphi=\partial_\zeta\bigl(\tfrac18\,\zeta^2\varphi^2\bigr)
+\frac{4zg(\theta)}{\zeta^{3}} .
\end{equation}
\end{lemma}

Equation~\eqref{eq:inner} is obtained by differentiating with respect to $\zeta$ the inner
equation $(\partial_\theta+\partial_\zeta)\psi=\tfrac18\zeta^2(\partial_\zeta\psi)^2-2zg(\theta)/\zeta^2$
for $\varphi=\partial_\zeta\psi$.

\begin{proof}
Let $\varphi_h(\zeta,\theta)=h^2\psi_0^2E^u$ at the point $u=ia+h\zeta$. Since
$\cosh(ia+h\zeta)=i\sinh h\zeta$, we have
$h^2\psi_0^2(ia+h\zeta)=-4h^2/\sinh^2(h\zeta)\to-4\zeta^{-2}$ together with the derivatives,
locally uniformly. For $\zeta\in S^u$ and small $h$ the point $ia+h\zeta$ lies in $D^{\rm out}$
(Section~\ref{ss:external}). Since $|\psi_0|^2=4/|\cosh u|^2$, the second
estimate~\eqref{eq:E1} gives
\[
|\varphi_h|\le\frac{4C_E\,|z|\,h^3}{|u^2+a^2|^3}\ \longrightarrow\ \frac{4C_E\,|z|}{(2a|\zeta|)^3}.
\]
Hence on every compact set $\mathcal K\subset S^u\times\TT_\sigma\times\{|z|<R_0\}$ the functions
$\varphi_h$ are defined and uniformly bounded for $h$ small enough. By (E2) they converge
pointwise to $\varphi^u$. By Vitali's theorem the convergence is locally uniform, and the limit is
analytic and satisfies~\eqref{eq:phi-bound}.

Let us derive the equation. From~\eqref{eq:HJ}, for $\tilde T=T^u-T_0$ we get
$h^{-1}\partial_\theta\tilde T+\partial_u\tilde T+(\partial_u\tilde T)^2/(2\psi_0^2)
=\tfrac12zg\,\psi_0^2$. Differentiate with respect to $u$ and substitute
$\partial_u\tilde T=h^{-2}\varphi_h$, $\partial_u=h^{-1}\partial_\zeta$:
\[
\partial_\theta\varphi_h+\partial_\zeta\varphi_h
+\partial_\zeta\Bigl[\frac{\varphi_h^2}{2h^2\psi_0^2}\Bigr]
=\frac{zg}2\,\partial_\zeta\bigl(h^2\psi_0^2\bigr).
\]
The limit $h\to0$ gives~\eqref{eq:inner}. For $b=s$ note that the transformation
$\varphi(\zeta,\theta;z)\mapsto\overline{\varphi(-\bar\zeta,-\bar\theta;\bar z)}$ preserves
equation~\eqref{eq:inner} for real even $g$ (both sides of the equation change sign) and
estimate~\eqref{eq:phi-bound}.
\end{proof}

\paragraph{Variables} Set
\begin{equation}\label{eq:gamma-def}
\zeta=-ix,\qquad\gamma^b=\tfrac14\,\zeta^2\varphi^b .
\end{equation}
Then $|\gamma^b|\le C_0|z|/(4|x|)$, and equation~\eqref{eq:inner} takes the form
\begin{equation}\label{eq:gamma-eq}
x\bigl[(1-\gamma)\,\gamma_x-i\gamma_\theta\bigr]=2\gamma-\gamma^2+zg .
\end{equation}
In the $x$-plane the sector $S^u$ has vertex $\kappa_6$ and directions
$\arg(x-\kappa_6)\in(-\tfrac\pi2-\beta_0,\ \tfrac\pi2-\beta_2)$, and the sector $S^s$ has
directions $(\beta_2-\tfrac\pi2,\ \tfrac\pi2+\beta_0)$. Both sectors contain the ray
$x\to+\infty$ and rays along which $\Re x\to-\infty$.

\begin{lemma}[uniqueness]\label{lem:S2}
Let $S$ be an open sector of the $x$-plane containing the ray $x\to+\infty$ and a ray along which
$\Re x\to-\infty$. Let a function $d(x,\theta)$ be analytic on $S\times\TT_{\sigma'}$,
$|d|\le C/|x|$ and $x(\partial_x-i\partial_\theta)d-2d=0$. Then $d\equiv0$.
\end{lemma}

\begin{proof}
For the mode $d_n(x)e^{in\theta}$ we get $x(d_n'+nd_n)=2d_n$, that is, $d_n=c_nx^2e^{-nx}$. For
$n\le0$ this function does not decay along the ray $x\to+\infty$, and for $n>0$ along the ray of
the sector with $\Re x\to-\infty$. Hence $c_n=0$.
\end{proof}

The sectors $S^b$, the strip $\TT_\sigma$ and the constant $C_0$ of Lemma~\ref{lem:S1} depend on
the profile. Below, for the profiles $g_t=M+tp$, we use a finite set of values of $t$ in each
order, and we consider all functions on the intersection of the corresponding sectors and on the
smallest of the strips. Such an intersection contains a sector of the same kind (with the largest
vertex $\kappa_6$ and the angles $\min\beta_0$, $\max\beta_2$), and Lemma~\ref{lem:S2} applies to
it. In what follows, $S^b$ denotes this common sector.

\paragraph{Taylor coefficients} Let $g=M+tp$ with real $t$, and let
\[
\gamma^b=\sum_{j\ge1}z^j\gamma^b_j(x,\theta;t).
\]
By~\eqref{eq:phi-bound} and Cauchy's inequality $|\gamma^b_j|\le C_j/|x|$ on $S^b$, and
from~\eqref{eq:gamma-eq}
\begin{equation}\label{eq:gamma-j}
x(\partial_x-i\partial_\theta)\gamma_j-2\gamma_j
=x\sum_{i+l=j}\gamma_i\,\partial_x\gamma_l-\sum_{i+l=j}\gamma_i\gamma_l+(M+tp)\,\delta_{j1}.
\end{equation}
By Lemma~\ref{lem:S2}, the decaying solution of equation~\eqref{eq:gamma-j} is unique once the
lower orders are given. Hence, by induction on $j$, $\gamma^b_j$ is a trigonometric polynomial in
$\theta$ and a polynomial in $t$ of degree at most $j$. Indeed, the right-hand side is a
polynomial in $t$. Taking $j+1$ values of $t$ and inverting the Vandermonde matrix, we obtain
decaying solutions for each coefficient of the right-hand side, and their combination coincides
with $\gamma_j(t)$ by uniqueness. The Fourier modes outside the spectrum of the right-hand side
satisfy the homogeneous equation and vanish. Write
\begin{equation}\label{eq:v-w}
\gamma^b_j=v^b_j+t\,w^b_j+O(t^2).
\end{equation}
Then $v^b=\sum_{j\ge1}z^jv^b_j$ is a formal solution of equation~\eqref{eq:gamma-eq} with $g=M$,
and $w^b=\sum_{j\ge1}z^{j-1}w^b_j$ is a formal solution of the linearized equation
\begin{equation}\label{eq:linearized}
x\bigl[(1-v)\,\partial_x-i\partial_\theta\bigr]w=(2-2v+xv_x)\,w+p .
\end{equation}
Set (termwise in $z$)
\[
B^b(x,\theta)=-\int_x^{+\infty}\frac{4\,w^b(x',\theta)}{x'^{2}}\,dx' .
\]
Then $\partial_xB^b=4w^b/x^2$, $B^b=O(x^{-2})$ and
\begin{equation}\label{eq:B-eq}
\bigl[(1-v^b)\,\partial_x-i\partial_\theta\bigr]B^b=-\frac{2p}{x^2} .
\end{equation}
Indeed, the derivative in $x$ of the difference of the two sides vanishes
by~\eqref{eq:linearized}, and both sides tend to zero as $x\to+\infty$.

\begin{lemma}[coefficients of $\chi$]\label{lem:S3}
For every $k\ge1$
\begin{equation}\label{eq:S3}
[t^1]\,[z^k]\,\chi\bigl[z(M+tp)\bigr]
=\lim_{x\to+\infty}\bigl(-i\,e^{x}\bigr)\,\partial_x\bigl[B^u_k-B^s_k\bigr]^{[1]}(x),
\end{equation}
where $B^b_k$ is the coefficient of $z^{k-1}$ in the series $B^b$, and $[\,\cdot\,]^{[1]}$
denotes the first harmonic in $\theta$.
\end{lemma}

\begin{proof}
By~\eqref{eq:phi-def} and~\eqref{eq:gamma-def}, at $\zeta=-ix$ we have
$\cJ=\varphi^u-\varphi^s=-4x^{-2}(\gamma^u-\gamma^s)$. As in the proof of Lemma~\ref{lem:stokes},
let $\chi_x(z)=i\,e^x\,\fc{\cJ}{1}(x;z)$. By item~1 of that lemma, $|\chi_x(z)|\le C_8$ for
$|z|<R_0$ and large $x$, and by item~2, $\chi_x(z)\to\chi[zg]$ pointwise. The functions $\chi_x$
are analytic in $z$, so by Vitali's theorem the convergence is locally uniform and the Taylor
coefficients in $z$ converge. For the profiles $g_t$ with $k+1$ values of $t$ take $x$ larger
than all the corresponding $\kappa_6$. The coefficient $[z^k]\,\chi_x$ is a polynomial in $t$ of
degree at most $k$ which converges for every real $t$. Hence its coefficient of $t$ converges as
well, and it equals
\[
-4i\,e^x\,x^{-2}\,\bigl[w^u_k-w^s_k\bigr]^{[1]}=-i\,e^x\,\partial_x\bigl[B^u_k-B^s_k\bigr]^{[1]}.
\qedhere
\]
\end{proof}

\subsection{The Borel class and Laplace integrals}\label{ss:borel-class}

Denote by $\mathbb A_N$ the class of functions $\hat f$ analytic in $\Omega_N$ and polynomially
bounded on each of the sets
\[
\Sigma_\delta=\{s:\ \delta\le|\arg s|\le\pi-\delta\}\cup\{s:\ |s|\le N/2\},\qquad\delta>0,
\]
that is, $|\hat f(s)|\le C_\delta\,(1+|s|)^{N_\delta}$ for $s\in\Sigma_\delta$.

\begin{lemma}\label{lem:S4}
The class $\mathbb A_N$ is closed under addition, multiplication by polynomials in $s$,
differentiation, integration $\int_0^s$, convolution, division by $(n-s)^m$ for real $n$ with
$|n|\ge N$, and division by $s^m$ if the function has a zero of order at least $m$ at zero.
\end{lemma}

\begin{proof}
The domain $\Omega_N$ is star-shaped with respect to zero, and the convolution
$(\hat f*\hat g)(s)=s\int_0^1\hat f(s\tau)\,\hat g(s(1-\tau))\,d\tau$ is analytic in it. The
points $s\tau$ lie in the same set $\Sigma_\delta$ as $s$, which gives the growth estimate. On
$\Sigma_\delta$ we have $|n-s|\ge\min\{|n|\sin\delta,\ N/2\}$. The rest is obvious. For the
derivative we use Cauchy's inequality with a smaller $\delta$.
\end{proof}

Let $\hat f\in\mathbb A_N$ and $\vartheta\in(0,\pi)$. The integral
$\int_0^{\infty e^{i\vartheta}}e^{-sx}\hat f(s)\,ds$ converges for $\Re(xe^{i\vartheta})>0$ and,
by Cauchy's theorem, does not depend on $\vartheta$ on the intersections of the half-planes. It
defines a function $\cL^+\hat f$ analytic in $\CC\setminus i[0,+\infty)$. Similarly, rays with
$\vartheta\in(-\pi,0)$ define a function $\cL^-\hat f$ analytic in
$\CC\setminus(-i)[0,+\infty)$. The same holds for functions of the form $\hat g(s)/(n-s)$ with
$\hat g\in\mathbb A_N$ and real $n$, $0<|n|<N$, because the rays do not pass through the pole.

\begin{lemma}\label{lem:S5}
Let $\hat f,\hat g\in\mathbb A_N$.
\begin{enumerate}[label=\textup{\arabic*.}]
\item $\cL^\pm(\hat f*\hat g)=\cL^\pm\hat f\cdot\cL^\pm\hat g$;\quad
$\partial_x\cL^\pm\hat f=\cL^\pm(-s\hat f)$;\quad
$x\,\cL^\pm\hat f=\hat f(0)+\cL^\pm(\hat f')$;\quad $\cL^\pm(s^m/m!)=x^{-m-1}$.
\item If $\hat f(s)=O(s^m)$ at zero, then $|\cL^\pm\hat f(x)|\le C|x|^{-m-1}$ for $|x|\ge1$ in
every closed subsector of the domain of definition.
\item As $x\to+\infty$, for any $\delta>0$ we have
$(\cL^+-\cL^-)\hat f(x)=O\bigl(e^{-(N-\delta)x}\bigr)$.
\item If $0<n<N$, then
\begin{equation}\label{eq:residue}
(\cL^+-\cL^-)\Bigl[\frac{\hat g(s)}{n-s}\Bigr](x)=2\pi i\,\hat g(n)\,e^{-nx}
+O\bigl(e^{-(N-\delta)x}\bigr).
\end{equation}
\end{enumerate}
The estimates of items 3 and 4 are preserved under differentiation in $x$. The rules of item~1
(except the convolution) and the estimate of item~2 also hold for transforms of the form
$\hat g(s)/(n-s)$ with $\hat g\in\mathbb A_N$, $0<|n|<N$. For them, item~3 holds when $n<0$, since
the pole lies on the negative half-axis and does not fall between the rays.
\end{lemma}

\begin{proof}
1. Fubini's theorem on the ray and integration by parts. 2. Choose a ray on which
$\Re(xe^{i\vartheta})\ge c|x|$ and estimate $|\hat f(s)|\le C|s|^m(1+|s|)^{N_\delta}$. 3. In the
sector $|\arg s|<\vartheta<\pi/2$ the function $\hat f$ is analytic outside the ray $[N,\infty)$,
so both rays can be translated parallel to themselves to start at the point $N-\delta$. On the
translated rays $|e^{-sx}|\le e^{-(N-\delta)x}e^{-\varrho x\cos\vartheta}$, where $\varrho$ is
the distance to the origin of the ray. 4. The same argument applies, taking the residue into
account. The closed contour ``lower ray outward, arc, upper ray inward'' goes around the pole
counterclockwise, and the residue of the function $e^{-sx}\hat g(s)/(n-s)$ at the point $s=n$
equals $-\hat g(n)e^{-nx}$. Therefore
$\cL^--\cL^+=-2\pi i\,\hat g(n)e^{-nx}+O(e^{-(N-\delta)x})$. Differentiation in $x$ multiplies the
transform by $-s$ and does not lead out of the class.
\end{proof}

\begin{remark}[sign check]\label{rem:sign-check}
In the first order equation~\eqref{eq:B-eq} reads $(\partial_x-i\partial_\theta)B_1=-2p/x^2$, and
the Borel transform of the mode $n$ equals $-2\fc pn\,s/(n-s)$. The unstable branch is the
integral along the characteristic from $t=-\infty$. The substitution
$(x+it)^{-2}=\int_0^\infty s\,e^{-s(x+it)}ds$ gives the pole $s=n-i0$, that is, the ray above
the real axis. Thus the unstable branch corresponds to $\cL^+$, and the stable one to $\cL^-$.
By~\eqref{eq:residue}, $[B_1^u-B_1^s]^{[1]}=-4\pi i\,\fc p1e^{-x}$, and Lemma~\ref{lem:S3} gives
$[t^1]\,[z^1]\,\chi=4\pi\fc p1$, in agreement with Proposition~\ref{prop:normalization}.
\end{remark}

\subsection{The Newtonian form of the inner equation}\label{ss:newton}

Below, all series in $z$ are formal, and in each order only finitely many terms and Fourier modes
are involved. The composition $f(\xi-Y,\theta)$ with a series $Y=O(z)$ is understood as
$\sum_m(-Y)^m\partial_\xi^mf/m!$. The chain rule and the inversion of the change of variables for
such compositions are identities between coefficients, since for small $|z|$ and truncated series
these are ordinary compositions of analytic functions.

\paragraph{Step 1: analyticity} By Lemma~\ref{lem:S4} and induction on $j$, all the functions
$\hY^{[n]}_{(j)}$ and $\hE^{[n]}_{(j)}$ belong to $\mathbb A_N$. For the zero mode we use that
$\hQ^{[0]}_{(j)}(t)=O(t)$, whence $\int_0^st^2\hQ^{[0]}_{(j)}\,dt=O(s^4)$. Moreover,
$\hE^{[n]}_{(j)}(s)=O(s^2)$ for $j\ge1$. The first statement of Proposition~\ref{thm:S} is proved.

\paragraph{Step 2: the equation} Set
\begin{equation}\label{eq:Ypm}
Y^\pm=\sum_{j\ge1}z^j\sum_ne^{in\theta}\,\cL^\pm\hY^{[n]}_{(j)}(\xi),\qquad
\sD=\partial_\xi-i\partial_\theta .
\end{equation}
Then
\begin{equation}\label{eq:Y-eq}
\xi\,\sD^2Y^\pm-2\,\sD Y^\pm=zM+Y^\pm\,\sD^2Y^\pm-(\sD Y^\pm)^2,
\end{equation}
which is equivalent to equation~\eqref{eq:newton} for $x=\xi-Y^\pm$:
\begin{equation}\label{eq:newton-pm}
(\sD x)^2-x\,\sD^2x=1+zM .
\end{equation}
Indeed, differentiating the definition \eqref{eq:Y1}, \eqref{eq:Yj}, we get
$\partial_s\bigl[(n-s)^4\hY^{[n]}\bigr]=(n-s)^2\hQ^{[n]}$, which is equivalent to the equality
\[
\partial_s\bigl[(n-s)^2\hY^{[n]}\bigr]-2(n-s)\,\hY^{[n]}=\hQ^{[n]},
\]
and $n^2\hY^{[n]}(0)=z\fc Mn$. Here $\hY^{[n]}=\sum_jz^j\hY^{[n]}_{(j)}$, and similarly for
$\hQ^{[n]}$. On the mode $n$ the operator $\sD$ corresponds to multiplication of the transform by
$n-s$, and multiplication by $\xi$ corresponds to the derivative $\partial_s$ plus the value at
zero (Lemma~\ref{lem:S5}(1)). Therefore
$\xi\,\sD^2Y^{[n]}-2\,\sD Y^{[n]}=z\fc Mn+\cL^\pm\hQ^{[n]}$. Finally, the bilinear form
$\tfrac12(f\,\sD^2g+g\,\sD^2f)-\sD f\,\sD g$ on the modes $n_1$, $n_2$ corresponds to the
convolution with the kernel
\[
\tfrac12(n_2-t+u)^2+\tfrac12(n_1-u)^2-(n_1-u)(n_2-t+u)=\tfrac12\bigl[(n_1-u)-(n_2-t+u)\bigr]^2,
\]
so $\cL^\pm\hQ^{[n]}$ is the mode $n$ of the function $Y\,\sD^2Y-(\sD Y)^2$.

\paragraph{Step 3: the background} Let $\mathcal Z^\pm(x,\theta)$ be the formal solution of the
equation $\mathcal Z=Y^\pm(x+\mathcal Z,\theta)$. Then $\xi=x+\mathcal Z^\pm(x,\theta)$ inverts
the change of variables $x=\xi-Y^\pm(\xi,\theta)$. Set
\begin{equation}\label{eq:v-tilde}
\tilde v^\pm(x,\theta)=(\sD Y^\pm)\bigl(x+\mathcal Z^\pm(x,\theta),\theta\bigr).
\end{equation}
The coefficients $\mathcal Z^\pm_j$, $\tilde v^\pm_j$ are polynomials in the functions
$\partial_\xi^m\cL^\pm\hY^{[n]}_{(i)}$, $i\le j$. They are analytic in
$\CC\setminus(\pm i)[0,\infty)$ and, by Lemma~\ref{lem:S5}(2), are $O(1/|x|)$ in closed
subsectors. In the variables $(\xi,\theta)$ we have
\[
\sD x=1-\tilde v\circ x,\qquad\sD^2x=-\bigl[(1-\tilde v)\,\tilde v_x-i\tilde v_\theta\bigr]\circ x,
\]
therefore equation~\eqref{eq:newton-pm} is equivalent to the equation
\[
x\bigl[(1-\tilde v)\,\tilde v_x-i\tilde v_\theta\bigr]=2\tilde v-\tilde v^2+zM,
\]
that is, to equation~\eqref{eq:gamma-eq} with $g=M$. The sector $S^u$ lies in the domain
$\CC\setminus i[0,\infty)$ together with a closed subsector, and the sector $S^s$ lies in the
domain $\CC\setminus(-i)[0,\infty)$. By Lemma~\ref{lem:S2} and induction on $j$
\begin{equation}\label{eq:v-identification}
v^u_j=\tilde v^+_j\ \ \text{on }S^u,\qquad v^s_j=\tilde v^-_j\ \ \text{on }S^s .
\end{equation}

\begin{remark}[origin of equation~\eqref{eq:newton}]\label{rem:newton-origin}
The inner equation is the Hamilton--Jacobi equation
$\partial_\theta\psi+\mathcal H(\zeta,\partial_\zeta\psi,\theta)=0$ with Hamiltonian
$\mathcal H=\Pi-\tfrac18\zeta^2\Pi^2+2zg(\theta)/\zeta^2$. Along its characteristics
$\dot\theta=1$, $\dot\zeta=1-\tfrac14\zeta^2\Pi$,
$\dot\Pi=\tfrac14\zeta\Pi^2+4zg/\zeta^3$, which implies the exact identity
\[
\dot\zeta^2-\zeta\,\ddot\zeta=1+zg(\theta).
\]
This is the inner limit of the pendulum equation. For $\zeta=\sqrt2\,e^{-u/2}$ it becomes the
equation $\ddot u=(1+zg(\theta))\,e^u$, which is the reference system
of~\cite[Section~4]{Gelfreich97}. If one introduces a coordinate that straightens
the flow of the characteristics on the graph $\Pi=\partial_\zeta\psi$ and passes to the variables
$\zeta=-ix$, the identity takes the form~\eqref{eq:newton}. In terms of~\cite{Gelfreich97},
equation~\eqref{eq:newton} is equation~(4.3) of that paper with $\omega=1$ and $g=zM$, written
for the function $x=i\sqrt2\,e^{-u/2}$ of the variables $\xi=is_1$ and $\theta=s_2$. The proof
does not use this argument.
We construct the functions explicitly, and the coincidence with the actual solutions is given by
equality~\eqref{eq:v-identification}.
\end{remark}

\subsection{Computation of the jump: proof of Theorem~\ref{thm:series}}\label{ss:jump}

\paragraph{Step 4: the response} Set
\[
\mathcal E^\pm=(\xi-Y^\pm)^{-2}=\sum_{r\ge0}(r+1)\,\xi^{-2-r}\,(Y^\pm)^r .
\]
By Lemma~\ref{lem:S5}(1) and definition~\eqref{eq:E}, $\mathcal E^\pm=\cL^\pm\hE$ mode by mode and
order by order. Let
\begin{equation}\label{eq:S-hat}
\hat S^{[n]}=2\sum_m\fc pm\,\hE^{[n-m]},\qquad
\tilde B^{\pm}=\sum_ne^{in\theta}\,\cL^\pm\Bigl[\frac{-\hat S^{[n]}(s)}{n-s}\Bigr](\xi).
\end{equation}
For $n=0$ the transform equals $\hat S^{[0]}(s)/s$. This is legitimate, since
$\hat S^{[0]}(s)=O(s^2)$. Indeed, the mean of $p$ is zero, and $\hE^{[n]}_{(j)}(s)=O(s^2)$ for
$(n,j)\ne(0,0)$. By Lemma~\ref{lem:S5}(1), $\sD\tilde B^\pm=-2p\,\mathcal E^\pm$, and
$\tilde B^\pm=O(|\xi|^{-2})$, since all the transforms are $O(s)$ at zero. Set
$B^\pm(x,\theta)=\tilde B^\pm(x+\mathcal Z^\pm,\theta)$. Then
$\sD\tilde B=\bigl[(1-\tilde v)\,\partial_xB-i\partial_\theta B\bigr]\circ x$ and
$\mathcal E=x^{-2}$, that is,
\[
\bigl[(1-\tilde v^\pm)\,\partial_x-i\partial_\theta\bigr]B^\pm=-\frac{2p}{x^2} .
\]
Differentiating in $x$, we find that the function $\tfrac14x^2\partial_xB^\pm$ satisfies
equation~\eqref{eq:linearized} with $v=\tilde v^\pm$ and decays as $O(1/|x|)$. The homogeneous
part of equation~\eqref{eq:linearized} in each order is the same operator
$x(\partial_x-i\partial_\theta)-2$. Therefore, by Lemma~\ref{lem:S2} and
equality~\eqref{eq:v-identification}, $\tfrac14x^2\partial_xB^+=w^u$ on $S^u$ and
$\tfrac14x^2\partial_xB^-=w^s$ on $S^s$. Since $B^\pm\to0$ as $x\to+\infty$, for large real $x$
\begin{equation}\label{eq:B-identification}
B^u=B^+,\qquad B^s=B^- .
\end{equation}

\paragraph{Step 5: the jump} Let $x\to+\infty$. By Lemma~\ref{lem:S5}(3), all the differences
$\partial^m_\xi(Y^+_{(j)}-Y^-_{(j)})$, and with them $\mathcal Z^+-\mathcal Z^-$, are of order
$O(e^{-(N-\delta)x})$. By items 3 and 4 of the same lemma
\begin{equation}\label{eq:B-tilde-jump}
\tilde B^+-\tilde B^-=-2\pi i\,\hat S^{[1]}(1)\,e^{-\xi+i\theta}+O\bigl(e^{-(2-\delta)\xi}\bigr).
\end{equation}
Indeed, the mode $n$ with $0<n<N$ gives $-2\pi i\,\hat S^{[n]}(n)\,e^{-n\xi}$, and the other
modes give $O(e^{-(N-\delta)\xi})$. The same holds after differentiation. Further,
\[
B^+-B^-=\sum_m\frac{(\mathcal Z^+)^m}{m!}\,\partial_\xi^m\bigl(\tilde B^+-\tilde B^-\bigr)
+\sum_m\frac{(\mathcal Z^+)^m-(\mathcal Z^-)^m}{m!}\,\partial_\xi^m\tilde B^- .
\]
The second sum is $O(e^{-(N-\delta)x})$. In the first sum the term with $m=0$ equals
$\tilde B^+-\tilde B^-$ at the point $(x,\theta)$, and the terms with $m\ge1$ contain the factor
$\mathcal Z^+=O(z/x)$ and are therefore $O(e^{-x}/x)$ in each order. Hence, for the coefficient of
$z^{k-1}$,
\[
\bigl[B^u_k-B^s_k\bigr]^{[1]}(x)=-2\pi i\,\hat S^{[1]}_{(k-1)}(1)\,e^{-x}+O\bigl(e^{-x}/x\bigr),
\]
and, after differentiation,
$\partial_x\bigl[B^u_k-B^s_k\bigr]^{[1]}(x)=2\pi i\,\hat S^{[1]}_{(k-1)}(1)\,e^{-x}+O(e^{-x}/x)$.
By Lemma~\ref{lem:S3}
\begin{align*}
[t^1]\,[z^k]\,\chi\bigl[z(M+tp)\bigr]
&=\lim_{x\to+\infty}\bigl(-i\,e^{x}\bigr)\cdot2\pi i\,\hat S^{[1]}_{(k-1)}(1)\,e^{-x}\\
&=2\pi\,\hat S^{[1]}_{(k-1)}(1)=4\pi\sum_m\fc pm\,\hE^{[1-m]}_{(k-1)}(1).
\end{align*}
Proposition~\ref{thm:S} is proved.

\begin{proof}[Proof of items (1) and (2) of Theorem~\ref{thm:series}]
Let $M$ have period $\pi$, and let $p$ contain only odd harmonics. Then Proposition~\ref{thm:S}
applies with $N=2$, and the analyticity of the functions $\hY^{[n]}_{(j)}$, $\hE^{[n]}_{(j)}$ in
$\Omega_2$ is proved. By Lemma~\ref{lem:variation-limit} the entire function $G_{M,p}$ is
defined, and by~\eqref{eq:G-taylor} and~\eqref{eq:S1}
\[
[z^{k-1}]\,G_{M,p}=4\pi\sum_m\fc pm\,\hE^{[1-m]}_{(k-1)}(1),\qquad k\ge1 .
\]
The Taylor series of an entire function converges everywhere. Consequently, the
series~\eqref{eq:G} converges for all $\mu\in\CC$, and its sum coincides with the function of
Lemma~\ref{lem:variation-limit}. This is item~(1). Item~(2) is equality~\eqref{eq:G-taylor}
together with Corollary~\ref{cor:chi-entire}.
\end{proof}

\subsection{Nonnegativity and consequences}\label{ss:positivity}

\begin{lemma}\label{lem:S6}
If $\fc Mn\ge0$ for all $n$, then $\hY^{[n]}_{(j)}(s)\ge0$ and $\hE^{[n]}_{(j)}(s)\ge0$ for
$0\le s<N$.
\end{lemma}

\begin{proof}
We argue by induction on $j$. For $n\in N\ZZ\setminus\{0\}$ the factors $(n-t)^2$ and
$(n-s)^{-4}$ in~\eqref{eq:Yj} are positive on $[0,N)$, and for $n=0$ these are $t^2$ and
$s^{-4}$. The kernel in~\eqref{eq:Qj} is a perfect square, and the convolution of nonnegative
functions is nonnegative. In the first order $\hY^{[n]}_{(1)}=n^2\fc Mn(n-s)^{-4}\ge0$.
Formula~\eqref{eq:E} involves convolutions of powers of $s$ and of the functions
$\hY^{[n]}_{(j)}$ with positive coefficients.
\end{proof}

The function $\hE^{[n]}_{(j)}$ is homogeneous of degree $j$ in $M$. Therefore, if
$\mu\fc Mn\ge0$ for all $n$, then $\mu^j\hE^{[n]}_{(j)}[M](s)=\hE^{[n]}_{(j)}[\mu M](s)\ge0$ for
$0\le s<N$.

\begin{proof}[Proof of item (3) of Theorem~\ref{thm:series}]
The first statement is Lemma~\ref{lem:S6} with $N=2$. If $\mu\fc Mn\ge0$ and $\fc pk\ge0$ for all
$n$, $k$, then by what was said above all terms of the series~\eqref{eq:G} are nonnegative, and
the term with $j=0$ equals $4\pi\fc p1$.
\end{proof}

\begin{proof}[Proof of Theorem~\ref{thm:example}]
Let $M=\cos4\theta$ and $p=\cos\theta$. The modes of the functions $\hE$ lie in $4\ZZ$, so
in~\eqref{eq:G} only the term with $k=1$ remains, that is,
$\sum_k\fc pk\,\hE^{[1-k]}_{(j)}=\tfrac12\hE^{[0]}_{(j)}$. The mode of a term of order $j$ is
congruent to $4j$ modulo $8$, so $\hE^{[0]}_{(j)}=0$ for odd $j$. Moreover,
$\hE^{[0]}_{(0)}(1)=1$. This gives~\eqref{eq:G-example} with $T_j=\hE^{[0]}_{(2j)}(1)$, and
$T_j\ge0$ by Lemma~\ref{lem:S6}. Consequently, $G(\mu)\ge2\pi(1+T_1\mu^2)\ge2\pi$ for real $\mu$.
The asymptotics~\eqref{eq:example} now follows from Theorem~\ref{thm:main} with $A_f=-1$,
$p_f=\cos\theta$. Since $G\ge2\pi$, the remainder $o(1)$ can be written as a factor $1+o(1)$.

The number $T_1$ is given by explicit integrals. From \eqref{eq:Y1}--\eqref{eq:E}
\begin{gather*}
\hY^{[\pm4]}_{(1)}(s)=\frac8{(4\mp s)^4},\qquad
\hQ^{[0]}_{(2)}(t)=64\int_0^t\frac{(8+t-2u)^2\,du}{(4-u)^4(4+t-u)^4},\\
\hY^{[0]}_{(2)}(s)=\frac1{s^{4}}\int_0^st^2\,\hQ^{[0]}_{(2)}(t)\,dt,\\
T_1=\int_0^1(1-t)^2\,\hY^{[0]}_{(2)}(t)\,dt
+\int_0^1(1-t)^3\,\bigl(\hY^{[4]}_{(1)}*\hY^{[-4]}_{(1)}\bigr)(t)\,dt=0.0053202218704\ldots
\qedhere
\end{gather*}
\end{proof}

\begin{proof}[Proof of Proposition~\ref{prop:nondeg}]
(i) The function $G_{M,p_f}$ is entire by Theorem~\ref{thm:series}, and
$G_{M,p_f}(0)=4\pi\fc{p_f}{1}$ by Corollary~\ref{cor:G-zero}. If this number is nonzero, then the
zeros of the entire function $G_{M,p_f}$ are isolated and do not lie near the origin. (ii) This
is item~(3) of Theorem~\ref{thm:series}. (iii) The equality $A_f=-\sum_jj^2\fc f{j,0}$ shows that
$A_f=0$ if $f$ has no harmonics with $\ell=0$. For $f=\cos\alpha_1+\cos3\alpha_2$ and
$m(\theta)=\cos4\theta$ we have $M=\cos8\theta$ and $p_f=\cos3\theta$. The modes of the functions
$\hE$ lie in $8\ZZ$, while the sum in~\eqref{eq:G} involves only $\hE^{[-2]}$ and $\hE^{[4]}$.
Therefore $G_{M,p_f}\equiv0$.
\end{proof}

\begin{remark}[first order in the background amplitude]\label{rem:first-order-mu}
From~\eqref{eq:Y1} and~\eqref{eq:E} we get
$\hE^{[n]}_{(1)}(1)=n^2\fc Mn\int_0^1(1-t)^2(n-t)^{-4}dt=\fc Mn/\bigl(3n(n-1)\bigr)$, and therefore
\[
G_{M,p}(\mu)=4\pi\,\fc p1-\frac{4\pi\mu}3\sum_{k\ne0,1}\frac{\fc pk\,\fc M{1-k}}{k(1-k)}+O(\mu^2).
\]
For example, $G_B(\mu)=2\pi\bigl(1+\tfrac{\mu}{12}+O(\mu^2)\bigr)$. Thus for the background
$\cos2\theta$ the linear term is nonzero, and for small negative $\mu$ the function $G_B$ is
smaller than $2\pi$.
\end{remark}
\section{The background \texorpdfstring{$\cos2\theta$}{cos 2theta} with negative amplitude}\label{sec:negative}

In this section we prove the second part of Theorem~\ref{thm:cos2}, namely that if $\mu<0$, then
\begin{equation}\label{eq:GB-negative}
\frac{G_B(\mu)}{2\pi}\ \ge\ \varkappa_0\,\frac{\sinh\lambda}{\lambda},\qquad
\lambda=\frac{|\mu|}{\sqrt8},\qquad\varkappa_0=\frac{6142921}{102400000}.
\end{equation}
Throughout what follows $\mu<0$ is fixed, and all functions are built from the \emph{positive}
background $|\mu|\cos2\theta$. That is, in formulas \eqref{eq:Y1}--\eqref{eq:E} we take
$M=\cos2\theta$ (the frequency lattice step is two), so by Lemma~\ref{lem:S6} the functions
$\hY^{[n]}_{(j)}$, $\hE^{[n]}_{(j)}$ are nonnegative on $[0,2)$. Note that
$|\mu|=2\sqrt2\,\lambda$.

The plan is as follows. The number $G_B(\mu)/2\pi$ equals the difference $E_0(1)-E_2(1)$ of two
positive quantities. We bound the first from below by comparison with an explicitly solvable
scalar equation. To bound the second from above, we write the system of equations in the Borel
plane as a fixed-point problem for a monotone operator. Its solution is the limit of increasing
iterates and therefore does not exceed any ``supersolution''. After the factor $e^{\lambda s}$ is
extracted, the system is majorized by a system of inequalities that does not contain the
amplitude, and a supersolution of this system is given explicitly by polynomials.

\paragraph{Conventions} We consider families $A=(A_k)_{k\in2\ZZ}$ of measurable functions on
$(0,1]$ with values in $[0,+\infty]$. The product of families is the convolution in $s$ and in the
mode index, that is, $(A*B)_k=\sum_{k_1+k_2=k}A_{k_1}*B_{k_2}$. We denote by $\delta$ the unit of
this multiplication, $A^{*0}=\delta$, and
\[
\exp_*A=\sum_{r\ge0}\frac{A^{*r}}{r!}.
\]
A function of one variable $s$ is identified with the family whose only nonzero mode is $k=0$. In
particular, the expressions $s*A$ and $1*A$ contain the functions $s\mapsto s$ and $s\mapsto1$.
Inequalities between families are understood mode by mode. Since all terms are nonnegative, sums
and integrals can be interchanged.

\subsection{Reduction to the difference of two modes}\label{ss:two-modes}

\paragraph{The functions $\hat H$} Along with $\hY$ and $\hE$ consider the series
\begin{equation}\label{eq:H-hat}
\hat H=2\sum_{\ell\ge1}\frac{(1*\hY)^{*\ell}}{\ell},\qquad
\hY=\sum_{j\ge1}z^j\sum_n\hY^{[n]}_{(j)}e^{in\theta},
\end{equation}
and denote by $\hat H^{[n]}_{(j)}$ its coefficient of $z^je^{in\theta}$. This is the Borel
transform of the series $\log(\xi^2/x^2)=-2\log(1-Y/\xi)$. The functions $\hat H^{[n]}_{(j)}$
belong to the class $\mathbb A_2$, are nonnegative on $[0,2)$ and are $O(s)$ at zero.

\begin{lemma}\label{lem:H-identities}
For all $n\in2\ZZ$ and $j\ge1$
\begin{equation}\label{eq:EH-orders}
\hE^{[n]}_{(j)}=\bigl(s*\exp_*\hat H\bigr)^{[n]}_{(j)},\qquad
(n-s)^2\,\hat H^{[n]}_{(j)}=2\,\hE^{[n]}_{(j)}+\hE^{[n-2]}_{(j-1)}+\hE^{[n+2]}_{(j-1)} .
\end{equation}
\end{lemma}

\begin{proof}
Let $Y=Y^+$ be the series~\eqref{eq:Ypm}, $x=\xi-Y$ and
\[
\mathsf h=\sum_{j\ge1}z^j\sum_ne^{in\theta}\,\cL^+\hat H^{[n]}_{(j)} .
\]
By Lemma~\ref{lem:S5}(1), $\mathsf h=2\sum_\ell(Y/\xi)^\ell/\ell=\log(\xi^2/x^2)$ as a series in
$z$. Therefore $\cL^+(s*\exp_*\hat H)=\xi^{-2}e^{\mathsf h}=x^{-2}=\cL^+\hE$. Further,
from~\eqref{eq:newton-pm}
\[
\sD^2\mathsf h=-2\,\sD\Bigl(\frac{\sD x}{x}\Bigr)-\frac2{\xi^2}
=\frac{2\bigl[(\sD x)^2-x\,\sD^2x\bigr]}{x^2}-\frac2{\xi^2}=\frac{2(1+zM)}{x^2}-\frac2{\xi^2},
\]
and the operator $\sD$ on the mode $n$ corresponds to multiplication of the transform by $n-s$.
For $M=\cos2\theta$ we find that the Laplace transforms of both sides of the second
equality~\eqref{eq:EH-orders} coincide. It remains to note that a function $\hat f$ of class
$\mathbb A_2$ with zero transform $\cL^+\hat f$ vanishes. Indeed, by Watson's lemma
$\cL^+\hat f(\xi)\sim\sum_m\hat f^{(m)}(0)\,\xi^{-m-1}$ as $\xi\to+\infty$, so all derivatives of
$\hat f$ at zero vanish, and the domain $\Omega_2$ is connected.
\end{proof}

\paragraph{Sums over the orders} For $s\in[0,1]$ set
\begin{equation}\label{eq:EH-sums}
E_k(s)=\sum_{j\ge0}|\mu|^j\,\hE^{[k]}_{(j)}(s),\qquad
H_k(s)=\sum_{j\ge1}|\mu|^j\,\hat H^{[k]}_{(j)}(s),\qquad k\in2\ZZ .
\end{equation}
These are sums of nonnegative terms, and until the contrary is proved, they may equal $+\infty$.
For a family $A$ set $\mathcal E[A]=s*\exp_*A$ and define the operator $\Phi$ by
\begin{equation}\label{eq:Phi}
\Phi[A]_k=\frac{2\bigl(\mathcal E[A]_k-s\,\mathbf 1_{k=0}\bigr)
+2\sqrt2\,\lambda\,\bigl(\mathcal E[A]_{k-2}+\mathcal E[A]_{k+2}\bigr)}{(k-s)^2},\qquad s\in(0,1].
\end{equation}
Here $\mathcal E[A]_0-s=\sum_{r\ge1}(s*A^{*r})_0/r!\ge0$. The operator $\Phi$ is monotone, that
is, if $A\le A'$, then $\Phi[A]\le\Phi[A']$.

\begin{lemma}\label{lem:two-modes}
\begin{enumerate}[label=\textup{(\alph*)}]
\item $E=\mathcal E[H]$ and $H=\Phi[H]$.
\item The numbers $E_0(1)$ and $E_2(1)$ are finite, and
\begin{equation}\label{eq:C3}
\frac{G_B(\mu)}{2\pi}=E_0(1)-E_2(1).
\end{equation}
\end{enumerate}
\end{lemma}

\begin{proof}
(a) Multiply equalities~\eqref{eq:EH-orders} by $|\mu|^j$ and sum over $j$. Since all terms are
nonnegative, the order of summation does not matter. In the second equality we take into account
that $\hE^{[k]}_{(0)}=s\,\mathbf 1_{k=0}$ and $|\mu|=2\sqrt2\,\lambda$.

(b) By Theorem~\ref{thm:series} the series
$G_B(|\mu|)/2\pi=\sum_j|\mu|^j\bigl[\hE^{[0]}_{(j)}(1)+\hE^{[2]}_{(j)}(1)\bigr]$ converges, and its
terms are nonnegative. A term of order $j$ contains only modes congruent to $2j$ modulo $4$, so
$\hE^{[0]}_{(j)}=0$ for odd $j$ and $\hE^{[2]}_{(j)}=0$ for even $j$. Consequently,
$E_0(1)+E_2(1)=G_B(|\mu|)/2\pi<\infty$, and replacing $|\mu|$ by $\mu=-|\mu|$ changes the sign of
the second term.
\end{proof}

The next lemma shows that $H$ is the smallest nonnegative solution of the equation $A=\Phi[A]$.

\begin{lemma}[minimality]\label{lem:minimal}
The sequence $\Phi^n[0]$ is nondecreasing and converges to $H$ pointwise on $(0,1]$.
Consequently, if $0\le A\le H$, then $\Phi^n[A]\to H$.
\end{lemma}

\begin{proof}
Since $\Phi[0]\ge0$ and the operator $\Phi$ is monotone, the sequence $\Phi^n[0]$ is
nondecreasing. From $0\le H$ and $H=\Phi[H]$ we get $\Phi^n[0]\le H$ by induction.

Let us prove the reverse inequality for the limit. Each iterate is a series in powers of $|\mu|$
with nonnegative coefficients, which we write as
$\Phi^n[0]=\sum_{j\ge1}|\mu|^j\,\mathsf h^{(n)}_j$. Separating in~\eqref{eq:Phi} the term linear
in $A$ of the same order, we obtain
\begin{equation}\label{eq:order-recursion}
\mathsf h^{(n+1)}_j=\mathsf L\,\mathsf h^{(n)}_j+\mathsf G_j\bigl(\mathsf h^{(n)}_1,\dots,\mathsf h^{(n)}_{j-1}\bigr),
\qquad(\mathsf L\,g)_k=\frac{2\,(s*g_k)}{(k-s)^2},
\end{equation}
where $\mathsf G_j$ is a finite sum of convolutions of its arguments with nonnegative
coefficients. By Lemma~\ref{lem:H-identities} the same relations connect the functions
$\hat H_{(j)}$, namely
$\hat H_{(j)}=\mathsf L\,\hat H_{(j)}+\mathsf G_j(\hat H_{(1)},\dots,\hat H_{(j-1)})$. The
operator $\mathsf L$ is of Volterra type. More precisely, if $0\le g\le C$ and $g_0(s)\le Cs$,
then $(\mathsf L^mg)_k\le C\,(2s^2)^m/(2m)!$ for $k\ne0$ (since $(k-s)^2\ge1$) and
$(\mathsf L^mg)_0\le Cs/3^m$. Therefore
$\hat H_{(j)}=\sum_{m\ge0}\mathsf L^m\mathsf G_j(\hat H_{(1)},\dots,\hat H_{(j-1)})$.

We show by induction on $j$ that $\mathsf h^{(n)}_j\uparrow\hat H_{(j)}$ as $n\to\infty$.
From~\eqref{eq:order-recursion} and $\mathsf h^{(0)}=0$
\[
\mathsf h^{(n)}_j=\sum_{m=0}^{n-1}\mathsf L^m\,\mathsf G_j\bigl(\mathsf h^{(n-1-m)}_1,\dots,\mathsf h^{(n-1-m)}_{j-1}\bigr).
\]
By the induction hypothesis, for fixed $m$ the term increases to
\[
\mathsf L^m\,\mathsf G_j\bigl(\hat H_{(1)},\dots,\hat H_{(j-1)}\bigr),
\]
and by the monotone convergence theorem the sum tends to $\hat H_{(j)}$. Applying the same theorem
to the series in $j$, we obtain
\[
\lim_{n\to\infty}\Phi^n[0]=\sum_{j\ge1}|\mu|^j\,\hat H_{(j)}=H .
\]
The second statement of the lemma follows from the inequalities
$\Phi^n[0]\le\Phi^n[A]\le\Phi^n[H]=H$.
\end{proof}

\subsection{Lower bound for the zero mode}\label{ss:zero-mode}

Set
\begin{equation}\label{eq:Xb-Zb}
X_b(s)=\frac{2\,(\cosh\lambda s-1)}{s},\qquad Z_b(s)=\frac{\sinh\lambda s}{\lambda}.
\end{equation}
These functions are related by the equalities
\begin{equation}\label{eq:scalar-eq}
Z_b=s*\exp_*X_b,\qquad s^2X_b=2\,(Z_b-s)+4\lambda^2\,Z_b*Z_b .
\end{equation}
Indeed, their Laplace transforms are $\cL Z_b=(\xi^2-\lambda^2)^{-1}=:z_b$ and
$\cL X_b=-\log(1-\lambda^2/\xi^2)=\log(\xi^2z_b)$. The first equality~\eqref{eq:scalar-eq} means
that $z_b=\xi^{-2}\exp(\log\xi^2z_b)$, and the second, after the Laplace transform, takes the form
$(\log z_b)''=2z_b+4\lambda^2z_b^2$ and is checked by differentiation.

\begin{lemma}\label{lem:zero-mode}
On $(0,1]$
\begin{equation}\label{eq:zero-mode}
H_0\ge X_b,\qquad E_0\ge Z_b,\qquad H_{\pm2}\ge\frac{2\sqrt2\,\lambda\,Z_b}{(2\mp s)^2}.
\end{equation}
\end{lemma}

\begin{proof}
Let $X=H_0$ and $Z=s*\exp_*X$. Keeping in the sum $\mathcal E[H]=s*\exp_*H$ only the terms in
which all modes are zero, or one mode equals $\pm2$ and the others are zero, we obtain
\begin{equation}\label{eq:E-lower}
E_0\ge Z,\qquad E_{\pm2}\ge Z*H_{\pm2}.
\end{equation}
By Lemma~\ref{lem:two-modes}(a),
$(2\mp s)^2H_{\pm2}\ge2\sqrt2\,\lambda\,E_0\ge2\sqrt2\,\lambda\,Z$, whence
$H_2+H_{-2}\ge2\sqrt2\,\lambda\,\bigl[(2-s)^{-2}+(2+s)^{-2}\bigr]Z\ge\sqrt2\,\lambda\,Z$. Now the
zero mode of the equality $H=\Phi[H]$ gives
\begin{equation}\label{eq:X-ineq}
s^2X=2\,(E_0-s)+2\sqrt2\,\lambda\,(E_2+E_{-2})\ \ge\ 2\,(Z-s)+4\lambda^2\,Z*Z .
\end{equation}
This inequality holds coefficientwise in $|\mu|$, since in its derivation nonnegative terms were
discarded in each order. Let $X=\sum_j|\mu|^jX_j$ and $X_b=\sum_j|\mu|^jX_{b,j}$. In order $j$
inequality~\eqref{eq:X-ineq} has the form $X_j-\mathsf L_0X_j\ge\mathsf g_j(X_1,\dots,X_{j-1})$,
where $\mathsf L_0g=2s^{-2}(s*g)$ and $\mathsf g_j$ is nondecreasing in each argument. For $X_b$,
by~\eqref{eq:scalar-eq}, the equality
$X_{b,j}-\mathsf L_0X_{b,j}=\mathsf g_j(X_{b,1},\dots,X_{b,j-1})$ holds. In the norm
$\sup_{0<s\le1}|g(s)|/s$ we have $\|\mathsf L_0\|\le\tfrac13$, so the operator
$(1-\mathsf L_0)^{-1}=\sum_m\mathsf L_0^m$ preserves nonnegativity. Induction on $j$ gives
$X_j\ge X_{b,j}$, that is, $X\ge X_b$. Then $E_0\ge Z\ge s*\exp_*X_b=Z_b$, and the third
inequality~\eqref{eq:zero-mode} follows from $(2\mp s)^2H_{\pm2}\ge2\sqrt2\,\lambda\,E_0$.
\end{proof}

\subsection{Equations for the remainders}\label{ss:remainders}

Set
\begin{equation}\label{eq:f-ell}
f(s)=\frac{1-e^{-2\lambda s}}2=\lambda\,e^{-\lambda s}Z_b(s),\qquad
\ell_\pm(s)=\frac{2\sqrt2}{(2\mp s)^2},
\end{equation}
and denote by $H^{(0)}$ the family with $H^{(0)}_0=X_b$,
$H^{(0)}_{\pm2}=e^{\lambda s}\ell_\pm f=2\sqrt2\,\lambda\,Z_b/(2\mp s)^2$ and zero other modes.
By Lemma~\ref{lem:zero-mode}, $H\ge H^{(0)}$.

A family $A$ is called \emph{admissible} if $A\ge H^{(0)}$ and
$\sup_{0<s\le1}\sum_kA_k(s)<\infty$. With an admissible family we associate the
\emph{remainders}
\begin{equation}\label{eq:remainders}
\begin{gathered}
W_k=e^{-\lambda s}\bigl(A_k-X_b\,\mathbf 1_{k=0}\bigr),\qquad U=\exp_*W-\delta,\qquad N=U-W,\\
S=W_0,\qquad D_{\pm2}=W_{\pm2}-\ell_\pm f,\qquad D_k=W_k\ \ (k\ne0,\pm2),\qquad D_0=0 .
\end{gathered}
\end{equation}
All these functions are nonnegative and bounded in the sum over the modes. Multiplication by $s$ is
a derivation of the convolution algebra, so $s\exp_*W=(sW)*\exp_*W$, that is,
\begin{equation}\label{eq:N-eq}
sN=(sW)*U=(sW)*(W+N).
\end{equation}
We also introduce the linear expression
\begin{equation}\label{eq:T-lambda}
\mathcal T_\lambda[X,V]=\frac f\lambda*X+\sqrt2\,f*V
\end{equation}
and the function
\begin{equation}\label{eq:d-tilde}
\tilde d(s)=\tfrac12\bigl(\ell_+(s)+\ell_-(s)-\sqrt2\bigr)=\frac{\sqrt2\,s^2(12-s^2)}{2\,(4-s^2)^2}\ \ge0 .
\end{equation}

\begin{lemma}\label{lem:transform}
Let a family $A$ be admissible, let $W$, $U$, $N$, $D$, $S$ be its remainders and $A'=\Phi[A]$.
Then
\begin{equation}\label{eq:C8}
\lambda\,e^{-\lambda s}\,\mathcal E[A]_k=f\,\mathbf 1_{k=0}+f*U_k,
\end{equation}
$A'\ge H^{(0)}$, and the remainders $D'$, $S'$ of the family $A'$ are
\begin{equation}\label{eq:C9}
\begin{aligned}
D'_k&=\frac2{(k-s)^2}\,\mathcal T_\lambda\bigl[U_k,\ U_{k-2}+U_{k+2}\bigr]\qquad(k\ne0),\\
s^2S'&=2\,\mathcal T_\lambda\bigl[U_0,\ 2\tilde d\,f+D_2+D_{-2}+N_2+N_{-2}\bigr].
\end{aligned}
\end{equation}
\end{lemma}

\begin{proof}
Multiplication by $e^{\lambda s}$ is an automorphism of the convolution algebra, so
$\exp_*A=\exp_*(X_b)*\bigl(e^{\lambda s}\exp_*W\bigr)$ and, by~\eqref{eq:scalar-eq},
$\mathcal E[A]=Z_b*\bigl(\delta+e^{\lambda s}U\bigr)$. Multiplying by $\lambda e^{-\lambda s}$,
we obtain~\eqref{eq:C8}. Let $F_k=f\,\mathbf 1_{k=0}+f*U_k$. For $k\ne0$, from~\eqref{eq:Phi}
and~\eqref{eq:C8}
\[
e^{-\lambda s}A'_k=\frac{2F_k+2\sqrt2\,\lambda\,(F_{k-2}+F_{k+2})}{\lambda\,(k-s)^2}
=\ell_\pm f\,\mathbf 1_{k=\pm2}+\frac2{(k-s)^2}\,\mathcal T_\lambda\bigl[U_k,\ U_{k-2}+U_{k+2}\bigr].
\]
This is the first formula~\eqref{eq:C9}. For $k=0$ subtract from~\eqref{eq:Phi} the second
equality~\eqref{eq:scalar-eq}:
\[
s^2\,(A'_0-X_b)=2\,\bigl(\mathcal E[A]_0-Z_b\bigr)+2\sqrt2\,\lambda\,\bigl(\mathcal E[A]_2+\mathcal E[A]_{-2}\bigr)-4\lambda^2\,Z_b*Z_b .
\]
Multiplying by $e^{-\lambda s}$ and using~\eqref{eq:C8}, we get
$s^2S'=2\,\tfrac f\lambda*U_0+2\sqrt2\,f*(U_2+U_{-2})-4\,f*f$. It remains to substitute
$U_{\pm2}=\ell_\pm f+D_{\pm2}+N_{\pm2}$ and to note that
$2\sqrt2\,(\ell_++\ell_-)-4=4\sqrt2\,\tilde d$. The right-hand sides of~\eqref{eq:C9} are
nonnegative, hence $A'\ge H^{(0)}$.
\end{proof}

Thus, in the variables $(D,S)$ the operator $\Phi$ is given by formulas~\eqref{eq:C9}, where $W$,
$U$, $N$ are expressed through $D$ and $S$ by~\eqref{eq:remainders}. The right-hand sides
of~\eqref{eq:C9} are nondecreasing in $D$ and $S$. It is essential that the negative term
$-4f*f$ is cancelled by the explicit part of the modes $\pm2$ and that the remaining coefficient
$\tilde d$ is nonnegative.

\begin{lemma}[sufficient condition]\label{lem:sufficient}
Let the family $H$ be admissible, let $U$ be its remainder, and let $\int_0^1U_2(s)\,ds\le m<1$.
Then
\[
\frac{G_B(\mu)}{2\pi}\ \ge\ (1-m)\,\frac{\sinh\lambda}{\lambda}.
\]
\end{lemma}

\begin{proof}
By~\eqref{eq:C3} and~\eqref{eq:C8} with $A=H$
\[
\frac{G_B(\mu)}{2\pi}=E_0(1)-E_2(1)=\frac{e^{\lambda}}\lambda\,\bigl[f(1)+(f*U_0)(1)-(f*U_2)(1)\bigr].
\]
The function $f$ is increasing and $U_2\ge0$, so
$(f*U_2)(1)=\int_0^1f(1-t)\,U_2(t)\,dt\le f(1)\,m$. Moreover, $f*U_0\ge0$. It remains to note that
$f(1)\,e^{\lambda}/\lambda=\sinh\lambda/\lambda$.
\end{proof}

\subsection{Eliminating the dependence on the amplitude}\label{ss:uniformity}

Formulas~\eqref{eq:C9} depend on $\lambda$ through the functions $f$ and $f/\lambda$. The
following identity allows one to replace them by their limits as $\lambda\to\infty$.

\begin{lemma}\label{lem:uniformity}
Let $e_\lambda(s)=e^{-2\lambda s}$. For functions $X$, $V$ on $(0,1]$ and any $\lambda>0$
\begin{equation}\label{eq:U1}
\mathcal T_\lambda[X,V]=\frac{\sqrt2}2\,1*V+e_\lambda*\Bigl[\,1*X-\frac{\sqrt2}2\,V\Bigr].
\end{equation}
Consequently, if $X,V\ge0$ and $1*X\le\tfrac{\sqrt2}2V$ on $(0,1]$, then
$\mathcal T_\lambda[X,V]\le\tfrac{\sqrt2}2\,1*V$ for all $\lambda>0$.
\end{lemma}

\begin{proof}
We have $f/\lambda=1*e_\lambda$, since $\int_0^se^{-2\lambda t}dt=(1-e^{-2\lambda s})/(2\lambda)$,
and $f=\tfrac12(1-e_\lambda)$. Therefore $(f/\lambda)*X=e_\lambda*(1*X)$ and
$\sqrt2\,f*V=\tfrac{\sqrt2}2\,1*V-\tfrac{\sqrt2}2\,e_\lambda*V$. Adding gives~\eqref{eq:U1}, and
the second statement follows from $e_\lambda\ge0$.
\end{proof}

The meaning of the lemma is as follows. At finite amplitude, \eqref{eq:C9} contains an extra
positive term with $f/\lambda$, but at the same time $f$ is smaller than its limit $\tfrac12$.
Identity~\eqref{eq:U1} shows that these two differences from the limit system combine into a
single term whose sign does not depend on $\lambda$.

\paragraph{Sums over classes of modes} We divide the modes $k\in2\ZZ$ into six classes according
to the residue of $k/2$ modulo $6$. For a family $A$ set
\[
A_{\langle i\rangle}=\sum_{k/2\,\equiv\,i\ (\mathrm{mod}\ 6)}A_k,\qquad i=0,\dots,5.
\]
Class indices are added modulo $6$. Passing to sums over classes commutes with convolution, that
is, $(A*B)_{\langle i\rangle}=\sum_hA_{\langle h\rangle}*B_{\langle i-h\rangle}$. Let
\begin{equation}\label{eq:K-i}
K_i(s)=\frac2{(k_i-s)^2},\qquad(k_0,\dots,k_5)=(12,\,2,\,4,\,6,\,-4,\,-2).
\end{equation}
The number $k_i$ is the nonzero frequency of class $i$ closest to the interval $[0,1]$, so
$2/(k-s)^2\le K_i(s)$ for all frequencies $k\ne0$ of class $i$ and $s\in[0,1]$.

\begin{lemma}\label{lem:classes}
Under the assumptions of Lemma~\ref{lem:transform}
\begin{equation}\label{eq:U2}
\begin{aligned}
W_{\langle i\rangle}&\le\frac{\sqrt2}{(2-s)^2}\,\mathbf 1_{i=1}+\frac{\sqrt2}{(2+s)^2}\,\mathbf 1_{i=5}
+D_{\langle i\rangle}+S\,\mathbf 1_{i=0},\\
D'_{\langle i\rangle}&\le K_i\,\mathcal T_\lambda\bigl[U_{\langle i\rangle},\ U_{\langle i-1\rangle}+U_{\langle i+1\rangle}\bigr],\\
s^2S'&\le2\,\mathcal T_\lambda\bigl[U_{\langle0\rangle},\ \tilde d+D_{\langle1\rangle}+D_{\langle5\rangle}+N_{\langle1\rangle}+N_{\langle5\rangle}\bigr],
\end{aligned}
\end{equation}
where $U_{\langle\cdot\rangle}=\exp_*W_{\langle\cdot\rangle}-\delta$ and
$N_{\langle\cdot\rangle}=U_{\langle\cdot\rangle}-W_{\langle\cdot\rangle}$.
\end{lemma}

\begin{proof}
The first line follows from definition~\eqref{eq:remainders} and the inequality
$\ell_\pm f\le\ell_\pm/2$. The expression $\mathcal T_\lambda$ is linear and nondecreasing in both
arguments, so the sum of the first formula~\eqref{eq:C9} over $k\ne0$ of class $i$ does not
exceed $K_i\,\mathcal T_\lambda\bigl[\sum_kU_k,\sum_k(U_{k-2}+U_{k+2})\bigr]$, and this is at most
the right-hand side of the second line. In the third line we used that $2\tilde d f\le\tilde d$,
and that the modes $U_0$, $D_{\pm2}$, $N_{\pm2}$ do not exceed the sums over the classes $0$, $1$,
$5$.
\end{proof}

\subsection{The comparison theorem}\label{ss:comparison}

Let $r_+=\tfrac{141422}{100000}$ and $r_-=\tfrac{141421}{100000}$, so that $r_-<\sqrt2<r_+$. Let
bounded functions $\hat B_1$, $\hat B_5$, $\hat d$, $\hat d_-$ be given on $(0,1]$ such that
\begin{equation}\label{eq:majorants}
\begin{gathered}
\hat B_1\ge\frac{r_+}{(2-s)^2},\qquad\hat B_5\ge\frac{r_+}{(2+s)^2},\\
\hat d\ge\frac{r_+\,s^2(12-s^2)}{2\,(4-s^2)^2},\qquad
0\le\hat d_-\le\frac{r_-\,s^2(12-s^2)}{2\,(4-s^2)^2},
\end{gathered}
\end{equation}
and set $\hat B_i=0$ for $i\ne1,5$. Let, further, nonnegative bounded measurable functions
$\bar D_i$, $\bar N_i$ ($i=0,\dots,5$) and $\bar S$ be given on $(0,1]$. Set
\begin{equation}\label{eq:upper-functions}
\bar W_i=\hat B_i+\bar D_i+\bar S\,\mathbf 1_{i=0},\qquad\bar U_i=\bar W_i+\bar N_i,\qquad
\bar R=\bar D_1+\bar D_5+\bar N_1+\bar N_5 .
\end{equation}
Consider two groups of conditions on $(0,1]$:
\begin{align}
&\bar D_i\ge K_i\,\frac{r_+}2\,1*\bigl(\bar U_{i-1}+\bar U_{i+1}\bigr),\qquad
s\,\bar N_i\ge\sum_{h=0}^{5}(s\bar W_h)*\bar U_{i-h},\notag\\
&s^2\bar S\ge r_+\,1*\bigl(\hat d+\bar R\bigr); \tag{I}\label{eq:cond-I}\\
&1*\bar U_i\le\frac{r_-}2\,\bigl(\bar U_{i-1}+\bar U_{i+1}\bigr)\quad(i=0,\dots,5),\qquad
1*\bar U_0\le\frac{r_-}2\,\bigl(\hat d_-+\bar R\bigr). \tag{II}\label{eq:cond-II}
\end{align}
Conditions~\eqref{eq:cond-I} are the inequalities of the limit system ($\lambda=\infty$), and
conditions~\eqref{eq:cond-II} make it possible to apply Lemma~\ref{lem:uniformity}. The amplitude
does not enter them.

\begin{proposition}[comparison]\label{prop:comparison}
Let conditions~\eqref{eq:cond-I} and~\eqref{eq:cond-II} hold. Then for every $\mu<0$ the family
$H$ is admissible, and its remainders satisfy on $(0,1]$
\begin{equation}\label{eq:comparison}
D_{\langle i\rangle}\le\bar D_i,\qquad N_{\langle i\rangle}\le\bar N_i,\qquad S\le\bar S,\qquad
U_{\langle i\rangle}\le\bar U_i\qquad(i=0,\dots,5).
\end{equation}
In particular, all the series~\eqref{eq:EH-sums} converge, and
$\int_0^1U_2(s)\,ds\le\int_0^1\bar U_1(s)\,ds$.
\end{proposition}

\begin{proof}
Let $H^{(n)}=\Phi^n[H^{(0)}]$. We prove by induction on $n$ that the family $H^{(n)}$ is
admissible and its remainders satisfy the inequalities $D^{(n)}_{\langle i\rangle}\le\bar D_i$,
$S^{(n)}\le\bar S$. For $n=0$ this holds, since $D^{(0)}=0$ and $S^{(0)}=0$.

\emph{Induction step.} Suppose the claim holds for $n$. Denote by $W$, $U$, $N$ the sums over
classes of the remainders of the family $H^{(n)}$. By the first line of~\eqref{eq:U2} and
inequalities~\eqref{eq:majorants}, $W_{\langle i\rangle}\le\bar W_i$. We show that
$N_{\langle i\rangle}\le\bar N_i$. Let
$(\mathsf Ku)_i=s^{-1}\sum_h(sW_{\langle h\rangle})*u_{i-h}$. By~\eqref{eq:N-eq},
$N=\mathsf K[W+N]$, whence $N=\sum_{m\ge1}\mathsf K^mW$, because the operator $\mathsf K$ is of
Volterra type with a bounded kernel, and $\mathsf K^mN\to0$. On the other hand, by the second
condition~\eqref{eq:cond-I} and the inequality $W\le\bar W$ we have
$\bar N\ge\mathsf K[W+\bar N]$, whence
$\bar N\ge\sum_{m=1}^{n'}\mathsf K^mW+\mathsf K^{n'}\bar N\ge\sum_{m=1}^{n'}\mathsf K^mW$ for every
$n'$. Hence $N_{\langle i\rangle}\le\bar N_i$ and $U_{\langle i\rangle}\le\bar U_i$.

Now we estimate the remainders of the family $H^{(n+1)}=\Phi[H^{(n)}]$. By the second line
of~\eqref{eq:U2}, the monotonicity of $\mathcal T_\lambda$, Lemma~\ref{lem:uniformity} and
conditions \eqref{eq:cond-II}, \eqref{eq:cond-I}
\[
D^{(n+1)}_{\langle i\rangle}\le K_i\,\mathcal T_\lambda\bigl[\bar U_i,\ \bar U_{i-1}+\bar U_{i+1}\bigr]
\le K_i\,\frac{\sqrt2}2\,1*\bigl(\bar U_{i-1}+\bar U_{i+1}\bigr)\le\bar D_i .
\]
Lemma~\ref{lem:uniformity} applies, since
$1*\bar U_i\le\tfrac{r_-}2(\bar U_{i-1}+\bar U_{i+1})\le\tfrac{\sqrt2}2(\bar U_{i-1}+\bar U_{i+1})$.
Similarly, by the third line of~\eqref{eq:U2}
\[
s^2S^{(n+1)}\le2\,\mathcal T_\lambda\bigl[\bar U_0,\ \tilde d+\bar R\bigr]
\le\sqrt2\,1*\bigl(\tilde d+\bar R\bigr)\le r_+\,1*\bigl(\hat d+\bar R\bigr)\le s^2\bar S .
\]
Here we used the inequalities $\tfrac{r_-}2\hat d_-\le\tfrac{\sqrt2}2\tilde d$ and
$\sqrt2\,\tilde d\le r_+\hat d$, which follow from~\eqref{eq:majorants}, \eqref{eq:d-tilde} and
$r_-^2<2<r_+^2$. Thus the remainders of the family $H^{(n+1)}$ are bounded, the family is
admissible by Lemma~\ref{lem:transform}, and the induction step is proved.

\emph{Passage to the limit.} By Lemma~\ref{lem:zero-mode}, $0\le H^{(0)}\le H$, so by
Lemma~\ref{lem:minimal} $H^{(n)}\to H$ pointwise. The sums over the modes of the families
$H^{(n)}$ are bounded uniformly in $n$. By Fatou's lemma for series, the family $H$ is admissible,
and its remainders satisfy $D_{\langle i\rangle}\le\bar D_i$, $S\le\bar S$. The inequalities for
$N$ and $U$ are obtained from them in the same way as in the induction step. Finally, the mode
$U_2$ does not exceed the sum $U_{\langle1\rangle}$.
\end{proof}

\begin{remark}\label{rem:comparison-proof}
The proof uses neither the continuity of the functions $H_k$, nor their behavior as $s\to0$, nor
the convergence of the series~\eqref{eq:EH-sums}, which is obtained as a consequence.
Conditions~\eqref{eq:cond-II} are imposed on the upper functions, not on the true ones. For the
true remainders at finite amplitude they fail near zero.
\end{remark}

\subsection{The certificate and the proof of Theorem~\ref{thm:cos2}}\label{ss:certificate}

It remains to exhibit functions satisfying conditions~\eqref{eq:cond-I} and~\eqref{eq:cond-II}.
All of them will be polynomials with rational coefficients.

\paragraph{Majorants} Set
\begin{equation}\label{eq:poly-majorants}
\begin{aligned}
\hat B_1&=\frac{r_+}4\Bigl[\sum_{k=0}^{6}(k+1)\Bigl(\frac s2\Bigr)^k+\frac9{64}\,s^7\Bigr],\\
\hat B_5&=\frac{r_+}4\sum_{k=0}^{6}(k+1)\Bigl(-\frac s2\Bigr)^k,\\
\hat d&=\frac{r_+}2\cdot\frac{s^2(12-s^2)}{16}\Bigl[1+\frac{s^2}2+\frac{3s^4}{16}+\frac{13}{144}\,s^6\Bigr],\\
\hat d_-&=\frac{r_-}2\cdot\frac{s^2(12-s^2)}{16}\Bigl[1+\frac{s^2}2+\frac{3s^4}{16}\Bigr].
\end{aligned}
\end{equation}
Inequalities~\eqref{eq:majorants} for them follow from the identities
\begin{equation}\label{eq:factorizations}
\begin{aligned}
\hat B_1\,(2-s)^2-r_+&=\frac{r_+}{256}\,s^7(1-s)(20-9s),\\
\hat B_5\,(2+s)^2-r_+&=\frac{r_+}{256}\,s^7(16+7s),\\
2\,(4-s^2)^2\,\hat d-r_+\,s^2(12-s^2)&=\frac{r_+}{2304}\,s^8(12-s^2)(1-s^2)(64-13s^2),\\
r_-\,s^2(12-s^2)-2\,(4-s^2)^2\,\hat d_-&=\frac{r_-}{256}\,s^8(12-s^2)(16-3s^2),
\end{aligned}
\end{equation}
whose right-hand sides are nonnegative on $[0,1]$.

\paragraph{Upper functions} Set
\begin{equation}\label{eq:poly-upper}
\bar D_i=s^{p_i}\,g^D_i(s),\qquad\bar N_i=s^{p_i}\,g^N_i(s),\qquad\bar S=s\,g^S(s),
\end{equation}
where $p_i=1$ for even $i$ and $p_i=2$ for odd $i$, and $g^D_i$, $g^N_i$, $g^S$ are cubic
polynomials whose coefficients are given in Table~\ref{tab:certificate}.

\begin{table}[ht]
\caption{Coefficients of $1$, $s$, $s^2$, $s^3$ in the polynomials $g^D_i$, $g^N_i$, $g^S$,
given as numerators of fractions with denominator $10^4$. For example,
$g^S(s)=0.5132+0.0409\,s+0.0071\,s^2+0.2285\,s^3$.}\label{tab:certificate}
\centering\small
\renewcommand{\arraystretch}{1.12}
\begin{tabular}{@{}c@{\qquad}rrrr@{\qquad\quad}rrrr@{}}
\toprule
$i$ & \multicolumn{4}{c}{$g^D_i$} & \multicolumn{4}{c}{$g^N_i$} \\
\midrule
0 & 73 & 24 & $-15$ & 70 & 1302 & 299 & $-152$ & 1647 \\
1 & 1334 & 3252 & $-5707$ & 9421 & 1147 & 502 & $-30$ & 803 \\
2 & 320 & 550 & $-635$ & 1326 & 642 & 1079 & $-1095$ & 2639 \\
3 & 40 & 21 & $-9$ & 49 & 202 & 40 & 4 & 231 \\
4 & 328 & $-313$ & 229 & $-56$ & 655 & $-610$ & 484 & $-46$ \\
5 & 1393 & $-1438$ & 1108 & $-304$ & 1150 & $-344$ & 303 & 110 \\
\midrule
\multicolumn{9}{@{}l}{$g^S$:\quad 5132,\ \ 409,\ \ 71,\ \ 2285} \\
\bottomrule
\end{tabular}
\end{table}

\begin{lemma}[certificate]\label{lem:certificate}
The polynomials $g^D_i$, $g^N_i$, $g^S$ are positive on $[0,1]$. The
functions~\eqref{eq:poly-upper} together with the majorants~\eqref{eq:poly-majorants} satisfy
conditions~\eqref{eq:cond-I} and~\eqref{eq:cond-II}, and inequalities~\eqref{eq:cond-I} are
strict. Moreover,
\begin{equation}\label{eq:m-prime}
\int_0^1\bar U_1(s)\,ds=\frac{96257079}{102400000}=1-\varkappa_0<0.9401 .
\end{equation}
\end{lemma}

\begin{proof}
All the functions involved are polynomials with rational coefficients, and the convolution of
polynomials is computed exactly by the formula $s^x*s^y=x!\,y!\,s^{x+y+1}/(x+y+1)!$. After
multiplying the first inequality~\eqref{eq:cond-I} by $(k_i-s)^2$, each condition takes the form
``a polynomial is nonnegative on $(0,1]$''. After division by a suitable power of $s$ one obtains
a statement on the positivity or nonnegativity of a polynomial on the interval $[0,1]$. There are
$13+13+7=33$ such statements: the positivity of the thirteen polynomials $g$, the thirteen
inequalities~\eqref{eq:cond-I} and the seven inequalities~\eqref{eq:cond-II}. Together with the
four inequalities~\eqref{eq:majorants} there are $37$ of them. They have been verified in exact
rational arithmetic in two independent ways: by Bernstein coefficients with bisection of the
interval and by counting real roots with Sturm's method. Equality~\eqref{eq:m-prime} is obtained
by exact integration. The list of inequalities, a description of the verification and its results
are given in Appendix~\ref{app:certificate}.
\end{proof}

\begin{proof}[Proof of Theorem~\ref{thm:cos2}]
Formula~\eqref{eq:GB} is obtained from~\eqref{eq:G} with $M=\cos2\theta$, $p=\cos\theta$. If
$\mu\ge0$, then $G_B(\mu)\ge2\pi$ by item~(3) of Theorem~\ref{thm:series}. Let $\mu<0$. By
Lemma~\ref{lem:certificate} and Proposition~\ref{prop:comparison} the family $H$ is admissible
and $\int_0^1U_2\,ds\le1-\varkappa_0$. By Lemma~\ref{lem:sufficient},
$G_B(\mu)/2\pi\ge\varkappa_0\sinh\lambda/\lambda$, which is~\eqref{eq:GB-negative}. Since
$\sinh\lambda/\lambda\ge1$, the function $G_B$ is bounded away from zero on the real axis, and the
asymptotics~\eqref{eq:example} with $G_B$ in place of $G$ follows from Theorem~\ref{thm:main} with
$A_f=-1$, $p_f=\cos\theta$.
\end{proof}

\begin{remark}[a consequence for the background $\cos4\theta$]\label{rem:cos4-growth}
The inequality $E_0\ge Z_b$ of Lemma~\ref{lem:zero-mode} is proved for a positive background and
does not use a computer. If $x'(\xi,\theta)$ is a formal solution of equation~\eqref{eq:newton}
with background $\cos2\theta$, then $x(\xi,\theta)=\tfrac12x'(2\xi,2\theta)$ is a solution with
background $\cos4\theta$ and the same amplitude. Therefore the zero mode of the Borel transform of
the function $x^{-2}$ for the background $\cos4\theta$ at the point $s$ equals twice the zero mode
for the background $\cos2\theta$ at the point $s/2$. This gives, for the function $G$ of
Theorem~\ref{thm:example}, the bound
\[
\frac{G(\mu)}{2\pi}=2E_0\bigl(\tfrac12\bigr)\ge2Z_b\bigl(\tfrac12\bigr)
=\frac{\sinh\lambda'}{\lambda'},\qquad\lambda'=\frac{|\mu|}{\sqrt{32}},
\]
which sharpens the inequality $G\ge2\pi$. Here $E_0$ and $Z_b$ are built from the background
$|\mu|\cos2\theta$.
\end{remark}
\section{Discussion}\label{sec:discussion}

This section collects remarks on the theorems and open questions. Some of the statements are
based on heuristic arguments. We say so each time, and such statements are not proofs.

\subsection{First-order perturbation theory}\label{ss:first-order}

\paragraph{The first-order formula} Suppose there is no strong term, that is, $\mu=0$. Then for
$\tau=0$ there is no perturbation, $\mathbf D_0=0$, and by~\eqref{eq:F-coeffs}, for fixed $h$,
\[
F=\tau^2\det\mathbf H+O(\tau^3),
\]
where $\mathbf H$ is the second derivative at zero of the first variation of the difference of the
generating functions with respect to $\tau$. The unperturbed trajectory with $q(0)=\pi$ is
$q^0(t)$, and $1-\cos q^0(t)=2\cosh^{-2}t$. Therefore, by~\eqref{eq:first-variation}, the first
variation on the section $u=0$ equals
\begin{equation}\label{eq:melnikov}
\begin{gathered}
\mathcal M_f(\alpha)=\int_{-\infty}^{\infty}f(\alpha+\omega t)\,\frac{2\,dt}{\cosh^2t}
=\sum_{j,\ell}\fc f{j,\ell}\,\hat K\bigl(j\nu+\ell/h\bigr)\,e^{i(jx+\ell\theta)},\\
\hat K(\sigma)=\frac{2\pi\sigma}{\sinh(\pi\sigma/2)} .
\end{gathered}
\end{equation}
This is the Melnikov potential. The function $\hat K$ is even, $\hat K(\sigma)=4+O(\sigma^2)$ as
$\sigma\to0$ and $\hat K(\sigma)=4\pi|\sigma|e^{-\pi|\sigma|/2}\bigl(1+O(e^{-\pi|\sigma|})\bigr)$
as $|\sigma|\to\infty$. Since $|\nu|=\nu_0h$, the harmonics with $\ell=0$ enter~\eqref{eq:melnikov}
with the factor $4+O(h^2)$, the harmonics with $\ell=\pm1$ with the factor $4\pi L_h(1+O(h))$, and
the harmonics with $|\ell|\ge2$ with factors $O(h^{-1}e^{-2a/h})$. The evenness of $f$ gives
$\sum_j\fc f{j,1}=\sum_j\fc f{j,-1}=\fc{p_f}1$. Differentiating~\eqref{eq:melnikov} twice, we
obtain
\begin{equation}\label{eq:H-first-order}
\mathbf H_{11}=4A_f+O(h^2),\qquad
\mathbf H_{22}=-8\pi\,\fc{p_f}1\,L_h+O(hL_h),\qquad
\mathbf H_{12}=O(L_h),
\end{equation}
whence $\det\mathbf H=-32\pi A_f\,\fc{p_f}1\,L_h+O(hL_h)$. The conversion to the original
variables is the same as in Section~\ref{ss:thmQ}:
\begin{equation}\label{eq:first-order-Delta}
\Delta\approx-\frac{32\pi|\bom_2|}{g^3}\,A_f\,\fc{p_f}1\;\eps^2\,\eta^{-\frac12}
\exp\Bigl[-\frac{\pi|\bom_2|}{2g\sqrt\eta}\Bigr].
\end{equation}
This is the prediction of first-order perturbation theory. It coincides with the right-hand side
of~\eqref{eq:main} at $\mu=0$, since $G_{M,p}(0)=4\pi\fc p1$. As noted in
Section~\ref{ss:mechanism}, it is valid for $\eps=\eta^c$, $c>c_0$, if $\epsp\to0$.

Formula~\eqref{eq:H-first-order} shows the mechanism. To leading order the determinant is the
product of the slow entry $\mathbf H_{11}$, which is not small, and the fast entry
$\mathbf H_{22}$, which is exponentially small. The factor $e^{-a/h}$ is produced only by the pair
``a harmonic of $f$ independent of the fast angle'' and ``a harmonic with $\ell=\pm1$''. Two fast
harmonics enter the determinant through $\mathbf H_{12}^2$ and give a quantity of order
$h^{-2}e^{-2a/h}$.

\paragraph{What the strong term changes} For $\mu\ne0$ the matrix $\mathbf H$ is still the second
derivative of a Melnikov-type integral, but taken along the trajectories of the problem with the
background. The slow entry remains the same to leading order, since
$(\mathbf H_E)_{11}=4A_f+O(h)$ by~\eqref{eq:HE}. In the fast entry the number $4\pi\fc{p_f}1$ is
replaced by the value of an entire function, namely
$(\mathbf H_O)_{22}=-2L_h\bigl(G_{M,p_f}(\mu)+o(1)\bigr)$ by~\eqref{eq:HO22}. Thus the whole
effect of the strong term on the leading term is the factor $G_{M,p_f}(\mu)/G_{M,p_f}(0)$.

\begin{remark}[the series in the amplitude of the strong term]\label{rem:leading-series}
The determinant $F$ is analytic in the pair $(\tau,z)$ (Section~\ref{ss:detdef}). Steps 1--5
of Section~\ref{ss:general-f} also go through for a complex background amplitude $z$,
$|z|\le R_0$, with constants independent of $z$ (Remark~\ref{rem:complex-mu}), and the convergence
$G_h\to G_{M,p_f}$ is locally uniform. Therefore $[\tau^2]\,F/L_h\to-8A_fG_{M,p_f}(z)$ uniformly
in the disc $|z|\le R_0$, and by Cauchy's formula
\[
\lim_{h\to0}\frac{[\tau^2z^j]\,F}{L_h}=-8A_f\,[\mu^j]\,G_{M,p_f}\qquad(j\ge0).
\]
In other words, the leading terms of orders $\eps^2\epsp^{\,j}$ of the perturbation series for
$\Delta$ add up to the Taylor series of an entire function. This series converges for all
$\epsp$, and its sum is the coefficient of the asymptotics~\eqref{eq:main}. At the same time, the
terms of all orders in $\epsp$ are of the same order $\eps^2\eta^{-1/2}e^{-b/\sqrt\eta}$,
$b=\pi|\bom_2|/(2g)$. Therefore, for $\epsp$ of order one, no finite part of the series gives the
asymptotics. This is precisely the difficulty pointed out in the problem
statement~\cite{Gallavotti98}. This remark is not used in the proofs, and we only sketch its
argument, because the uniformity of the estimates of Section~\ref{ss:general-f} in complex $z$
is not written out separately there.
\end{remark}

\paragraph{The power of $\eta$ in the problem statement} In the text of~\cite{Gallavotti98} the
leading term for small amplitudes is written as $C(f)\,\eps^2\eta^{-1}e^{-b\eta^{-1/2}}$.
Formula~\eqref{eq:first-order-Delta} gives the power $\eta^{-1/2}$ for the
Hamiltonian~\eqref{eq:H} and any $f$. Indeed, the factor $e^{-b/\sqrt\eta}$ appears together
with a single factor $1/h$, while the product of two fast entries would give $\eta^{-1}$, but
with a doubled exponent.

In~\cite{GGM99}, to which the problem statement refers, the splitting is defined in the same
way. Formula~(4.2) of that paper considers the determinant of the matrix
$\partial_\alpha\bigl(A^+(\pi,\alpha)-A^-(\pi,\alpha)\bigr)$ at
$\alpha=0$. For the example with the perturbation
$\cos(\alpha+\varphi)+\cos(\lambda+\varphi)$ and $g=1$, formula~(6.14) there gives
$\det D=32\pi\,\eta^{-1/2}\eps^2e^{-\frac\pi2\eta^{-1/2}}$, that is, the power $\eta^{-1/2}$.
In~\cite{GGM99} the power $\eta^{-1}$ occurs only for an auxiliary quantity in formula~(6.11).
Therefore $\eta^{-1}$ in the text of the problem statement is apparently a misprint rather than a
different normalization of $\Delta$. The model of~\cite{GGM99} is more general. There the
perturbation is an arbitrary even trigonometric polynomial in $\varphi$ and two angles, and the
power of $\eta$ in the leading term depends on which harmonics in $\varphi$ enter the
perturbation.

\subsection{On the exponent \texorpdfstring{$c$}{c}}\label{ss:exponent}

\paragraph{Where $c_0$ comes from} In Proposition~\ref{thm:Qprime} the weak amplitude is bounded
from above by the condition $|\tau|\le h^{3P+2}$, where $P=\Lambda+5$ and $\Lambda=2C_S+1$
(Remark~\ref{rem:P-value}). It results from two requirements. First, we continue the perturbed
manifolds to a distance of order $h$ from the singularities of the separatrix by comparison with
the background ones, and the Gronwall estimate in Lemma~\ref{lem:continuation} requires the sum of
the weak amplitudes not to exceed $c_3(\kappa h)^{\Lambda+2}$. Here $C_S=46$ is the crude
constant of Lemma~\ref{lem:geometry}. Second, we estimate the terms of order $\tau^3$ and higher
by Cauchy's inequalities on the disc $|\tau|\le h^P$, which gives the relative remainder
$O(h^{-1-3P}|\tau|)$. Neither is optimal, and the value $c_0=(3P+2)/2=148$ has no significance of
its own.

\paragraph{A lower bound for $c$: a heuristic argument} The condition that the weak amplitude
be small compared with $h$ is apparently essential. Consider the example of Proposition~\ref{thm:Q},
in which $M=\cos4\theta$ and $V=s\cos x+r\cos\theta$, and let $s$ be real. The trajectory passes
along the separatrix in a time of order $\log(1/h)$. The slow angle changes during this time by
$O(h\log(1/h))$, and to leading order it can be frozen. For $x=0$ we obtain the periodic problem
with $K=(1+s)+\mu M(\theta)+r\cos\theta$. The change of time $t\mapsto t\sqrt{1+s}$ brings it to
the form~\eqref{eq:A} with the parameters
\[
h'=h\sqrt{1+s},\qquad\mu'=\frac\mu{1+s},\qquad r'=\frac r{1+s},
\]
and multiplies the generating functions by $\sqrt{1+s}$. The slow entry of the matrix $\mathbf D$
then equals $-s\int_{-\infty}^\infty2\cosh^{-2}(t\sqrt{1+s})\,dt=-4s/\sqrt{1+s}$, and the fast
entry equals $-2\sqrt{1+s}\;r'L_{h'}G(\mu')$. Multiplying, we obtain instead of~\eqref{eq:Q}
\begin{equation}\label{eq:Theta}
\begin{gathered}
F\approx8\,G(\mu)\,s\,r\,L_h\cdot\Theta(s),\\
\Theta(s)=(1+s)^{-3/2}\exp\Bigl[\frac ah\Bigl(1-\frac1{\sqrt{1+s}}\Bigr)\Bigr]\,
\frac{G\bigl(\mu/(1+s)\bigr)}{G(\mu)} .
\end{gathered}
\end{equation}
As $s\to0$ we have $\Theta(s)=\exp\bigl[as/(2h)\bigr]\bigl(1+O(s+s^2/h)\bigr)$.

In the original variables, for the example of Theorem~\ref{thm:example}, where $s=r=\tau$, the
exponent in~\eqref{eq:Theta} for $\eps=\eta^c$ equals
\[
\frac{a\tau}{2h}=\frac{\pi|\bom_2|}{4Jg^3}\,\frac{\eps}{\sqrt\eta}
=\frac{\pi|\bom_2|}{4Jg^3}\;\eta^{\,c-\frac12}.
\]
Hence one should expect the following. For $c>\tfrac12$ the factor $\Theta$ tends to one and
formula~\eqref{eq:example} keeps its form, for $c=\tfrac12$ it tends to
$\exp\bigl[\pi|\bom_2|/(4Jg^3)\bigr]\ne1$, and for $c<\tfrac12$ it grows without bound. So for
$c\le\tfrac12$ the ratio of the left- and right-hand sides of~\eqref{eq:example} apparently does
not tend to one. This agrees with the remark in~\cite{GGM99} that the correct condition on the
exponent may be $c>\tfrac12$. For a general perturbation the same argument gives the exponent
$a\tau\bar f_0/(2h)$, where $\bar f_0=\sum_j\fc f{j,0}$ is the mean of the function $f(0,\theta)$
over $\theta$.

Thus the statements are proved for $c>c_0$, and the heuristic argument above suggests that they
are false in the stated form for $c\le\tfrac12$. The range $\tfrac12<c\le c_0$ remains open. To close it, one would apparently have
to continue the manifolds to the singularities taking the weak perturbation into account, rather
than by comparison with the background ones, and to estimate the higher variations in $\tau$
directly.

\subsection{Growing amplitude of the strong term}\label{ss:growing}

The problem statement~\cite{Gallavotti98} suggests, to fix ideas, setting $\epsp=\eta^T$ and notes
that there is a critical value of $T$ below which the asymptotics changes, the transition taking
place inside the domain of convergence of the perturbation series. According
to~\cite{Gallavotti98}, the series for $\Delta$ in powers of $\epsp$ converges for
$|\epsp|<O(\eta^{-b'})$ with any $b'<\tfrac12$ (see also~\cite[\S\,8]{GGM99}).

\paragraph{What is proved} The remainder in Theorem~\ref{thm:main} is uniform in $\epsp$ on
compact sets. Therefore, for $\eps=\eta^c$, $c>c_0$, the following holds. If $T>0$, then
$\epsp\to0$ and the asymptotics of $\Delta$ is the one given by the first order, that
is,~\eqref{eq:first-order-Delta} (Section~\ref{ss:mechanism}). If $T=0$, that is, $\epsp$ does
not depend on $\eta$, then the coefficient of this asymptotics is multiplied by
$G_{M,p_f}(\mu)/G_{M,p_f}(0)$. Thus, for this choice of $\eps$ the asymptotics changes at $T=0$.
For $T<0$ the amplitude $\mu$ grows without bound as $\eta\to0$, and the theorems do not apply,
because the constants of the lemmas depend on the bound $R_0$ for $|\mu|$, and we have not
traced this dependence.

\paragraph{Growth of the coefficient} In both readings of the example the coefficient grows
exponentially with the amplitude. The following bounds are proved:
\begin{equation}\label{eq:growth}
\begin{gathered}
\frac{G(\mu)}{2\pi}\ge\frac{\sinh\lambda'}{\lambda'},\qquad\lambda'=\frac{|\mu|}{\sqrt{32}};\\
\frac{G_B(\mu)}{2\pi}\ge\frac{\sinh\lambda}{\lambda}\quad(\mu\ge0),\qquad
\frac{G_B(\mu)}{2\pi}\ge\varkappa_0\,\frac{\sinh\lambda}{\lambda}\quad(\mu<0),\qquad
\lambda=\frac{|\mu|}{\sqrt8}.
\end{gathered}
\end{equation}
The first is Remark~\ref{rem:cos4-growth}, and the third is Theorem~\ref{thm:cos2}. The second
follows from Lemmas~\ref{lem:two-modes} and~\ref{lem:zero-mode}, since for $\mu>0$ we have
$G_B(\mu)/2\pi=E_0(1)+E_2(1)\ge E_0(1)\ge Z_b(1)$. The asymptotics of $G_{M,p}(\mu)$ as
$|\mu|\to\infty$ is unknown.

A function of the same form is known from the periodic problem with the single-harmonic profile
$\mu\cos\theta$, that is, from the problem of a pendulum whose suspension point oscillates
rapidly in the vertical direction. There the leading coefficient of the splitting is $\mu$
times an entire function of $\mu^2$. Nine Taylor coefficients of this function were computed
numerically in~\cite{Treschev96}. It is noted there, with attribution to Sim\'o, that they are
close to the Taylor coefficients of a function which in our notation is
$2\sinh\lambda/\lambda$ with $\lambda=|\mu|/\sqrt2$ (see
also~\cite[Section~2.2]{Gelfreich97}). The relative differences of the nine numbers do not
exceed seven percent.

\paragraph{A heuristic explanation} Let $P_M$ be the primitive of the background with zero mean
(Section~\ref{ss:separatrix}), so that $P_M'=M$ and $\langle P_M\rangle_\theta=0$. We set
\[
\lambda=|\mu|\,\langle P_M^2\rangle_\theta^{1/2}.
\]
For the backgrounds $\cos2\theta$ and $\cos4\theta$ these are the numbers $\lambda$ and
$\lambda'$ from~\eqref{eq:growth}, and for the profile $\cos\theta$ this is the number
$|\mu|/\sqrt2$ mentioned above. The change $v=q'-h\mu P_M(\theta)\sin q$ in the equation
$q''=\bigl(1+\mu M(t/h)\bigr)\sin q$ gives
$v'=\sin q-h\mu P_M\,v\cos q-h^2\mu^2P_M^2\sin q\cos q$ (cf.~\eqref{eq:qv}). After averaging over
the fast angle one obtains the system with Hamiltonian
\begin{equation}\label{eq:averaged}
\tfrac12\,v^2+(\cos q-1)+\tfrac12\,(\lambda h)^2\sin^2q .
\end{equation}
Its separatrix can be found explicitly. If $\sigma_\lambda=\sqrt{1-\lambda^2h^2}$, then
\[
\cos\frac q2=-\frac{\tanh\sigma_\lambda t}{\sqrt{\sigma_\lambda^2+\lambda^2h^2\tanh^2\sigma_\lambda t}} .
\]
The pole $t=i\pi/2$ of the unperturbed separatrix splits into two branch points
\begin{equation}\label{eq:averaged-sing}
t=\frac{i}{\sigma_\lambda}\Bigl(\frac\pi2\mp\arcsin\lambda h\Bigr)
=i\Bigl(\frac\pi2\mp\lambda h\Bigr)+O(\lambda^2h^2).
\end{equation}
The distance from the real axis to the nearest singularity decreases by $\lambda h$, and it is
natural to expect that the exponent in the splitting changes accordingly, that is, that
$e^{-a/h}$ is replaced by $e^{-(a-\lambda h)/h}=e^{-a/h}e^{\lambda}$. In the inner variable
$\xi$ the function $x^{-2}$ then becomes $(\xi^2-\lambda^2)^{-1}$, and its Borel transform
becomes $\sinh(\lambda s)/\lambda$, which is the comparison function $Z_b$ of
Section~\ref{ss:zero-mode}. In other words, the lower bounds~\eqref{eq:growth} are the
contribution of the averaged system, and the exponential growth of the coefficient with the
amplitude is a trace of the shift of the singularity. The averaging describes the range
$1\ll|\mu|\ll h^{-1}$.

If this argument is correct, then for $\epsp=\eta^T$ with $-\tfrac12<T<0$ already the exponent
changes, because the term $\lambda\,(1+o(1))$, growing like $\eta^{T}$, is added to
$-b/\sqrt\eta$. For $T=-\tfrac12$ the quantity $\lambda h$ is of order one, and the
singularities~\eqref{eq:averaged-sing} move by a distance of order one. This is the same
boundary as that of the domain of convergence of the series in~\cite{Gallavotti98}. Everything
said in the last two paragraphs is heuristic. There is no proof, and it is unknown whether
formula~\eqref{eq:main} persists when the amplitude grows together with $1/\eta$.

\subsection{Open questions}\label{ss:open}

\begin{enumerate}[label=\textup{\arabic*.}]
\item \emph{The exponent $c$.} Is Theorem~\ref{thm:main} true for all $c>\tfrac12$? What is the
asymptotics for $c\le\tfrac12$? Is it true that it differs by a factor of the
form~\eqref{eq:Theta}?
\item \emph{Growing amplitude.} What is the asymptotics of $\Delta$ for $\epsp=\eta^T$,
$-\tfrac12<T<0$? What is the asymptotics of $G_{M,p}(\mu)$ as $|\mu|\to\infty$?
\item \emph{Degenerate cases.} Find the leading term when $A_fG_{M,p_f}(\mu)=0$.
\item \emph{The sign of the coefficient for general backgrounds.} The nonnegativity of the terms of
the series~\eqref{eq:G} is proved when the Fourier coefficients $\mu\fc Mn$ are nonnegative. For
the background $\cos2\theta$ with $\mu<0$ positivity is obtained by a separate argument with a
computer check. There is no general statement on the sign of $G_{M,p}$ when the coefficients
$\mu\fc Mn$ have different signs.
\item \emph{Positivity of $G_B$ without a computer.} The proof of
inequality~\eqref{eq:GB-negative} relies on the verification of $37$ polynomial inequalities. A
proof of the positivity of $G_B(\mu)$ for $\mu<0$ that does not use a computer is unknown.
\item \emph{Assumptions on the perturbation.} We use the evenness of $f$ and $m$ essentially,
since it gives reversibility, the homoclinic point $\alpha=0$ and the relation between $T^u$ and
$T^s$. We have not considered odd perturbations, analytic perturbations with infinite spectrum
or non-isochronous systems.
\end{enumerate}

\backmatter

\section*{Statements and Declarations}

\begin{itemize}
\item Funding: No funding was received for conducting this study.
\item Competing interests: The author declares no competing interests.
\item Data availability: All data supporting the results are contained in the paper. The
coefficients of the polynomials of Lemma~\ref{lem:certificate}, on which the computer-assisted
part of the proof relies, are given in full in Table~\ref{tab:certificate}, and the inequalities
to be verified are listed in Appendix~\ref{app:certificate}. A machine-readable copy of the
coefficients is included in the repository~\cite{GulkanovCode}. No external datasets were used.
\item Code availability: The programs that verify the inequalities of
Appendix~\ref{app:certificate}, the test with deliberate errors and the program that produced the
coefficients of Table~\ref{tab:certificate} are available in the Zenodo
repository~\cite{GulkanovCode}, \url{https://doi.org/10.5281/zenodo.24615549}.
\item Use of AI tools: The AI language models Claude Opus~5.0 and Claude Opus~5.5 (Anthropic)
were used to write the code. The author has checked all of the code and takes full
responsibility for it.
\end{itemize}

\begin{appendices}
\renewcommand{\theequation}{\thesection.\arabic{equation}}%
\counterwithin*{table}{section}%
\renewcommand{\theHtable}{\thetable}%
\section{Statements (E1) and (E2): relation to the literature and proof of (E1)}\label{app:bfgs}

This appendix concerns Propositions~\ref{prop:E1} and~\ref{prop:E2}, that is, statements (E1)
and (E2) of Section~\ref{ss:external}. We prove (E1) in Section~\ref{ss:A-direct}. In
Sections~\ref{ss:A-model}--\ref{ss:A-E1} we describe the setting of~\cite{BFGS12} and the
relation of (E1) to the results of that paper, and in Section~\ref{ss:A-E2} we explain how (E2)
follows from the results of~\cite{BFGS12} and~\cite{Baldoma06}. We do not reproduce the
proofs of the results of~\cite{BFGS12}. We check the conditions under which they apply to the
pendulum. In formulas written in the notation of~\cite{BFGS12}, the letters $\eps$, $\mu$,
$\eta$, $\beta$, $r$, $\ell$, $\rho$, $\kappa$, $\lambda$ have the meaning of that paper. The
correspondence with the notation of this paper is given in Table~\ref{tab:bfgs-dictionary}.

\subsection{The model}\label{ss:A-model}

The paper~\cite{BFGS12} considers systems with the Hamiltonian (formulas (7), (8) of that paper)
\begin{equation}\label{eq:bfgs-H}
H\Bigl(x,y,\frac t\eps;\eps\Bigr)=H_0(x,y)+\mu\,\eps^{\eta}H_1\Bigl(x,y,\frac t\eps;\eps\Bigr),
\qquad H_0(x,y)=\frac{y^2}2+V(x),
\end{equation}
where $V$ is a polynomial or a trigonometric polynomial, and the perturbation
$H_1(x,y,\tau;\eps)$ is $2\pi$-periodic in the fast time $\tau=t/\eps$. The periodic
problem~\eqref{eq:per} with profile $zg$ is obtained for
\begin{equation}\label{eq:bfgs-our}
\begin{gathered}
V(x)=\cos x-1,\qquad H_1(x,y,\tau)=g(\tau)\,(\cos x-1),\\
\eta=0,\qquad\eps=h,\qquad\mu=z .
\end{gathered}
\end{equation}
This is the Hamiltonian $p^2/2+(\cos q-1)\bigl(1+zg(\theta)\bigr)$ with $\theta'=1/h$.

\begin{table}[ht]
\caption{Correspondence of notation. The letter $z$ denotes the inner variable
in~\cite{BFGS12} and the background amplitude in this paper.}\label{tab:bfgs-dictionary}
\centering\small
\renewcommand{\arraystretch}{1.15}
\begin{tabular}{@{}>{\raggedright\arraybackslash}p{0.47\linewidth}%
  >{\raggedright\arraybackslash}p{0.45\linewidth}@{}}
\toprule
Paper~\cite{BFGS12} & This paper \\
\midrule
coordinate and momentum $x$, $y$ (after the shift $q$, $p$) & $q$, $p$ \\
small parameter $\eps$ & $h$ \\
fast time $\tau=t/\eps$ & fast angle $\theta$ \\
amplitude $\mu$, disc $B(\mu_0)$ & amplitude $z$, disc $|z|<R_0$ \\
separatrix $(q_0(u),p_0(u))$, singularities $\pm ia$ & $(q^0(u),\psi_0(u))$, $a=\pi/2$ \\
generating functions $T_0$, $T^u$; $T_1=T^u-T_0$ & the same; $\partial_uT_1=\psi_0^2E^u$ \\
domain constants $\rho_1$, $\rho_2$, $\kappa_1$ & $\varrho_1$, $\varrho_2$, $\kappa_1$ \\
angles $\beta_0<\beta_1<\beta_2$ & $\beta_0<\beta_1=\pi/4<\beta_2$ \\
inner variable $z=(u-ia)/\eps$ & $\zeta=(u-ia)/h$ \\
strip $\TT_\sigma$ & $\TT_\sigma$ \\
\bottomrule
\end{tabular}
\end{table}

\subsection{Verification of the hypotheses}\label{ss:A-hyp}

The hypotheses of~\cite{BFGS12} are denoted there by HP1--HP5 (pp.~3313, 3314). For
system~\eqref{eq:bfgs-our} they are satisfied.

\begin{itemize}
\item \emph{HP1.1} (hyperbolic point). $V(x)=-\tfrac12x^2+O(x^4)$, the eigenvalues are
$\pm\lambda$ with $\lambda=1$.
\item \emph{HP2} (separatrix). The parametrization $(q^0(u),\psi_0(u))$ from~\eqref{eq:q0} is
real analytic, analytic in the strip $|\Im u|<a$, $a=\pi/2$, and its only singularities on the
lines $\Im u=\pm a$ are the points $\pm ia$. The potential is a trigonometric polynomial of degree
$M=1$. If $w=u-ia$, then $\cosh u=i\sinh w$, and therefore, as $w\to0$,
\begin{gather*}
\psi_0(u)=-\frac{2i}{w}\bigl(1+O(w^2)\bigr),\\
\cos q^0(u)=\frac2{w^2}\bigl(1+O(w^2)\bigr),\qquad
\sin q^0(u)=\frac{2i}{w^2}\bigl(1+O(w^2)\bigr).
\end{gather*}
These are expansions (14) of~\cite{BFGS12} with $M=1$ and $C_+=-2i$. By the convention adopted
there, in the trigonometric case $r=1$ and $\beta=M=1$. The function $p_0=\psi_0$ has no zeros.
\item \emph{HP3} (perturbation). The function $H_1=g(\tau)(\cos x-1)$ is real analytic,
$2\pi$-periodic in $\tau$, is a trigonometric polynomial in $x$ of the form (10)
of~\cite{BFGS12} and has zero mean in $\tau$, since $\langle g\rangle=0$. It depends neither on
$y$ nor on $\eps$.
\item \emph{HP4.1} (order of the perturbation). $H_1=-\tfrac12g(\tau)x^2+\dots$, the order is
$n=2\ge1$.
\item \emph{HP5} (order of the singularity). By formula (15) of~\cite{BFGS12}, the order $\ell$ of
the singularity of the function $H_1(q_0(u),p_0(u),\tau)$ at the points $\pm ia$ equals $2|k|/M$
for the highest harmonic $e^{\pm ikx}$, that is, $\ell=2$. One needs
$\eta\ge\eta^*=\max\{0,\ell-2r\}$, and here $\eta=\eta^*=0$.
\end{itemize}

Thus $\ell-2r=0$ and $\eta=\ell-2r$. This case is called singular in~\cite{BFGS12}. In it the
perturbation is not small, and near the singularities the perturbed manifolds are not
approximated by the unperturbed separatrix. We note three more points.

\begin{itemize}
\item The origin remains a fixed point for all $z$, since $H_1$ and its gradient vanish at
$x=y=0$. Therefore the periodic orbit of Theorem 4.1 of~\cite{BFGS12} is the origin, the shift
(38) is the identity, and the Hamiltonian (39) coincides with the original one.
\item The width $\sigma$ of the strip in~\cite{BFGS12} is chosen smaller than the width of the
strip of analyticity of the perturbation in $\tau$. For a trigonometric polynomial $g$ we can
take any $\sigma>0$.
\item In Section 4.1 of~\cite{BFGS12} (p.~3330) an arbitrary $\mu_0>0$ is fixed, and it is
stipulated that all functions arising in the proofs are analytic in $\mu\in B(\mu_0)$ and the
constants in the estimates do not depend on $\mu$ and $\eps$. The disc $B(\mu_0)$ is complex.
Indeed, complex values of $\mu$ are used explicitly in the proof of Theorem 4.19 (p.~3349), and
this is also stated in~\cite{BMS24} (p.~1182), where Theorem 4.4 is restated. We take
$\mu_0=R_0+1$.
\end{itemize}

\subsection{The Hamilton--Jacobi equation and statement (E1)}\label{ss:A-E1}

\paragraph{The equation} The paper~\cite{BFGS12} makes the symplectic change of variables
$q=q_0(u)$, $p=w/p_0(u)$ (formula (45)) and considers the Hamiltonian
$\overline H(u,w,\tau)=\widehat H\bigl(q_0(u),w/p_0(u),\tau\bigr)$. For
system~\eqref{eq:bfgs-our}
\[
\overline H(u,w,\tau)=\frac{w^2}{2\psi_0(u)^2}-\bigl(1+zg(\tau)\bigr)\frac{\psi_0(u)^2}2 ,
\]
since $\cos q^0-1=-\psi_0^2/2$. The manifolds are sought in the form $w=\partial_uT(u,\tau)$,
where $T$ satisfies equation (47), which reads
$\partial_\tau T+\eps\,\overline H(u,\partial_uT,\tau)=0$. After division by $\eps=h$ this is
equation~\eqref{eq:HJ}. The parametrization (48) of the manifold, $p=\partial_uT/p_0(u)$,
coincides with the formula $p=\psi_0(1+E)$ of Section~\ref{ss:external}. The unperturbed
solution (57) is the function $T_0$ with $\partial_uT_0=p_0^2=\psi_0^2$. The unstable manifold
is given by the solution $T^u$ of Theorems 4.3 and 4.4 of~\cite{BFGS12}, which is close to $T_0$.

\paragraph{The statements of~\cite{BFGS12}} The domains (34) and (35) of~\cite{BFGS12} are
\begin{gather*}
D^u_{\infty,\rho_1}=\{u:\ \Re u<-\rho_1\},\\
D^{{\rm out},u}_{\rho,\kappa}=\{u:\ |\Im u|<-\tan\beta_1\,\Re u+a-\kappa\eps,\ \ \Re u>-\rho\},
\end{gather*}
where $\beta_1\in(0,\pi/2)$ is any fixed angle (p.~3326). For $\beta_1=\pi/4$ these are the
domains $D^{\rm loc}$ and $D^{\rm out}$ from~\eqref{eq:domains}. By Theorem 4.3 (p.~3334), if
$\rho_1$ is large enough, then for small $\eps$ and $\mu\in B(\mu_0)$ equation (47) has an
analytic solution $T^u$ in $D^u_{\infty,\rho_1}\times\TT_\sigma$ with
$|\partial_uT^u-\partial_uT_0|\le b_1|\mu|\eps^{\eta+1}$. Proposition 6.4 (p.~3369), from which
Theorem 4.3 is derived, gives for $T_1=T^u-T_0$ the weighted estimate
\begin{equation}\label{eq:bfgs-prop64}
\|\partial_uT_1\|_{2\lambda,\rho_1,\sigma}\le b_1|\mu|\eps^{\eta+1},\qquad
\|f\|_{\alpha,\rho,\sigma}=\sum_{k\in\ZZ}\ \sup_{u\in D^u_{\infty,\rho}}\bigl|e^{-\alpha u}\fc fk(u)\bigr|\,e^{|k|\sigma}
\end{equation}
(the norm is defined on p.~3367). By Theorem 4.4 (p.~3335), which treats the case $p_0\ne0$, if
$\rho_2>\rho_1$, $\kappa_1$ is large enough and $\eps_0$ is small enough, then for
$\mu\in B(\mu_0)$ and $\eps\in(0,\eps_0)$ the function $T^u$ extends analytically to
$D^{{\rm out},u}_{\rho_2,\kappa_1}\times\TT_\sigma$, and
\begin{equation}\label{eq:bfgs-thm44}
|\partial_uT^u(u,\tau)-\partial_uT_0(u)|\le\frac{b_2|\mu|\eps^{\eta+1}}{|u^2+a^2|^{\ell+1}}
\end{equation}
with constants $b_1$ and $b_2$ independent of $\eps$ and $\mu$.

\paragraph{Relation to (E1)} Let $\varrho_1=\rho_1$ and $\varrho_2=\rho_2$. By
definition~\eqref{eq:E-def}, $E^u=\partial_uT_1/\psi_0^2$.
\begin{itemize}
\item On $D^{\rm loc}$, \eqref{eq:bfgs-prop64} with $\lambda=1$, $\eta=0$ gives
$|\partial_uT_1(u,\theta)|\le b_1|z|h\,e^{2\Re u}$ for $\theta\in\TT_\sigma$. On the other hand,
$|\cosh u|\le\cosh\Re u\le e^{-\Re u}$ for $\Re u<0$, whence $|\psi_0(u)|^2\ge4e^{2\Re u}$.
Therefore $|E^u|\le\tfrac14b_1|z|h$.
\item On $D^{\rm out}$, \eqref{eq:bfgs-thm44} with $\eta=0$, $\ell=2$ and the equality
$|\psi_0|^{-2}=\tfrac14|\cosh u|^2$ give
$|E^u|\le\tfrac14b_2|z|h\,|\cosh u|^2/|u^2+a^2|^3$.
\end{itemize}
These are the estimates~\eqref{eq:E1} with $C_E=\tfrac14\max\{b_1,b_2\}$. The proof of
Proposition~\ref{prop:E1} in Section~\ref{ss:A-direct} does not use the results
of~\cite{BFGS12}.

\subsection{The inner equation and statement (E2)}\label{ss:A-E2}

\paragraph{Inner variables} The paper~\cite{BFGS12} introduces the inner variable (66) and the
scaled generating function (67). In our notation and with $r=1$, $C_+^2=-4$ these are
\begin{equation}\label{eq:bfgs-scaling}
\zeta=\frac{u-ia}h,\qquad\psi^u(\zeta,\theta)=h^{2r-1}C_+^{-2}\,T^u(ia+h\zeta,\theta)
=-\frac h4\,T^u(ia+h\zeta,\theta),
\end{equation}
so that
\begin{equation}\label{eq:bfgs-Phi}
h^2\,\partial_uT^u(ia+h\zeta,\theta)=-4\,\partial_\zeta\psi^u(\zeta,\theta).
\end{equation}
The function $\psi^u$ satisfies the Hamilton--Jacobi equation (68) with Hamiltonian (69). For
system~\eqref{eq:bfgs-our} its limit as $h\to0$ equals
$\tfrac12\zeta^2w^2-\bigl(1+zg(\theta)\bigr)/(2\zeta^2)$, since
$\psi_0(ia+h\zeta)^2=-4/\sinh^2(h\zeta)$. Therefore the inner equation (72) of~\cite{BFGS12} takes
the form
\begin{equation}\label{eq:bfgs-inner}
\partial_\theta\psi_0+\frac12\,\zeta^2(\partial_\zeta\psi_0)^2-\frac{1+zg(\theta)}{2\zeta^2}=0 .
\end{equation}

\paragraph{Domains} The inner domain (36) of~\cite{BFGS12} in the variable $\zeta$ (formula (70))
is given by the inequalities
\begin{equation}\label{eq:bfgs-inner-domain}
\begin{gathered}
\Im\zeta>-\tan\beta_1\bigl(\Re\zeta+c\,h^{\gamma-1}\bigr),\\
\Im\zeta<-\tan\beta_2\,\Re\zeta-\kappa,\qquad\Im\zeta<-\tan\beta_0\,\Re\zeta-\kappa ,
\end{gathered}
\end{equation}
where $\kappa>0$, $c>0$, $\gamma\in(0,1)$. According to pp.~3326, 3329 of that paper, $\beta_0$
is any angle for which the difference $\beta_1-\beta_0$ is bounded away from zero, and
$\beta_2>\beta_1$ is chosen so that there are no zeros of $p_0$ on the lines
$|\Im u|=\tan\beta_2\,\Re u+a$. For the pendulum $p_0$ has no zeros, and we take arbitrary angles
$\beta_0\in(0,\pi/4)$ and $\beta_2\in(\pi/4,\pi/2)$. The second and third
inequalities~\eqref{eq:bfgs-inner-domain} with $\kappa=\kappa_6$ define the sector $S^u$
from~\eqref{eq:sector}. The first inequality ceases to be a restriction as $h\to0$, since
$\gamma<1$. In other words, every point $\zeta\in S^u$ belongs to the
domain~\eqref{eq:bfgs-inner-domain} for $h$ small enough.

\paragraph{The statements used} \emph{Theorem 4.12, item 1} (p.~3343). It is quoted
in~\cite{BFGS12} from~\cite{Baldoma06}, where it is Theorem 3.1, and we state it in the form
given there. For any fixed $\theta_0>0$ there exists $\kappa_4>0$, depending on
$\theta_0$, $\sigma$, $\mu_0$ and $g$ but not on $\mu$, such that
equation~\eqref{eq:bfgs-inner} has a solution $\psi_0^u$, analytic in $(\zeta,\theta,\mu)$ on
$\{\zeta:\ |\Im\zeta|>\theta_0\Re\zeta+\kappa_4\}\times\TT_\sigma\times B(\mu_0)$, of the form
\begin{equation}\label{eq:bfgs-psi0}
\psi_0^u(\zeta,\theta)=-\frac1\zeta+\mu\,\overline\psi{}_0^{\,u}(\zeta,\theta)+K^u,\qquad
\sup_{\zeta,\theta,\mu}\bigl|\zeta^{\ell+1}\,\partial_\zeta\overline\psi{}_0^{\,u}\bigr|<\infty ,
\end{equation}
and for $\mu\ne0$ the derivative $\partial_\zeta\overline\psi{}_0^{\,u}$ is uniquely determined
by this condition.

\emph{Theorem 4.16} (p.~3345). Let $\gamma\in(0,1)$, $c_1>0$, let $\eps_0$ be small enough and
$\kappa_6$ large enough. Then for $\eps\in(0,\eps_0)$ and $\mu\in B(\mu_0)$, on the
domain~\eqref{eq:bfgs-inner-domain} with $\kappa=\kappa_6$, $c=c_1$ and on $\TT_\sigma$,
\begin{equation}\label{eq:bfgs-thm416}
|\partial_\zeta\psi^u(\zeta,\theta)-\partial_\zeta\psi_0^u(\zeta,\theta)|\le
\frac{b_{10}\,\eps^{1/\beta}}{|\zeta|^{2r-1/\beta}} ,
\end{equation}
where $\psi_0^u$ is the solution of Theorem 4.12 taken with $\theta_0=\tfrac12\tan\beta_2$. We
take $\kappa_6\ge2\kappa_4$, so that the domain~\eqref{eq:bfgs-inner-domain} lies in the domain
of $\psi_0^u$.

\paragraph{Derivation of (E2)} For $\beta=r=1$ the right-hand side of~\eqref{eq:bfgs-thm416}
equals $b_{10}h/|\zeta|$. Let $\zeta\in S^u$ and $\theta\in\TT_\sigma$ be fixed. For small $h$ the
point $\zeta$ lies in the domain~\eqref{eq:bfgs-inner-domain}, and \eqref{eq:bfgs-Phi},
\eqref{eq:bfgs-thm416} imply
\[
\lim_{h\to0}h^2\,\partial_uT^u(ia+h\zeta,\theta)=-4\,\partial_\zeta\psi_0^u(\zeta,\theta).
\]
This is statement (E2) with $\Phi_0^u=-4\,\partial_\zeta\psi_0^u$. Four clarifications are needed
here.

\begin{enumerate}[label=\textup{\arabic*.}]
\item \emph{Which function $T^u$.} Theorem 4.16 is stated for the general case, where $p_0$ may
have zeros. The function $T^u$ is then constructed by Theorem 4.8 on domains of a different shape,
and the constant $\kappa_3$ of that theorem enters the hypothesis. Remark 4.9 of~\cite{BFGS12}
(p.~3339) notes that for $p_0\ne0$ the function $T^u$ on these domains is already given by
Theorem 4.4 and that from then on the two cases are not distinguished. For small $h$ the
derivative $\partial_uT^u$ of this function coincides with that of the function $T^u$ of
Proposition~\ref{prop:E1}, as shown at the end of Section~\ref{ss:A-direct}. We need only that
the points $ia+h\zeta$ lie in the
domain $D^{\rm out}$ of statement (E1). For this we take $\kappa_6$
not smaller than $\kappa_1$. Then $S^u\subset\{\Im\zeta<-\Re\zeta-\kappa_1\}$, since
$\beta_0<\pi/4<\beta_2$.
\item \emph{Independence of $b_{10}$ from $\eps$.} To pass to the limit we need the constant
$b_{10}$ in~\eqref{eq:bfgs-thm416} to be independent of $\eps$. This is seen from the proof
(Section 8.2 of~\cite{BFGS12}, pp.~3406--3415). The theorem is derived from Proposition 8.3, and
the latter, for $\ell-2r=0$, from Lemmas 8.2 and 8.5--8.9 after the change of variable (269),
which in our case is the identity. The constants in these lemmas reduce to the constant $b_2$ of
Theorem 4.4 (Proposition 7.4), to the uniform bound of Proposition 4.8 of~\cite{Baldoma06}, where
the supremum is taken also over $\mu$, and to the geometry of the domains (Lemma 8.2). None of
them depends on $\eps$ or on $\mu\in B(\mu_0)$.
\item \emph{The exponent $\gamma$.} Proposition 8.3 is proved for $\gamma<\gamma_2$, where in our
case $\gamma_2=\tfrac12$, and Remark 8.4 there extends the result to all $\gamma\in(0,1)$. For
(E2) a single value of $\gamma$ suffices, and we take $\gamma\in[\tfrac13,\tfrac12)$. The proof
of Lemma 8.5 (p.~3411) uses the a priori bound
$|\zeta|^{\ell+1}|\partial_\zeta\psi^u-\partial_\zeta\psi_0^u|\le K$. For the pendulum and
$\gamma\ge\tfrac13$ this bound holds, because the expansion of $p_0$ at the point $ia$ has no
first-order correction.
\item \emph{What is taken from~\cite{Baldoma06}.} Theorem 4.12 is quoted in~\cite{BFGS12}
from~\cite{Baldoma06}, and the proof of Theorem 4.16 also uses the bound of Proposition 4.8
of~\cite{Baldoma06}, stated there in Fourier norms. Thus (E2) relies on that paper as well.
\end{enumerate}

Uniformity of the limit in $\zeta$, $\theta$ and $z$ in (E2) is neither claimed nor needed. In
Sections~\ref{ss:stokes} and~\ref{ss:inner} we obtain it from the estimates (E1) and Vitali's
theorem.

\paragraph{Relation to the normalization of Section~\ref{sec:borel}} If we write
$\psi_0^u=-1/\zeta-\tfrac14\psi+\mathrm{const}$, then equation~\eqref{eq:bfgs-inner} becomes
\begin{equation}\label{eq:bms-inner}
(\partial_\theta+\partial_\zeta)\psi=\tfrac18\,\zeta^2(\partial_\zeta\psi)^2-\frac{2zg(\theta)}{\zeta^2},
\end{equation}
and $\Phi_0^u=-4\zeta^{-2}+\partial_\zeta\psi$. This is representation~\eqref{eq:phi-def} with
$\varphi^u=\partial_\zeta\psi=-4z\,\partial_\zeta\overline\psi{}_0^{\,u}$, and
equation~\eqref{eq:inner} is the derivative of~\eqref{eq:bms-inner} with respect to $\zeta$. The
condition in~\eqref{eq:bfgs-psi0} for $\ell=2$ means $|\varphi^u|\le C|z|\,|\zeta|^{-3}$. In
Lemma~\ref{lem:S1} such an estimate is derived independently, from (E1). Of Theorem 4.12, only the
existence of the solution is used in this paper.

\subsection{What is used and what is not}\label{ss:A-unused}

The statements of~\cite{BFGS12} related to (E1) and (E2) are listed in
Table~\ref{tab:bfgs-used}. The proofs of this paper use them only through statement (E2). We do
not use the following
statements: Theorems 4.5--4.8 (continuation of the manifolds when $p_0$ has zeros); item 2 of
Theorem 4.12 and Proposition 4.15 (the difference of solutions of the inner equation and the
constants $\chi^{[k]}$); Theorems 4.21, 4.23 and their corollaries (asymptotics of the difference
$T^u-T^s$); the main Theorem 2.7. Nor do we use the statements on the stable manifold, because in
this paper it is obtained from the unstable one by reversibility. We study the difference of the
manifolds of the periodic problem anew in Sections~\ref{sec:straightening}
and~\ref{sec:stokes}, because we need estimates uniform in the complex amplitude and the first
variation with respect to the profile.

\begin{table}[ht]
\caption{Statements of~\cite{BFGS12} related to (E1) and (E2).}\label{tab:bfgs-used}
\centering\small
\renewcommand{\arraystretch}{1.15}
\begin{tabular}{@{}>{\raggedright\arraybackslash}p{0.30\linewidth}%
  >{\raggedright\arraybackslash}p{0.14\linewidth}%
  >{\raggedright\arraybackslash}p{0.48\linewidth}@{}}
\toprule
Statement of~\cite{BFGS12} & Page & Content \\
\midrule
Section 4.1; proof of Theorem 4.19 & 3330, 3349 & complex amplitude in a disc of arbitrary
radius; the constants depend neither on the amplitude nor on $\eps$ \\
Theorem 4.3, Proposition 6.4 & 3334, 3369 & the solution $T^u$ on $D^{\rm loc}$ and the weighted
estimate~\eqref{eq:bfgs-prop64}; cf.~(E1) \\
Theorem 4.4, domain (35) & 3327, 3335 & the continuation of $T^u$ to $D^{\rm out}$ and
estimate~\eqref{eq:bfgs-thm44}; cf.~(E1) \\
Theorem 4.12, item 1 (Theorem 3.1 of~\cite{Baldoma06}) & 3343 & the existence of the solution
$\psi_0^u$ of the inner equation \\
Theorem 4.16; formulas (36), (66), (67), (70) & 3328--3329, 3341, 3342, 3345 &
estimate~\eqref{eq:bfgs-thm416}; hence the limit (E2) \\
\bottomrule
\end{tabular}
\end{table}

\begin{remark}[agreement with the main theorem of~\cite{BFGS12}]\label{rem:bfgs-main}
Theorem 2.7 of~\cite{BFGS12} (p.~3319) gives, in the singular case $\ell-2r=0$, the area of the
lobe between two adjacent homoclinic points of the periodic problem. For a perturbation
independent of the momentum there is no logarithmic correction in the exponent, and the formula
reads
\[
\mathcal A=\frac{4|\mu|}{\eps}\,e^{-a/\eps}\Bigl(|f(\mu)|+O\bigl(1/|\ln\eps|\bigr)\Bigr),\qquad
f(\mu)=C_+^2\,\chi^{[-1]}(\mu),
\]
where $\chi^{[-1]}$ is the constant from item 2 of Theorem 4.12 and $f$ is an entire function
(formula (114) of~\cite{BFGS12}). In the notation of this paper this is
$\mathcal A=4L_h\bigl(|\chi[zg]|+o(1)\bigr)$, since the Stokes constant of
Definition~\ref{def:stokes} is related to the function $f$ by the equality $\chi[zg]=z\,f(z)$.
Indeed, the difference of the unstable and stable solutions of equation~\eqref{eq:bms-inner}
equals the difference of solutions of equation~\eqref{eq:bfgs-inner} multiplied by $-4=C_+^2$.
In~\cite{BFGS12} the factor $\mu$ is taken out of the latter, and the coefficient of the first
exponential in the former difference equals $\chi[zg]$ (Section~\ref{ss:A-bms}). This relation is
not used in the proofs.
\end{remark}

\subsection{Comparison with the work of Baldom\'a, M-Seara and Moreno}\label{ss:A-bms}

The paper~\cite{BMS24} studies the same periodic problem, since its Hamiltonian (1.1) is
$\tfrac12y^2+(\cos x-1)+\mu(\cos x-1)\,g(t/\eps)$, that is,~\eqref{eq:bfgs-our}. It is not needed
for the proofs of this paper, but the comparison makes it possible to check the normalizations.

\begin{itemize}
\item \emph{The inner equation.} Equation (1.16) of~\cite{BMS24} coincides
with~\eqref{eq:bms-inner}. By Theorem 2.3 of that paper (taken from~\cite{Baldoma06}) it has
solutions $\psi^\pm$ in domains of the form $\{|\Im\zeta|>\mp\vartheta\,\Re\zeta+\varrho\}$, whose
derivatives are uniquely determined by the estimate $|\partial_\zeta\psi^\pm|<M|\mu|/|\zeta|^3$.
The sign ``$-$'' corresponds to the unstable manifold. The sector $S^u$ is contained in such a
domain if $\tan\beta_0\le\vartheta\le\tan\beta_2$ and $\varrho\le\kappa_6$. Therefore the Taylor
coefficients in the amplitude of the functions $\varphi^u$ and $\partial_\zeta\psi^-$ coincide by
Lemma~\ref{lem:S2}, applied in the same way as in Section~\ref{ss:inner}, and
$\varphi^u=\partial_\zeta\psi^-$. Similarly $\varphi^s=\partial_\zeta\psi^+$.
\item \emph{The Stokes constant.} In~\cite{BMS24} the constants $\chi^{[k]}(\mu)$ are defined by
the expansion (2.5) of the difference $\psi^--\psi^+$ in the functions
$e^{ik(\zeta-\theta+\mu\mathfrak g)}$, $k<0$, where $\mathfrak g=O(1/\zeta)$. The first harmonic in
$\theta$ of this difference equals $\chi^{[-1]}(\mu)e^{-i\zeta}(1+o(1))$ as $\zeta\to-i\infty$.
Differentiating with respect to $\zeta$ and setting $\zeta=-ix$, we obtain
$\fc{\cJ}1(x;z)=-i\chi^{[-1]}(z)e^{-x}(1+o(1))$. Comparison with~\eqref{eq:chi-def} shows that
$\chi[zg]$ coincides with the function $\chi^{[-1]}(z)$ of~\cite{BMS24}.
\item \emph{The first Taylor coefficients.} By Theorem 2.4 of~\cite{BMS24}
\[
\chi_1=4\pi\,\fc g1,\qquad
\chi_2=-\frac{4\pi}3\sum_{k>1}\frac{\fc gk\,\fc g{1-k}}{k(1-k)} ,
\]
where the second formula is stated there for the case $\fc g1=0$. The first equality is
Proposition~\ref{prop:normalization}. The second agrees with Proposition~\ref{thm:S}. Indeed,
substituting $g=M+tp$ with $\fc p1=0$ and taking the coefficient of $t$, after the change
$k\mapsto1-k$ in one of the two sums one obtains
$-\tfrac{4\pi}3\sum_{k\ne0,1}\fc pk\fc M{1-k}/\bigl(k(1-k)\bigr)$, that is, the coefficient of
$\mu$ in Remark~\ref{rem:first-order-mu}. For the example $g=20\cos3\theta+16\cos2\theta$ from
Section 2.3.1 of~\cite{BMS24}, both formulas give $\chi_2=160\pi/9$.
\item \emph{A formula for the coefficients.} Formula (2.8) of~\cite{BMS24} expresses the first
nonzero Taylor coefficient of the function $\chi^{[-1]}$ through the limit as
$\zeta\to-i\infty$ of the derivative with respect to the amplitude of the first harmonic of the
difference $\psi^--\psi^+$, multiplied by $e^{i\zeta}$. Lemma~\ref{lem:S3} is an analogue of this
formula for all coefficients and for the derivative with respect to the profile.
\item \emph{The outer domain and the splitting.} Theorem 4.2 of~\cite{BMS24} restates Theorem 4.4
of~\cite{BFGS12} for the rhombus that is the intersection of the outer domains of the two
manifolds, and there the analyticity in $(u,\tau,\mu)$ and the fact that the disc of amplitudes is
complex are stated explicitly. Formula (4.6) of that paper reads, for real $u$,
\[
\partial_u\cD=2h^{-2}e^{-a/h}\bigl[\Im\bigl(\chi^{[-1]}e^{i(\theta-u/h)}\bigr)+o(1)\bigr].
\]
It agrees with item~4 of Lemma~\ref{lem:stokes}. Indeed, the transport equation gives
$\partial_\theta\cD\approx-h\,\partial_u\cD$, whence
$\partial_\theta^2\cD(0,0)\approx-2L_h\Re\chi^{[-1]}$.
\end{itemize}

\subsection{A direct proof of statement (E1)}\label{ss:A-direct}

In this section we prove Proposition~\ref{prop:E1}, that is, statement (E1). We obtain the
solution by one fixed point argument on the whole union $D=D^{\rm loc}\cup D^{\rm out}$. This is
possible for
two reasons: $\psi_0=2/\cosh u$ is analytic in the half-plane $\Re u<0$, and the rays going from a
point of $D$ to the left, or to the left at the angle $\pm45^\circ$, do not leave $D$. We invert
the operator $L=h^{-1}\partial_\theta+\partial_u$ by integration along these rays up to infinity,
namely along the horizontal ray for the zero harmonic in $\theta$, and along the rays parallel to
the sides of the wedge $D^{\rm out}$ for the harmonics with $k>0$ and $k<0$.

In this section $\sigma>0$ and $0<\varrho_1<\varrho_2$ are arbitrary, the domain $D^{\rm out}$ is
taken with a parameter $\kappa$ in place of $\kappa_1$, and $\kappa\ge1$, $\kappa h\le1$. The
letters $C_j$, $M_1$, $\mu_1$ denote constants that depend only on $\varrho_1$ and $\varrho_2$.

\paragraph{Weight and norm} We write $f=\sum_kf_k(u)e^{ik\theta}$ and
$g=\sum_{k\ne0}g_ke^{ik\theta}$, and set $\|g\|_*=\sum_{k\ne0}e^{|k|\sigma}|g_k|/|k|$. Let
\begin{equation}\label{eq:A-weight}
m(u)=|\psi_0(u)|^2\bigl(1+|\tanh u|\bigr),\qquad
\|f\|=\sum_ke^{|k|\sigma}\sup_{u\in D}\frac{|f_k(u)|}{m(u)},
\end{equation}
and let $X$ be the Banach space of functions $f$ with coefficients $f_k$ analytic on $D$ and
$\|f\|<\infty$. For $|\Im\theta|\le\sigma$ we have $|f(u,\theta)|\le m(u)\|f\|$, and every $f\in X$
is analytic on $D\times\TT_\sigma$. Let $r(u)=\min\{|u-ia|,|u+ia|\}$. The weight has the following
properties:
\begin{enumerate}[label=\textup{(m\arabic*)},leftmargin=*]
\item $m(s)\le M_1\max\{1,r(s)^{-3}\}$ on $D$, and $m(s)\le M_1e^{2\Re s}$ for $\Re s\le-\varrho_1$;
\item $m(u)\ge\mu_1\max\{1,r(u)^{-3}\}$ on $D^{\rm out}$, and $m(u)\ge\mu_1e^{2\Re u}$ for
$\Re u\le-\varrho_1$;
\item $r(u)\ge\kappa h/\sqrt2$ on $D^{\rm out}$, and $1+|\tanh u|\le C_3/(\kappa h)$ on $D$;
\item $1+|\tanh u|\le C_4|\cosh u|^2/|u^2+a^2|^3$ and $m(u)\le C_5|u^2+a^2|^{-3}$ on $D^{\rm out}$.
\end{enumerate}
Indeed, the zeros $i(\tfrac\pi2+n\pi)$ of $\cosh u$ do not belong to $D$. Near $\pm ia$, with
$w=u\mp ia$ and $|w|\le1$, we have $\cosh u=\pm i\sinh w$, $\sinh u=\pm i\cosh w$,
$0.82|w|\le|\sinh w|\le1.18|w|$ and $\cos1\le|\cosh w|\le\cosh1$. Hence $m\asymp r^{-3}$ and
$|\tanh u|\le1.9/r$ there. For $\Re u\le-\varrho_1$ the numbers $|\cosh u|$ and $|\sinh u|$ lie
between $\sinh|\Re u|$ and $\cosh\Re u$, so $m\asymp e^{2\Re u}$ and $|\tanh u|\le\coth\varrho_1$.
The rest of $D^{\rm out}$ lies in the compact set obtained from the closure of the wedge with
$\kappa h=0$ by removing the unit discs around $\pm ia$. On this set neither $\cosh u$ nor the
entire function $\cosh u/(u^2+a^2)$ vanishes (the zeros of the latter satisfy $|\Im u|\ge3a$).
Finally, $D^{\rm out}$ lies in the half-planes $\Re u\pm\Im u<a-\kappa h$, whose distance from
$\pm ia$ is $\kappa h/\sqrt2$.

\paragraph{Reduction} If $T=T_0+T_1$, then
$(\partial_uT)^2/(2\psi_0^2)=\psi_0^2/2+\partial_uT_1+(\partial_uT_1)^2/(2\psi_0^2)$, and
equation~\eqref{eq:HJ} becomes
\begin{equation}\label{eq:A-R}
LT_1=A-N(\partial_uT_1),\qquad A=\tfrac12zg(\theta)\psi_0^2,\qquad N(w)=\frac{w^2}{2\psi_0^2}.
\end{equation}
Let $d_0=1$, $d_k=(1-i)/\sqrt2$ for $k>0$, $d_k=(1+i)/\sqrt2$ for $k<0$, and
$\gamma_k(u)=\{u-d_kt:\ t\ge0\}$. Define the operators
\begin{equation}\label{eq:A-G}
\begin{gathered}
G[f]_k(u)=d_k\int_0^\infty e^{-ikd_kt/h}f_k(u-d_kt)\,dt,\\
B[f]=f-h^{-1}\partial_\theta G[f],
\end{gathered}
\end{equation}
so that $B[f]_k=f_k-(ik/h)G[f]_k$ and $B[f]_0=f_0$. For $k\ne0$ we have
$|e^{-ikd_kt/h}|=e^{-|k|t/(\sqrt2h)}$.

\begin{lemma}\label{lem:A-rays}
If $u\in D$, then $\gamma_k(u)\subset D$ for every $k$.
\end{lemma}

\begin{proof}
Along every ray $\Re s$ decreases, so the rays do not leave $D^{\rm loc}$. Let $u\in D^{\rm out}$.
The condition $|\Im s|<a-\kappa h-\Re s$ is the pair of inequalities $\Re s+\Im s<a-\kappa h$
and $\Re s-\Im s<a-\kappa h$. Along the horizontal ray both left-hand sides decrease. Along the ray
of direction $(-1+i)/\sqrt2$ the first one is constant and the second decreases, and along the ray
of direction $(-1-i)/\sqrt2$ it is the other way round. Hence the ray stays in $D^{\rm out}$ while
$\Re s>-\varrho_2$, and after that it lies in $D^{\rm loc}$, since $\varrho_2>\varrho_1$.
\end{proof}

\begin{lemma}\label{lem:A-weight}
For $u\in D$, $k\ne0$ and $t\ge0$ we have $m(u-d_kt)\le C_1m(u)\,e^{t/(2\sqrt2h)}$.
\end{lemma}

\begin{proof}
By complex conjugation it suffices to take $k>0$. Let $s=u+t(-1+i)/\sqrt2$. If $u\in D^{\rm loc}$,
then $\Re s=\Re u-t/\sqrt2$, and (m1), (m2) give $m(s)\le(M_1/\mu_1)m(u)$. Let $u\in D^{\rm out}$.
By (m1), $m(s)\le M_1\bigl(1+|s+ia|^{-3}+|s-ia|^{-3}\bigr)$, and by (m2)
$\max\{1,|u\pm ia|^{-3}\}\le m(u)/\mu_1$. While $s\in D^{\rm out}$, we have
$\tfrac{d}{dt}|s+ia|^2=\sqrt2(\Im s+a-\Re s)>0$, and afterwards $|s+ia|\ge|\Re s|>\varrho_1$.
Hence $|s+ia|\ge\min\{|u+ia|,\varrho_1\}$. The ray lies on the line
$\Re s+\Im s=\Re u+\Im u<a-\kappa h$, so $|s-ia|\ge\kappa h/\sqrt2$, and also
$|s-ia|\ge|u-ia|-t$. If $t\le|u-ia|/2$, then $|s-ia|^{-3}\le8|u-ia|^{-3}$. If $t>|u-ia|/2$, then
with $x=|u-ia|/h$ and $C_2=\sup_{x>0}x^3e^{-x/(4\sqrt2)}$
\[
|s-ia|^{-3}\le\Bigl(\frac{\sqrt2\,x}{\kappa}\Bigr)^3|u-ia|^{-3}\le
2\sqrt2\,C_2\,e^{t/(2\sqrt2h)}\,|u-ia|^{-3}.
\]
Together these bounds give the lemma.
\end{proof}

\begin{lemma}\label{lem:A-G}
Let $f\in X$.
\begin{enumerate}[label=\textup{\arabic*.}]
\item For $k\ne0$ we have $|G[f]_k(u)|\le C_P\,(h/|k|)\,m(u)\sup_D|f_k/m|$, where
$C_P=2\sqrt2C_1$. Hence $\|h^{-1}\partial_\theta G[f]\|\le C_P\|f\|$ and
$\|B[f]\|\le(1+C_P)\|f\|$.
\item The coefficients $G[f]_k$ are analytic on $D$, $LG[f]=f$, and therefore
$\partial_uG[f]=B[f]$.
\item If also $\partial_uf\in X$, then $\partial_uG[f]=G[\partial_uf]$.
\end{enumerate}
\end{lemma}

\begin{proof}
Item~1 follows from Lemma~\ref{lem:A-weight}, since $|k|-\tfrac12\ge|k|/2$. For item~2 write
$G[f]_k(u)=e^{-iku/h}\Psi(u)$, where $\Psi(u)$ is the integral of $Y(s)=e^{iks/h}f_k(s)$ along
$\gamma_k(u)$ from infinity to $u$. For $u'$ close to $u$, the half-strip bounded by
$\gamma_k(u)$, $\gamma_k(u')$ and the segment $[u,u']$ consists of rays starting on $[u,u']$, so it
lies in $D$ by Lemma~\ref{lem:A-rays}. Far along the half-strip $Y$ is exponentially small,
because of the factor $e^{iks/h}$ for $k\ne0$ and by (m1) for $k=0$. By Cauchy's theorem
$\Psi(u')-\Psi(u)=\int_{[u,u']}Y\,ds$, hence $\partial_uG[f]_k=f_k-(ik/h)G[f]_k$. By item~1 the
Fourier series of $G[f]$ and of $h^{-1}\partial_\theta G[f]$ converge absolutely and uniformly on
compact subsets of $D\times\TT_\sigma$, so $LG[f]=f$. Item~3 is differentiation under the
integral sign, which is justified by Lemma~\ref{lem:A-weight} and (m1).
\end{proof}

\begin{lemma}\label{lem:A-fixed}
Let $R=R_0+1$, $\beta=(1+C_P)C_3$ and $\kappa\ge\kappa_0:=8\beta C_PR\|g\|_*$. Then for $|z|<R$
the equation
\begin{equation}\label{eq:A-fixed}
w=w_0-B[N(w)],\qquad w_0=G[\partial_uA],
\end{equation}
has a solution $w\in X$, analytic in $(u,\theta,z)$, with $\|w\|\le2C_P\|g\|_*|z|h$.
\end{lemma}

\begin{proof}
Since $\partial_uA=-zg\psi_0^2\tanh u$, we have $|(\partial_uA)_k|\le|z||g_k|\,m$ and
$(\partial_uA)_0=0$. By Lemma~\ref{lem:A-G}, $\|w_0\|\le C_P\|g\|_*|z|h$. The product of Fourier
series with the weights $e^{|k|\sigma}$ is submultiplicative, so by (m3)
\begin{gather*}
\|N(w)\|\le\frac{C_3}{2\kappa h}\|w\|^2,\\
\|N(w)-N(w')\|\le\frac{C_3}{2\kappa h}\bigl(\|w\|+\|w'\|\bigr)\|w-w'\| .
\end{gather*}
Together with Lemma~\ref{lem:A-G} this shows that on the ball $\|w\|\le2C_PR\|g\|_*h$ the
right-hand side of~\eqref{eq:A-fixed} maps the ball into itself with Lipschitz constant at most
$\tfrac14$. Its fixed point lies in the ball $\|w\|\le2\|w_0\|$ as well. The iterations starting
from $w=0$ are polynomials in $z$ and converge uniformly for $|z|<R$, hence the solution is
analytic in $z$.
\end{proof}

\begin{proof}[Proof of statement (E1)]
Take $\kappa=\kappa_1:=\max\{1,\kappa_0\}$ and $h<h_1:=1/\kappa_1$, let $w$ be the solution of
Lemma~\ref{lem:A-fixed}, and set $T_1=G[A-N(w)]$. By Lemma~\ref{lem:A-G},
$\partial_uT_1=G[\partial_uA]-B[N(w)]=w$ and $LT_1=A-N(\partial_uT_1)$. Thus $T^u=T_0+T_1$
solves~\eqref{eq:HJ}. For fixed $(u,\theta)$ the map $f\mapsto G[f](u,\theta)$ is a bounded linear
functional on $X$. Indeed, for $k=0$ its norm is at most $\int_0^\infty m(u-t)\,dt$, and for
$k\ne0$ one uses Lemma~\ref{lem:A-G}. Hence $T_1$ is analytic in $z$. It is analytic in
$(u,\theta)$ by Lemma~\ref{lem:A-G} and locally bounded, so it is jointly analytic by Hartogs'
theorem. Finally, $E^u=w/\psi_0^2$ and
\[
|E^u|\le\bigl(1+|\tanh u|\bigr)\|w\|\le2C_P\|g\|_*\,|z|h\,\bigl(1+|\tanh u|\bigr).
\]
On $D^{\rm loc}$ we have $1+|\tanh u|\le1+\coth\varrho_1$, and on $D^{\rm out}$ estimate (m4)
applies. This gives~\eqref{eq:E1}, even for all $|z|<R_0+1$, with
$C_E=2C_P\|g\|_*\max\{1+\coth\varrho_1,C_4\}$. The constants $\kappa_1$, $h_1$, $C_E$ depend on
$\varrho_2$ and are chosen after it, as in Section~\ref{ss:constants}.
\end{proof}

\paragraph{Relation to~\cite{BFGS12}} The graph of $T^u$ is the unstable manifold. Indeed, the
proof of item~4 of Lemma~\ref{lem:local} uses only that $E^u$ comes from a solution
of~\eqref{eq:HJ} that is analytic on $D^{\rm loc}\times\TT_\sigma$ and satisfies the first
estimate~\eqref{eq:E1}. The function given by Theorems 4.3 and 4.4 of~\cite{BFGS12} is also such
a solution, and it satisfies the first estimate~\eqref{eq:E1} by Section~\ref{ss:A-E1}. We take
$\varrho_1=\rho_1$, $\varrho_2=\rho_2$ and the same $\sigma$ as there.
Proposition~\ref{prop:E1} remains valid when $\kappa_1$ is increased, so we may assume that the
constant $\kappa_1$ is the same for both functions. For $|z|\le R_0$ and small $h$,
equality~\eqref{eq:coincidence} holds for both functions with the same right-hand side. Hence
the two functions $E^u$, and with them the derivatives $\partial_uT^u=\psi_0^2(1+E^u)$, coincide
for $\Re u<-\varrho_1'$. Both derivatives are analytic in $u$ on the connected domain
$D^{\rm loc}\cup D^{\rm out}$, and therefore they coincide on it. Hence statement (E2), derived
from Theorem 4.16 of~\cite{BFGS12}, applies to $T^u$.
By (m4) we also have $\sum_ke^{|k|\sigma}\sup_{D^{\rm out}}|(u^2+a^2)^3w_k|\le C_5\|w\|$. For the
pendulum and $\beta_1=\pi/4$ this is the estimate of Proposition 7.4 of~\cite{BFGS12} in the norm
of its Section 7.1.1. In summary, the proof of (E1) given in this section does not use the
results of~\cite{BFGS12}, while statement (E2) relies on Theorem 4.16 of~\cite{BFGS12} and on
Theorem 3.1 and Proposition 4.8 of~\cite{Baldoma06} (Section~\ref{ss:A-E2}).
\section{Certificate for Theorem 3: data and reproduction}\label{app:certificate}

The computer-assisted part of the proof of Theorem~\ref{thm:cos2} is Lemma~\ref{lem:certificate}.
It states that $37$ explicitly given polynomials with rational coefficients are nonnegative or
positive on the interval $[0,1]$, and gives the exact value of one integral. The other steps of
Section~\ref{sec:negative} do not depend on the computer. In this appendix we write out the
inequalities explicitly, describe the two programs by which they were verified, and explain how to
repeat the verification. The whole verification is carried out in exact rational arithmetic and
involves no floating-point numbers.

\subsection{Data}\label{ss:C-data}

The input data are the $52$ integers of Table~\ref{tab:certificate}, namely the numerators of the
coefficients of the thirteen cubic polynomials $g^D_i$, $g^N_i$ ($i=0,\dots,5$) and $g^S$ with the
common denominator $10^4$. How these numbers were found is described in
Section~\ref{ss:C-origin}. This plays no role in the verification.

Everything else is given by the formulas of Section~\ref{sec:negative}: the numbers
$r_+=\tfrac{141422}{100000}$, $r_-=\tfrac{141421}{100000}$; the frequencies
$(k_0,\dots,k_5)=(12,2,4,6,-4,-2)$; the exponents $p_i=1$ for even $i$ and $p_i=2$ otherwise; the
majorants~\eqref{eq:poly-majorants}; the upper functions~\eqref{eq:poly-upper} and the functions
$\bar W_i$, $\bar U_i$, $\bar R$ composed of them by~\eqref{eq:upper-functions}. All of them are
polynomials. The polynomials $\bar D_i$ and $\bar N_i$ have degree $4$ for even $i$ and $5$
otherwise, $\bar S$ has degree $4$, and $\bar U_0,\dots,\bar U_5$ have degrees $4$, $7$, $4$,
$5$, $4$, $6$. For example,
\begin{equation}\label{eq:C-U1}
\begin{aligned}
\bar U_1(s)=\frac1{12\,800\,000}\bigl(&4525504+4525504\,s+6569808\,s^2+7067872\,s^3\\
&-5929140\,s^4+13935252\,s^5+494977\,s^6+636399\,s^7\bigr),
\end{aligned}
\end{equation}
and termwise integration gives $\int_0^1\bar U_1\,ds=\tfrac{96257079}{102400000}$, that is,
equality~\eqref{eq:m-prime}.

\subsection{The list of inequalities}\label{ss:C-list}

Recall that $1*a$ denotes $\int_0^sa(t)\,dt$, and $a*b$ denotes the convolution
$\int_0^sa(t)\,b(s-t)\,dt$. For powers, $s^x*s^y=x!\,y!\,s^{x+y+1}/(x+y+1)!$. Class indices are
added modulo $6$. Set
\begin{align}
\mathsf M_1&=\hat B_1\,(2-s)^2-r_+, &
\mathsf M_d&=s^{-2}\bigl[2\,(4-s^2)^2\,\hat d-r_+\,s^2(12-s^2)\bigr],\notag\\*
\mathsf M_5&=\hat B_5\,(2+s)^2-r_+, &
\mathsf M_{d-}&=s^{-2}\bigl[r_-\,s^2(12-s^2)-2\,(4-s^2)^2\,\hat d_-\bigr],\label{eq:C-M}
\end{align}
\begin{equation}\label{eq:C-Q}
\begin{aligned}
\mathsf Q^D_i&=s^{-p_i}\Bigl[(k_i-s)^2\,\bar D_i-r_+\,1*\bigl(\bar U_{i-1}+\bar U_{i+1}\bigr)\Bigr],\\
\mathsf Q^N_i&=s^{-p_i-1}\Bigl[s\,\bar N_i-\sum_{h=0}^{5}(s\bar W_h)*\bar U_{i-h}\Bigr],\\
\mathsf Q^S&=s^{-3}\Bigl[s^2\,\bar S-r_+\,1*\bigl(\hat d+\bar R\bigr)\Bigr],
\end{aligned}
\end{equation}
\begin{equation}\label{eq:C-T}
\begin{gathered}
\mathsf T_i=s^{1-p_i}\Bigl[\frac{r_-}2\bigl(\bar U_{i-1}+\bar U_{i+1}\bigr)-1*\bar U_i\Bigr],\\
\mathsf T_S=s^{-2}\Bigl[\frac{r_-}2\bigl(\hat d_-+\bar R\bigr)-1*\bar U_0\Bigr],
\end{gathered}
\end{equation}
where $i=0,\dots,5$. Each expression in square brackets is a polynomial divisible by the indicated
power of $s$, so all $37$ functions are polynomials. The inequalities of
Lemma~\ref{lem:certificate} are equivalent to the following statements on the interval $[0,1]$:
\begin{enumerate}[label=\textup{(\alph*)}]
\item $\mathsf M_1\ge0$, $\mathsf M_5\ge0$, $\mathsf M_d\ge0$, $\mathsf M_{d-}\ge0$ (four
inequalities for the majorants);
\item $g^D_i>0$, $g^N_i>0$ for $i=0,\dots,5$, and $g^S>0$ (thirteen inequalities);
\item $\mathsf Q^D_i>0$, $\mathsf Q^N_i>0$ for $i=0,\dots,5$, and $\mathsf Q^S>0$ (thirteen
inequalities; these are conditions~\eqref{eq:cond-I});
\item $\mathsf T_i\ge0$ for $i=0,\dots,5$, and $\mathsf T_S\ge0$ (seven inequalities; these are
conditions~\eqref{eq:cond-II}).
\end{enumerate}
Conditions~\eqref{eq:cond-I} and~\eqref{eq:cond-II} are obtained from (c) and (d) by multiplication
by a power of $s$, and the first of conditions~\eqref{eq:cond-I} also by division by
$(k_i-s)^2>0$. The inequalities (a) also follow from the factorizations~\eqref{eq:factorizations},
which can be checked by hand:
\[
\mathsf M_1=\tfrac{r_+}{256}\,s^7(1-s)(20-9s),\qquad
\mathsf M_d=\tfrac{r_+}{2304}\,s^6(12-s^2)(1-s^2)(64-13s^2),
\]
and similarly for the other two. As an example, we write out one of the inequalities (c) in full:
\begin{equation}\label{eq:C-QD1}
\begin{aligned}
\mathsf Q^D_1(s)=\frac1{5\cdot10^{9}}\bigl(&27297705+3279504430\,s-16928204285\,s^2\\
&+30755290926\,s^3-21695500000\,s^4\\
&+4710500000\,s^5\bigr)>0\quad\text{on }[0,1].
\end{aligned}
\end{equation}

The degrees of the $37$ polynomials do not exceed $14$.

\subsection{Two verification programs}\label{ss:C-programs}

The inequalities were verified by two programs written independently of each other and using
different positivity criteria.

\paragraph{Bernstein coefficients} The first program uses only the Python standard library (the
class \texttt{Fraction}). Polynomials are stored as lists of coefficients. Addition,
multiplication, the operation $1*$ and convolution are performed exactly, and in a division by
$s^q$ it is checked that the lowest $q$ coefficients vanish. Positivity of a polynomial
$a(s)=\sum_{j\le n}a_js^j$ on $[0,1]$ is established as follows. We write it in the Bernstein
basis:
\[
a(s)=\sum_{k=0}^{n}b_k\binom nk s^k(1-s)^{n-k},\qquad
b_k=\sum_{j=0}^{k}\binom kj\binom nj^{-1}a_j .
\]
The basis polynomials are nonnegative on $[0,1]$ and sum to one. Therefore, if all $b_k>0$, then
$a>0$ on $[0,1]$, and if all $b_k\ge0$, then $a\ge0$. If the criterion fails while the numbers
$b_0=a(0)$ and $b_n=a(1)$ are nonnegative, the interval is bisected, and the Bernstein
coefficients on the halves are computed exactly by de Casteljau's algorithm. At most twelve
bisections are allowed. A positive answer of the program proves the inequality, while a negative
answer only means that the criterion has failed.

\paragraph{The Sturm method} The second program was written separately, does not use the code of
the first one and relies on the library \texttt{sympy}. The polynomials are built by symbolic
integration, and the inequality $q>0$ on $[0,1]$ is checked as the pair of conditions
$q(\tfrac12)>0$ and ``the number of real roots of $q$ on $[0,1]$ is zero''. The roots are counted
by Sturm's method. For the four inequalities (a), the factors $s$ and $1-s$ corresponding to the
roots at the endpoints are first taken out. The seven inequalities (d) are checked by this program
as strict ones, that is, in a stronger form than required.

\subsection{Results and test modifications}\label{ss:C-results}

Both programs report that all $37$ inequalities hold, and both give the same value
$\int_0^1\bar U_1\,ds=\tfrac{96257079}{102400000}$. Table~\ref{tab:certificate} in the text of the
paper was compared automatically with the data used by the programs.

To make sure that the programs really distinguish correct tables from incorrect ones, we introduced
deliberate errors into the table. Both programs reject each of
the five modifications listed in Table~\ref{tab:C-mutations} and report the same violated
condition.

\begin{table}[ht]
\caption{Test modifications of Table~\ref{tab:certificate} and the response of both
programs.}\label{tab:C-mutations}
\centering\small
\renewcommand{\arraystretch}{1.12}
\begin{tabular}{@{}l@{\qquad}l@{\qquad}l@{}}
\toprule
Modified coefficient & Change & Violated condition \\
\midrule
constant term of $g^S$ & $5132\to4500$ & $\mathsf Q^S>0$ \\
constant term of $g^D_1$ & $1334\to1100$ & $\mathsf Q^D_1>0$ \\
constant term of $g^N_0$ & $1302\to1200$ & $\mathsf Q^N_0>0$ \\
coefficient of $s^3$ in $g^D_4$ & $-56\to-400$ & $g^D_4>0$ and $\mathsf Q^D_4>0$ \\
coefficient of $s^2$ in $g^N_2$ & $-1095\to-1400$ & $\mathsf Q^N_2>0$ \\
\bottomrule
\end{tabular}
\end{table}

\subsection{How to repeat the verification}\label{ss:C-repro}

The programs are available in the repository~\cite{GulkanovCode} (see the code availability
statement). The first program runs for about $0.2$ seconds, the second for about $5$ seconds, and
the test with deliberate errors for about $30$ seconds. We used Python~3.10.0 and sympy~1.13.3.

Statements (a)--(d) do not depend on these programs. They can be checked in any computer algebra
system. This requires addition, multiplication and integration of polynomials with rational
coefficients and isolation of the real roots of a polynomial of degree at most $14$.

What is taken on trust here is the correctness of the integer and rational arithmetic of Python,
on which the rigor of the verification relies, and, for the second program, also the correctness
of the library \texttt{sympy}. That the inequalities (a)--(d) suffice for Theorem~\ref{thm:cos2}
is proved in Section~\ref{sec:negative} without a computer.

\subsection{How the polynomials were found}\label{ss:C-origin}

This section is not part of the proof. The polynomials of Table~\ref{tab:certificate} were
obtained as follows. Consider conditions~\eqref{eq:cond-I} with the
majorants~\eqref{eq:poly-majorants}, with equality signs and with the right-hand sides multiplied
by $1.05$. The smallest solution of this system is computed by power series. Then the functions
$D_i/s^{p_i}$, $N_i/s^{p_i}$ and $S/s$ are fitted by cubic polynomials at Chebyshev nodes, and the
coefficients are rounded to four decimal places. The program that carries out this construction
is available together with the verification programs.
\end{appendices}

\bibliography{references}


\begin{thebibliography}{31}
\ifx \bisbn   \undefined \def \bisbn  #1{ISBN #1}\fi
\ifx \binits  \undefined \def \binits#1{#1}\fi
\ifx \bauthor  \undefined \def \bauthor#1{#1}\fi
\ifx \batitle  \undefined \def \batitle#1{#1}\fi
\ifx \bjtitle  \undefined \def \bjtitle#1{#1}\fi
\ifx \bvolume  \undefined \def \bvolume#1{\textbf{#1}}\fi
\ifx \byear  \undefined \def \byear#1{#1}\fi
\ifx \bissue  \undefined \def \bissue#1{#1}\fi
\ifx \bfpage  \undefined \def \bfpage#1{#1}\fi
\ifx \blpage  \undefined \def \blpage #1{#1}\fi
\ifx \burl  \undefined \def \burl#1{\textsf{#1}}\fi
\ifx \doiurl  \undefined \def \doiurl#1{\url{https://doi.org/#1}}\fi
\ifx \betal  \undefined \def \betal{\textit{et al.}}\fi
\ifx \binstitute  \undefined \def \binstitute#1{#1}\fi
\ifx \binstitutionaled  \undefined \def \binstitutionaled#1{#1}\fi
\ifx \bctitle  \undefined \def \bctitle#1{#1}\fi
\ifx \beditor  \undefined \def \beditor#1{#1}\fi
\ifx \bpublisher  \undefined \def \bpublisher#1{#1}\fi
\ifx \bbtitle  \undefined \def \bbtitle#1{#1}\fi
\ifx \bedition  \undefined \def \bedition#1{#1}\fi
\ifx \bseriesno  \undefined \def \bseriesno#1{#1}\fi
\ifx \blocation  \undefined \def \blocation#1{#1}\fi
\ifx \bsertitle  \undefined \def \bsertitle#1{#1}\fi
\ifx \bsnm \undefined \def \bsnm#1{#1}\fi
\ifx \bsuffix \undefined \def \bsuffix#1{#1}\fi
\ifx \bparticle \undefined \def \bparticle#1{#1}\fi
\ifx \barticle \undefined \def \barticle#1{#1}\fi
\bibcommenthead
\ifx \bconfdate \undefined \def \bconfdate #1{#1}\fi
\ifx \botherref \undefined \def \botherref #1{#1}\fi
\ifx \url \undefined \def \url#1{\textsf{#1}}\fi
\ifx \bchapter \undefined \def \bchapter#1{#1}\fi
\ifx \bbook \undefined \def \bbook#1{#1}\fi
\ifx \bcomment \undefined \def \bcomment#1{#1}\fi
\ifx \oauthor \undefined \def \oauthor#1{#1}\fi
\ifx \citeauthoryear \undefined \def \citeauthoryear#1{#1}\fi
\ifx \endbibitem  \undefined \def \endbibitem {}\fi
\ifx \bconflocation  \undefined \def \bconflocation#1{#1}\fi
\ifx \arxivurl  \undefined \def \arxivurl#1{\textsf{#1}}\fi
\csname PreBibitemsHook\endcsname

\bibitem[\protect\citeauthoryear{Gallavotti}{}]{Gallavotti98}
\begin{botherref}
\oauthor{\bsnm{Gallavotti}, \binits{G.}}:
Separatrix splitting for a quasi periodically forced pendulum.
In: Aizenman, M. (ed.) Some Problems in Mathematical Physics, problem 9805.
  \urlGallavottiProblem\ (1998). Accessed 17 November 2026
\end{botherref}
\endbibitem

\bibitem[\protect\citeauthoryear{Gallavotti et~al.}{1999}]{GGM99}
\begin{barticle}
\bauthor{\bsnm{Gallavotti}, \binits{G.}},
\bauthor{\bsnm{Gentile}, \binits{G.}},
\bauthor{\bsnm{Mastropietro}, \binits{V.}}:
\batitle{Separatrix splitting for systems with three time scales}.
\bjtitle{Commun. Math. Phys.}
\bvolume{202}(\bissue{1}),
\bfpage{197}--\blpage{236}
(\byear{1999})
\doiurl{10.1007/s002200050579} .
\bcomment{Preprint version: arXiv:chao-dyn/9709004, titled ``Pendulum:
  separatrix splitting''}
\end{barticle}
\endbibitem

\bibitem[\protect\citeauthoryear{Baldom\'a et~al.}{2012}]{BFGS12}
\begin{barticle}
\bauthor{\bsnm{Baldom\'a}, \binits{I.}},
\bauthor{\bsnm{Fontich}, \binits{E.}},
\bauthor{\bsnm{Guardia}, \binits{M.}},
\bauthor{\bsnm{Seara}, \binits{T.M.}}:
\batitle{Exponentially small splitting of separatrices beyond {M}elnikov
  analysis: rigorous results}.
\bjtitle{J. Differ. Equ.}
\bvolume{253}(\bissue{12}),
\bfpage{3304}--\blpage{3439}
(\byear{2012})
\doiurl{10.1016/j.jde.2012.09.003}
{\href{https://arxiv.org/abs/1201.5152}{{arXiv:1201.5152}}}
\end{barticle}
\endbibitem

\bibitem[\protect\citeauthoryear{Chierchia and Gallavotti}{1994}]{CG94}
\begin{barticle}
\bauthor{\bsnm{Chierchia}, \binits{L.}},
\bauthor{\bsnm{Gallavotti}, \binits{G.}}:
\batitle{Drift and diffusion in phase space}.
\bjtitle{Ann. Inst. Henri Poincar\'e Phys. Th\'eor.}
\bvolume{60}(\bissue{1}),
\bfpage{1}--\blpage{144}
(\byear{1994})
\end{barticle}
\endbibitem

\bibitem[\protect\citeauthoryear{Chierchia and Gallavotti}{1998}]{CG98}
\begin{barticle}
\bauthor{\bsnm{Chierchia}, \binits{L.}},
\bauthor{\bsnm{Gallavotti}, \binits{G.}}:
\batitle{Erratum: ``{D}rift and diffusion in phase space''}.
\bjtitle{Ann. Inst. Henri Poincar\'e Phys. Th\'eor.}
\bvolume{68}(\bissue{1}),
\bfpage{135}
(\byear{1998})
\end{barticle}
\endbibitem

\bibitem[\protect\citeauthoryear{Gallavotti}{1994}]{Gallavotti94}
\begin{barticle}
\bauthor{\bsnm{Gallavotti}, \binits{G.}}:
\batitle{Twistless {KAM} tori, quasi flat homoclinic intersections, and other
  cancellations in the perturbation series of certain completely integrable
  {H}amiltonian systems. {A} review}.
\bjtitle{Rev. Math. Phys.}
\bvolume{6}(\bissue{3}),
\bfpage{343}--\blpage{411}
(\byear{1994})
\doiurl{10.1142/S0129055X9400016X}
\end{barticle}
\endbibitem

\bibitem[\protect\citeauthoryear{Gallavotti et~al.}{1999}]{GGM99b}
\begin{barticle}
\bauthor{\bsnm{Gallavotti}, \binits{G.}},
\bauthor{\bsnm{Gentile}, \binits{G.}},
\bauthor{\bsnm{Mastropietro}, \binits{V.}}:
\batitle{A field theory approach to {L}indstedt series for hyperbolic tori in
  three time scales problems}.
\bjtitle{J. Math. Phys.}
\bvolume{40}(\bissue{12}),
\bfpage{6430}--\blpage{6472}
(\byear{1999})
\doiurl{10.1063/1.533101}
{\href{https://arxiv.org/abs/chao-dyn/9812004}{{arXiv:chao-dyn/9812004}}}
\end{barticle}
\endbibitem

\bibitem[\protect\citeauthoryear{Gallavotti et~al.}{2000}]{GGM00}
\begin{barticle}
\bauthor{\bsnm{Gallavotti}, \binits{G.}},
\bauthor{\bsnm{Gentile}, \binits{G.}},
\bauthor{\bsnm{Mastropietro}, \binits{V.}}:
\batitle{{H}amilton--{J}acobi equation, heteroclinic chains and {A}rnol'd
  diffusion in three time scale systems}.
\bjtitle{Nonlinearity}
\bvolume{13}(\bissue{2}),
\bfpage{323}--\blpage{340}
(\byear{2000})
\doiurl{10.1088/0951-7715/13/2/301}
\end{barticle}
\endbibitem

\bibitem[\protect\citeauthoryear{Gallavotti}{1998}]{GallavottiMPAG}
\begin{barticle}
\bauthor{\bsnm{Gallavotti}, \binits{G.}}:
\batitle{Arnold's diffusion in isochronous systems}.
\bjtitle{Math. Phys. Anal. Geom.}
\bvolume{1}(\bissue{4}),
\bfpage{295}--\blpage{312}
(\byear{1998})
\doiurl{10.1023/A:1009893118532}
{\href{https://arxiv.org/abs/chao-dyn/9709011}{{arXiv:chao-dyn/9709011}}}
\end{barticle}
\endbibitem

\bibitem[\protect\citeauthoryear{Berti and Bolle}{2002}]{BertiBolle02}
\begin{barticle}
\bauthor{\bsnm{Berti}, \binits{M.}},
\bauthor{\bsnm{Bolle}, \binits{P.}}:
\batitle{Fast {A}rnold diffusion in systems with three time scales}.
\bjtitle{Discrete Contin. Dyn. Syst.}
\bvolume{8}(\bissue{3}),
\bfpage{795}--\blpage{811}
(\byear{2002})
\doiurl{10.3934/dcds.2002.8.795}
{\href{https://arxiv.org/abs/math/0103065}{{arXiv:math/0103065}}}
\end{barticle}
\endbibitem

\bibitem[\protect\citeauthoryear{Procesi}{2003}]{Procesi03}
\begin{barticle}
\bauthor{\bsnm{Procesi}, \binits{M.}}:
\batitle{Exponentially small splitting and {A}rnold diffusion for multiple time
  scale systems}.
\bjtitle{Rev. Math. Phys.}
\bvolume{15}(\bissue{4}),
\bfpage{339}--\blpage{386}
(\byear{2003})
\doiurl{10.1142/S0129055X03001655}
\end{barticle}
\endbibitem

\bibitem[\protect\citeauthoryear{Mastropietro}{2001}]{Mastropietro01}
\begin{barticle}
\bauthor{\bsnm{Mastropietro}, \binits{V.}}:
\batitle{Arnold diffusion and the d'{A}lembert precession problem}.
\bjtitle{Regul. Chaotic Dyn.}
\bvolume{6}(\bissue{4}),
\bfpage{355}--\blpage{375}
(\byear{2001})
\doiurl{10.1070/RD2001v006n04ABEH000183}
\end{barticle}
\endbibitem

\bibitem[\protect\citeauthoryear{Sim\'o}{1994}]{Simo94}
\begin{bchapter}
\bauthor{\bsnm{Sim\'o}, \binits{C.}}:
\bctitle{Averaging under fast quasiperiodic forcing}.
In: \beditor{\bsnm{Seimenis}, \binits{J.}} (ed.)
\bbtitle{Hamiltonian Mechanics: Integrability and Chaotic Behavior}.
\bsertitle{NATO ASI Series B: Physics},
vol. \bseriesno{331},
pp. \bfpage{13}--\blpage{34}.
\bpublisher{Plenum Press},
\blocation{New York}
(\byear{1994}).
\doiurl{10.1007/978-1-4899-0964-0_2}
\end{bchapter}
\endbibitem

\bibitem[\protect\citeauthoryear{Delshams et~al.}{1997}]{DGJS97}
\begin{barticle}
\bauthor{\bsnm{Delshams}, \binits{A.}},
\bauthor{\bsnm{Gelfreich}, \binits{V.}},
\bauthor{\bsnm{Jorba}, \binits{{\`A}.}},
\bauthor{\bsnm{Seara}, \binits{T.M.}}:
\batitle{Exponentially small splitting of separatrices under fast quasiperiodic
  forcing}.
\bjtitle{Commun. Math. Phys.}
\bvolume{189}(\bissue{1}),
\bfpage{35}--\blpage{71}
(\byear{1997})
\doiurl{10.1007/s002200050190}
\end{barticle}
\endbibitem

\bibitem[\protect\citeauthoryear{Delshams and Seara}{1992}]{DS92}
\begin{barticle}
\bauthor{\bsnm{Delshams}, \binits{A.}},
\bauthor{\bsnm{Seara}, \binits{T.M.}}:
\batitle{An asymptotic expression for the splitting of separatrices of the
  rapidly forced pendulum}.
\bjtitle{Commun. Math. Phys.}
\bvolume{150}(\bissue{3}),
\bfpage{433}--\blpage{463}
(\byear{1992})
\doiurl{10.1007/BF02096956}
\end{barticle}
\endbibitem

\bibitem[\protect\citeauthoryear{Gelfreich and Lazutkin}{2001}]{GL01}
\begin{barticle}
\bauthor{\bsnm{Gelfreich}, \binits{V.G.}},
\bauthor{\bsnm{Lazutkin}, \binits{V.F.}}:
\batitle{Splitting of separatrices: perturbation theory and exponential
  smallness}.
\bjtitle{Russ. Math. Surv.}
\bvolume{56}(\bissue{3}),
\bfpage{499}--\blpage{558}
(\byear{2001})
\doiurl{10.1070/RM2001v056n03ABEH000394}
\end{barticle}
\endbibitem

\bibitem[\protect\citeauthoryear{Lazutkin}{2005}]{Lazutkin84}
\begin{barticle}
\bauthor{\bsnm{Lazutkin}, \binits{V.F.}}:
\batitle{Splitting of separatrices for the {C}hirikov standard map}.
\bjtitle{J. Math. Sci.}
\bvolume{128}(\bissue{2}),
\bfpage{2687}--\blpage{2705}
(\byear{2005})
\doiurl{10.1007/s10958-005-0219-7} .
\bcomment{Translated from Zap. Nauchn. Semin. POMI \textbf{300}, 25--55 (2003);
  the original manuscript dates from 1984}
\end{barticle}
\endbibitem

\bibitem[\protect\citeauthoryear{Gelfreich}{1999}]{Gelfreich99}
\begin{barticle}
\bauthor{\bsnm{Gelfreich}, \binits{V.G.}}:
\batitle{A proof of the exponentially small transversality of the separatrices
  for the standard map}.
\bjtitle{Commun. Math. Phys.}
\bvolume{201}(\bissue{1}),
\bfpage{155}--\blpage{216}
(\byear{1999})
\doiurl{10.1007/s002200050553}
\end{barticle}
\endbibitem

\bibitem[\protect\citeauthoryear{Treschev}{1996}]{Treschev96}
\begin{barticle}
\bauthor{\bsnm{Treschev}, \binits{D.V.}}:
\batitle{An averaging method for {H}amiltonian systems, exponentially close to
  integrable ones}.
\bjtitle{Chaos}
\bvolume{6}(\bissue{1}),
\bfpage{6}--\blpage{14}
(\byear{1996})
\doiurl{10.1063/1.166149}
\end{barticle}
\endbibitem

\bibitem[\protect\citeauthoryear{Gelfreich}{1997}]{Gelfreich97}
\begin{barticle}
\bauthor{\bsnm{Gelfreich}, \binits{V.G.}}:
\batitle{Reference systems for splitting of separatrices}.
\bjtitle{Nonlinearity}
\bvolume{10}(\bissue{1}),
\bfpage{175}--\blpage{193}
(\byear{1997})
\doiurl{10.1088/0951-7715/10/1/012}
\end{barticle}
\endbibitem

\bibitem[\protect\citeauthoryear{Guardia et~al.}{2010}]{GOS10}
\begin{barticle}
\bauthor{\bsnm{Guardia}, \binits{M.}},
\bauthor{\bsnm{Oliv\'e}, \binits{C.}},
\bauthor{\bsnm{Seara}, \binits{T.M.}}:
\batitle{Exponentially small splitting for the pendulum: a classical problem
  revisited}.
\bjtitle{J. Nonlinear Sci.}
\bvolume{20}(\bissue{5}),
\bfpage{595}--\blpage{685}
(\byear{2010})
\doiurl{10.1007/s00332-010-9068-8}
\end{barticle}
\endbibitem

\bibitem[\protect\citeauthoryear{Neishtadt}{1984}]{Neishtadt84}
\begin{barticle}
\bauthor{\bsnm{Neishtadt}, \binits{A.I.}}:
\batitle{The separation of motions in systems with rapidly rotating phase}.
\bjtitle{J. Appl. Math. Mech.}
\bvolume{48}(\bissue{2}),
\bfpage{133}--\blpage{139}
(\byear{1984})
\doiurl{10.1016/0021-8928(84)90078-9} .
\bcomment{Translated from Prikl. Mat. Mekh. \textbf{48}(2), 197--204 (1984)}
\end{barticle}
\endbibitem

\bibitem[\protect\citeauthoryear{Baldom\'a}{2006}]{Baldoma06}
\begin{barticle}
\bauthor{\bsnm{Baldom\'a}, \binits{I.}}:
\batitle{The inner equation for one and a half degrees of freedom rapidly
  forced {H}amiltonian systems}.
\bjtitle{Nonlinearity}
\bvolume{19}(\bissue{6}),
\bfpage{1415}--\blpage{1445}
(\byear{2006})
\doiurl{10.1088/0951-7715/19/6/011}
\end{barticle}
\endbibitem

\bibitem[\protect\citeauthoryear{Oliv\'e et~al.}{2003}]{OSS03}
\begin{barticle}
\bauthor{\bsnm{Oliv\'e}, \binits{C.}},
\bauthor{\bsnm{Sauzin}, \binits{D.}},
\bauthor{\bsnm{Seara}, \binits{T.M.}}:
\batitle{Resurgence in a {H}amilton--{J}acobi equation}.
\bjtitle{Ann. Inst. Fourier}
\bvolume{53}(\bissue{4}),
\bfpage{1185}--\blpage{1235}
(\byear{2003})
\doiurl{10.5802/aif.1977}
\end{barticle}
\endbibitem

\bibitem[\protect\citeauthoryear{Gelfreich and Sauzin}{2001}]{GS01}
\begin{barticle}
\bauthor{\bsnm{Gelfreich}, \binits{V.}},
\bauthor{\bsnm{Sauzin}, \binits{D.}}:
\batitle{Borel summation and splitting of separatrices for the {H}\'enon map}.
\bjtitle{Ann. Inst. Fourier}
\bvolume{51}(\bissue{2}),
\bfpage{513}--\blpage{567}
(\byear{2001})
\doiurl{10.5802/aif.1831}
\end{barticle}
\endbibitem

\bibitem[\protect\citeauthoryear{Baldom\'a et~al.}{2024}]{BMS24}
\begin{barticle}
\bauthor{\bsnm{Baldom\'a}, \binits{I.}},
\bauthor{\bsnm{M-Seara}, \binits{T.}},
\bauthor{\bsnm{Moreno}, \binits{R.}}:
\batitle{Splitting of separatrices for rapid degenerate perturbations of the
  classical pendulum}.
\bjtitle{SIAM J. Appl. Dyn. Syst.}
\bvolume{23}(\bissue{2}),
\bfpage{1159}--\blpage{1198}
(\byear{2024})
\doiurl{10.1137/23M1550992}
\end{barticle}
\endbibitem

\bibitem[\protect\citeauthoryear{Gomide et~al.}{2025}]{GGSZ25}
\begin{barticle}
\bauthor{\bsnm{Gomide}, \binits{O.M.L.}},
\bauthor{\bsnm{Guardia}, \binits{M.}},
\bauthor{\bsnm{Seara}, \binits{T.M.}},
\bauthor{\bsnm{Zeng}, \binits{C.}}:
\batitle{On small breathers of nonlinear {K}lein--{G}ordon equations via
  exponentially small homoclinic splitting}.
\bjtitle{Invent. Math.}
\bvolume{240}(\bissue{2}),
\bfpage{661}--\blpage{777}
(\byear{2025})
\doiurl{10.1007/s00222-025-01327-y}
{\href{https://arxiv.org/abs/2107.14566}{{arXiv:2107.14566}}}
\end{barticle}
\endbibitem

\bibitem[\protect\citeauthoryear{Baldom\'a et~al.}{2023}]{BCGS23}
\begin{barticle}
\bauthor{\bsnm{Baldom\'a}, \binits{I.}},
\bauthor{\bsnm{Capi{\'n}ski}, \binits{M.J.}},
\bauthor{\bsnm{Guardia}, \binits{M.}},
\bauthor{\bsnm{Seara}, \binits{T.M.}}:
\batitle{Breakdown of heteroclinic connections in the analytic {H}opf-zero
  singularity: rigorous computation of the {S}tokes constant}.
\bjtitle{J. Nonlinear Sci.}
\bvolume{33}(\bissue{2}),
\bfpage{28}
(\byear{2023})
\doiurl{10.1007/s00332-022-09882-x}
{\href{https://arxiv.org/abs/2206.13840}{{arXiv:2206.13840}}}
\end{barticle}
\endbibitem

\bibitem[\protect\citeauthoryear{Baldom\'a et~al.}{2025}]{BCGG25}
\begin{barticle}
\bauthor{\bsnm{Baldom\'a}, \binits{I.}},
\bauthor{\bsnm{Capi{\'n}ski}, \binits{M.J.}},
\bauthor{\bsnm{Giralt}, \binits{M.}},
\bauthor{\bsnm{Guardia}, \binits{M.}}:
\batitle{Breakdown of homoclinic orbits to {$L_3$}: nonvanishing of the
  {S}tokes constant}.
\bjtitle{Discrete Contin. Dyn. Syst.}
\bvolume{45}(\bissue{1}),
\bfpage{56}--\blpage{88}
(\byear{2025})
\doiurl{10.3934/dcds.2024086}
{\href{https://arxiv.org/abs/2312.13138}{{arXiv:2312.13138}}}
\end{barticle}
\endbibitem

\bibitem[\protect\citeauthoryear{Sauzin}{2001}]{Sauzin01}
\begin{barticle}
\bauthor{\bsnm{Sauzin}, \binits{D.}}:
\batitle{A new method for measuring the splitting of invariant manifolds}.
\bjtitle{Ann. Sci. \'Ec. Norm. Sup\'er. (4)}
\bvolume{34}(\bissue{2}),
\bfpage{159}--\blpage{221}
(\byear{2001})
\doiurl{10.1016/S0012-9593(00)01063-6}
\end{barticle}
\endbibitem

\bibitem[\protect\citeauthoryear{Gulkanov}{2026}]{GulkanovCode}
\begin{botherref}
\oauthor{\bsnm{Gulkanov}, \binits{A.}}:
Verification of the certificate for {T}heorem 3 in ``{A}symptotics of
  separatrix splitting in {G}allavotti's problem''.
Zenodo
(2026).
\doiurl{10.5281/zenodo.24615549}
\end{botherref}
\endbibitem

\end{thebibliography}

\end{document}